\documentclass{gtart}

\def\ifplaintex{\expandafter\ifx\csname documentclass\endcsname\relax}

\def\gtp{{\mathsurround=0pt\it $\cal G\mskip-2mu$eometry \&\ 
$\cal T\!\!$opology $\cal P\!$ublications}}  

\def\recd{{\small Received:\qua\receiveddate\ifx\reviseddate\relax
\else\qquad Revised:\qua\reviseddate\fi\par}} 

\def\lognumber#1{\def\thelognumber{#1}}
\def\volumenumber#1{\def\thevolumenumber{#1}}
\def\volumeyear#1{\def\thevolumeyear{#1}}
\def\papernumber#1{\def\thepapernumber{#1}}
\def\pagenumbers#1#2{\def\startpage{#1}\def\finishpage{#2}}
\def\published#1{\def\publishdate{#1}}

\def\received#1{\def\receiveddate{#1}}
\def\revised#1{\def\reviseddate{#1}}
\def\accepted#1{\def\accepteddate{#1}}

\def\asciiaddress#1{\def\theasciiaddress{#1}}

\long\def\asciiabstract#1{\long\def\theasciiabstract{#1}}

\let\\\par\let\thelognumber\relax\let\thevolumenumber\relax
\let\thepapernumber\relax\let\thevolumeyear\relax\let\startpage\relax
\let\finishpage\relax\let\publishdate\relax\let\receiveddate\relax
\let\reviseddate\relax\let\accepteddate\relax\let\theasciititle\relax
\let\theasciiauthors\relax\let\theasciiaddress\relax
\let\theasciiabstract\relax

\let\theasciiemail\relax

\ifplaintex
\font\logobig=cmssbx10 scaled 3836
\font\logomed=cmssbx10 scaled 2557
\else
\font\logobig=cmssbx10 scaled 4200
\font\logomed=cmssbx10 scaled 2800
\fi

\long\def\makeagttitle{   
\count0=\startpage
\agt\hfill      
\hbox to 45truept{\vbox to 0pt{\vglue -13truept{\logomed A\kern -.37em{\logobig 
T}\kern -.38em G}\vss}\hss}
\break
{\small Volume \thevolumenumber\ (\thevolumeyear)
\startpage--\finishpage\nl
Published: \publishdate}

\vglue .25truein

{\parskip=0pt\leftskip 0pt plus
1fil\def\\{\par\smallskip}{\Large\bf\thetitle}\par\medskip} \vglue
0.05truein

{\parskip=0pt\leftskip 0pt plus 1fil\def\\{\par}{\sc\theauthors}
\par\medskip}%
 
\vglue 0.03truein 

{\small\leftskip 25truept\rightskip 25truept{\bf Abstract}\stdspace\theabstract

{\bf AMS Classification}\stdspace\theprimaryclass
\ifx\thesecondaryclass\relax\else; \thesecondaryclass\fi\par
{\bf Keywords}\stdspace \thekeywords\par}\vglue 7truept

}   

\ifplaintex
\font\phead=cmsl9 scaled 950
\font\pnum=cmbx10 scaled 913
\font\pfoot=cmsl9 scaled 950
\headline{\vbox to 0pt{\vskip -4.5mm\line{\small\phead\ifnum
\count0=\startpage ISSN 1472-2739 (on-line) 1472-2747 (printed)
\hfill {\pnum\folio}\else\ifodd\count0\def\\{ }%
\ifx\theshorttitle\relax\thetitle\else\theshorttitle\fi\hfill{\pnum\folio}
\else\def\\{ and }{\pnum\folio}\hfill\ifx\theshortauthors\relax\theauthors
\else\theshortauthors\fi\fi\fi}\vss}}
\footline{\vbox to 0pt{\vglue 0mm\line{\small\pfoot\ifnum\count0=\startpage
\copyright\ \gtp\hfill\else
\agt, Volume \thevolumenumber\ (\thevolumeyear)\hfill\fi}\vss}}
\else
\font=cmsl9 scaled 1050
\font\lnum=cmbx10 
\font=cmsl9 scaled 1050
\makeatletter
\def\@oddhead{{\small\lhead\ifnum\count0=\startpage ISSN 1472-2739 
(on-line) 1472-2747 (printed)\hfill {\lnum\number\count0}\else\ifodd\count0
\def\\{ }\ifx\theshorttitle\relax \thetitle \else\theshorttitle\fi\hfill
{\lnum\number\count0}\else\def\\{ and }{\lnum\number\count0}
\hfill\ifx\theshortauthors\relax 
\theauthors\else\theshortauthors\fi\fi\fi}}\def\@evenhead{\@oddhead}
\def\@oddfoot{\small\lfoot\ifnum\count0=\startpage\copyright\ \gtp\hfill\else
\agt, Volume \thevolumenumber\ (\thevolumeyear)\hfill\fi}
\def\@evenfoot{\@oddfoot}
\makeatother
\fi
\let\maketitlepage\makeagttitle
\let\makeshorttitle\maketitlepage
\let\maketitle\maketitlepage

\newwrite\gtoutfile
\long\gdef\makeheadfile{  
{\def\\{, }\def\s{ }
\immediate\openout\gtoutfile head.xxx
\immediate\write\gtoutfile{To: math@arxiv.org}
\immediate\write\gtoutfile{Subject: put OR rep NNNNN:ppppp}
\immediate\write\gtoutfile{--text follows this line--}
\immediate\write\gtoutfile{Proxy-for: \ifx\theasciiauthors\relax
\theauthors\else\theasciiauthors\fi\s<\ifx\theasciiemail\relax\theemail\else\theasciiemail\fi>}
\immediate\write\gtoutfile{\noexpand\\}
\immediate\write\gtoutfile{Authors: \ifx\theasciiauthors\relax
\theauthors\else\theasciiauthors\fi}
{\def\\{ }\immediate\write\gtoutfile{Title: \ifx\theasciititle\relax
\thetitle\else\theasciititle\fi}}
\immediate\write\gtoutfile{Subj-class: GT or SG, GR etc}
\immediate\write\gtoutfile{MSC-class: \theprimaryclass\ifx\thesecondaryclass\relax\else, \thesecondaryclass\fi}
\immediate\write\gtoutfile{Journal-ref: Algebr. Geom. Topol. \thevolumenumber\s
(\thevolumeyear) \startpage-\finishpage}
\immediate\write\gtoutfile{Comments: Published by Algebraic and
Geometric Topology at}
\immediate\write\gtoutfile{\s\s\s  http://www.maths.warwick.ac.uk/agt/AGTVol\thevolumenumber/agt-\thevolumenumber-\thepapernumber.abs.html}
\immediate\write\gtoutfile{\noexpand\\}
\immediate\write\gtoutfile{}
\ifx\theasciiabstract\relax
\immediate\write\gtoutfile{\theabstract}\else
\immediate\write\gtoutfile{\theasciiabstract}\fi
\immediate\write\gtoutfile{}
\immediate\write\gtoutfile{\noexpand\\}
\immediate\write\gtoutfile{}
\immediate\closeout\gtoutfile}}  

\def\maketitlepage{\makeagttitle\makeheadfile}
\let\makeshorttitle\maketitlepage
\let\maketitle\maketitlepage

\lognumber{12}
\volumenumber{3}
\volumeyear{2003}
\papernumber{12}
\published{21 March 2003}
\pagenumbers{335}{398}
\received{18 April 2002}
\revised{5 February 2003}
\accepted{12 March 2003}

\usepackage{amsmath,amssymb,rlepsf}
\usepackage[all]{xy}
\nocolon

\theoremstyle{plain}
\newtheorem{thh}{Theorem}[section]
\newtheorem{lem}[thh]{Lemma}
\newtheorem{prop}[thh]{Proposition}

\newtheorem{cor}[thh]{Corollary}
\theoremstyle{definition}
\newtheorem{df}[thh]{Definition} 
\newtheorem{rem}{Remark}

\let\tilde\widetilde
\renewcommand{\t}{\tilde}
\renewcommand{\b}{\partial}
\renewcommand{\H}[1]{H_1(#1,{\Z})}
\newcommand{\z}[1]{{\Z}/#1{\Z}}
\newcommand{\Z}{\bf Z}

\newcommand{\N}{\bf N}
\newcommand{\A}{\bf A}
\newcommand{\Q}{\bf Q}
\newcommand{\C}{\bf C}
\newcommand{\F}{\bf F}
\renewcommand{\l}{\langle}
\renewcommand{\r}{\rangle}
\renewcommand{\o}{\overline}
\renewcommand{\phi}{\varphi}

\newcommand{\et}{_{\ast}}
\newcommand{\h}{{\Z})}
\newcommand{\e}{\varepsilon}
\renewcommand{\epsilon}{\varepsilon}

\begin{document}
\title{A criterion for homeomorphism between\\closed Haken manifolds}
\author{Pierre Derbez}
\address{Laboratoire de Topologie, UMR 5584 du CNRS\\
Universit\'e de Bourgogne, 9, avenue Alain Savary -- BP 47870\\
21078 Dijon CEDEX -- France}
\asciiaddress{Laboratoire de Topologie, UMR 5584 du CNRS\\
Universite de Bourgogne, 9, avenue Alain Savary -- BP 47870\\
21078 Dijon CEDEX -- France}
\email{derbez@topolog.u-bourgogne.fr}

\begin{abstract} 
In this paper we consider two connected closed Haken manifolds denoted
by $M^3$ and $N^3$, with the same Gromov simplicial volume. We give a
simple homological criterion to decide when a given map $f\co M^3\to
N^3$ between $M^3$ and $N^3$ can be changed by a homotopy to a
homeomorphism.  We then give a convenient process for constructing
maps between $M^3$ and $N^3$ satisfying the homological hypothesis of
the map $f$.
\end{abstract}

\asciiabstract{In this paper we consider two connected closed Haken 
manifolds denoted by M^3 and N^3, with the same Gromov simplicial
volume.  We give a simple homological criterion to decide when a given
map f: M^3-->N^3 between M^3 and N^3 can be changed by a homotopy to a
homeomorphism.  We then give a convenient process for constructing
maps between M^3 and N^3 satisfying the homological hypothesis of the
map f.}

\maketitle

\let\\\par

\section{Introduction}
\subsection{The main result}
Let $N^3$ be an orientable connected, compact three-manifold without
boundary. We denote by $\|N^3\|$ the Gromov simplicial volume (or
Gromov Invariant) of $N^3$, see Gromov \cite[paragraph 0.2]{Gro} and
Thurston \cite[paragraph 6.1]{Th2} for definitions. Then, our main
result is stated as follows.
\begin{thh}\label{MT}  
Let $M^3$ and $N^3$ be two closed Haken manifolds with the same Gromov simplicial volume. Let $f\co M^3\to N^3$ be a map such that for any finite covering $\t{N}$ of $N^3$ (regular or not) the induced map $\t{f}\co \t{M}\to\t{N}$ is a homology equivalence (with coefficients ${\Z}$). Then $f$ is homotopic to a homeomorphism.
\end{thh}
Note that the homological hypothesis on the map $f$ required by the Theorem \ref{MT} is usually not easy to check. The following result,  \cite[Proposition 0.2 and Lemma 0.6]{Pe-S}, gives a convenient process which allows us to construct such a map between $M^3$ and $N^3$.
\begin{prop}
 Let $M^3$, $N^3$ be two closed Haken manifolds and assume that there is a cobordism $W^4$ between $M^3$ and $N^3$ such that:\\
{\rm(i)}\qua the map $\pi_1(N^3)\to\pi_1(W^4)$ is an epimorphism,\\
{\rm(ii)}\qua $W^4$ is obtained from $N^3$ adding  handles of index $\leq 2$,\\
{\rm(iii)}\qua the inclusions $M^3\hookrightarrow W^4$ and $N^3\hookrightarrow W^4$ are ${\Z}$-homological equivalences.\\
Then there exists a map $f\co M^3\to N^3$ satisfying the homological hypothesis of Theorem \ref{MT} and thus if $\|M\|=\|N\|$ then $f$ is homotopic to a homeomorphism.
\end{prop}
\subsection{The motivation}
The aim of Theorem \ref{MT} is to extend a main result of B. Perron and P. Shalen which gives a homological criterion for deciding when a given map between two closed, irreducible, graph manifolds, with infinite fundamental group, can be homotoped to a homeomorphism (see  \cite[Proposition 0.1]{Pe-S}). Thus, in this paper we want to find a larger class of three-manifolds for which Proposition 0.1 of B. Perron and P. Shalen holds. Obviously their result does not hold for any closed three-manifold. Consider for example a (closed) ${\Z}$-homology sphere $M^3$ such that $\|M\|=0$ and $M^3\not\simeq{\bf S}^3$. Then it is easy to construct a map $f\co M^3\to {\bf S}^3$ which satisfies the hypothesis of Theorem \ref{MT}. In order to generalize the result of B. Perron and P. Shalen, a ``good" class of closed three-manifolds seems to be the Haken manifolds. This class allows us to avoid the above type of obvious conter-example and strictly contains the class of irreducible graph manifolds with infinite fundamental group considered by B. Perron and P. Shalen. Indeed, it follows from Thurston \cite{Th2} and \cite[paragraph IV.11]{Ja} that irreducible graph manifolds with infinite fundamental group correspond exactly to Haken manifolds with zero Gromov Invariant. Thus when the given manifolds $M^3$ and $N^3$ have their Gromov Invariant equal to zero  (i.e.\ if $\|M^3\|=\|N^3\|=0$) then Theorem \ref{MT} is equivalent to \cite[Proposition 0.1]{Pe-S}. Therefore, the result of \cite{Pe-S} allows us, from now on, to assume  that the given manifolds satisfy $\|M^3\|=\|N^3\|\not=0$.

Finally note that the hypothesis on the Gromov Invariant of the given manifolds is necessary in Theorem \ref{MT}. Indeed in \cite{Bo-W}, M. Boileau and S. Wang construct two closed Haken manifolds $M^3$ and $N^3$ satisfying $\|M\|>\|N\|$ and a map $f\co M\to N$ satisfying the homological hypothesis of Theorem \ref{MT}. 

\subsection{Preliminaries and notations}
We first state the following terminology which will be convenient. Let ${\cal T}$ be a 2-manifold whose components are all tori and let $m$ be a positive integer. A covering space $\t{\cal T}$ of ${\cal T}$ will be termed $m\times m- characteristic$  if each component of $\t{\cal T}$ is equivalent to the covering space of some component $T$ of ${\cal T}$ associated to the characteristic subgroup $H_m$ of index $m^2$ in $\pi_1(T)$ (if we identify $\pi_1(T)$ with ${\Z}\times{\Z}$, we have $H_m=m{\Z}\times m{\Z}$).
  
 Recall that for a closed Haken manifold $M^3$, the torus decomposition of Jaco-Shalen and Johannson (\cite{Ja-S} and \cite{Jo}) together with the uniformization Theorem of Thurston (\cite{Th1}) say that there is a collection of incompressible tori $W_M\subset M$, unique up to ambiant isotopy, which cuts $M$ into Seifert fibered manifolds and hyperbolic manifolds of finite volume. Denote the regular neighborhood of $W_M$ by $W_M\times [-1,1]$ with $W_M\times\{0\}=W_M$. We write $M\setminus W_M\times(-1,1)=H_M\cup S_M$, where $H_M$ is the union of the finite volume hyperbolic manifold components and $S_M$ is the union of the Seifert fibered  manifold components. Note that since we assume that $\|M\|\not=0$ we always have $H_M\not=\emptyset$.
 
The hypothesis on the Gromov simplicial volume of the given manifolds allows us to apply the following rigidity Theorem of Soma:
\begin{thh}\label{soma}{\rm\cite[Theorem 1]{So}}\qua Let $f\co M\to N$ be a proper, continuous map of strictly positive degree between two Haken manifolds with (possibly empty) toral boundary. Then $f$ is properly homotopic to a map $g$ such that $g(H_M)\subset H_N$ and $g|H_M\co H_M\to H_N$ is a $deg(f)$-fold covering if and and only if $\|M\|=deg(f)\|N\|$.
\end{thh}
 In our case this result implies that the map $f\co M\to N$ is homotopic to a map $g$ which induces a homeomorphism between $H_M$ and $H_N$. But this result does not say anything about the behavior of $f$ on the Seifert components $S_M$ of $M$. Even if we knew that $f(S_M)\subset S_N$ we can not have a reduction to the Perron-Shalen case ``with boundaries" (which is not anyway treated in their article). This comes from the fact that one does not know how to extend a given finite covering of $S_N$ to the whole manifold $N$, see \cite[Lemma 4.1]{He2}. More precisely, in \cite{He2}, J. Hempel shows that if $S$ is a 3-manifold with non-empty boundary which admits either a Seifert fibration or a complete hyperbolic structure of finite volume then for all but finitely many  primes $q$ there is a finite covering $p\co\t{S}\to S$ such that for each component $T$ of $\b 
 S$ and for each component $\t{T}$  of $p^{-1}(T)$ the induced map $p|\t{T}\co\t{T}\to T$ is the $q\times q$-characteristic covering of $T$. In particular, we can show that Hempel's Lemma is true for any prime $q$ in the case of  Seifert fibered spaces without exceptional fiber and with orientable base whose boundary contains at least two boundary components (i.e.\ $S\simeq F\times{\bf S}^1$ where $F$ is an orientable compact surface with at least two boundary components). This fact is crucially used in \cite{Pe-S} (see  proofs of Propositions 0.3 and 0.4) to construct their finite coverings. But in the hyperbolic manifolds case we must  exclude a finite collection  of primes, thus we cannot extend the coverings of \cite{Pe-S} in our case. So we have to develop some other techniques  to avoid these main difficulties.
 \subsection{Main steps in the proof of Theorem \ref{MT} and statement of the intermediate results}    
It follows from Waldhausen, see \cite[Corollary 6.5]{Wa}, that to prove Theorem \ref{MT} it is sufficient to show that the map $f$ induces an isomorphism $f_{\ast}\co\pi_1(M)\to\pi_1(N)$. Note that since $f$ is a ${\Z}$-homology equivalence then it is a degree one map  so it is sufficient to see that  $f_{\ast}\co\pi_1(M)\to\pi_1(N)$ is injective. On the other hand it follows from the hypothesis of Theorem \ref{MT} that to prove Theorem \ref{MT} it is  sufficient to find a finite covering $\t{N}$ of $N$ such that the induced map $\t{f}\co\t{M}\to\t{N}$ is homotopic to a homeomorphism (i.e.\ is $\pi_1$-injective). Hence, we can replace $M$, $N$ and $f$ by $\t{M}$, $\t{N}$ and $\t{f}$ (for an appropriate choice of the finite sheeted covering of $N$).

{\bf First step: Simplification of $N^3$}\qua The first step consists in finding some finite covering $\t{N}$ of $N$ which is more ``convenient" than $N$. More precisely, the first step  is to show the following result whose proof will occupy Section \ref{Srev}.
\begin{prop}\label{rev} 
Let $N^3$ be a non geometric closed Haken manifold.   Then there is a finite covering $\t{N}$ of $N$ satisfying the following property: $\t{N}^3$ has large first Betti number ($\beta_1(\t{N})\geq 3)$, each component of $\t{N}\setminus W_{\t{N}}$ contains at least two components in its boundary and each Seifert fibered space of $\t{N}$ is homeomorphic to a product of type $F\times S^1$ where $F$ is an orientable surface of genus $\geq 3$.
\end{prop}
{\rem In view of the  above paragraph we assume  now that $N^3$ always satisfies the conclusion of Proposition \ref{rev}.}
\medskip

{\bf Second step: The obstruction}\qua This step will show that to prove Theorem \ref{MT} it is sufficient to see that the canonical tori of $M$ do not degenerate (i.e.\ the map $f|W_M\co W_M\to N$ is $\pi_1$-injective). More precisely we state here the following result which will be proved in Section \ref{Storesinjectifs}.
\begin{thh}\label{toresinjectifs} Let $f\co M^3\to N^3$ be a map between two closed Haken manifolds with the same Gromov Invariant and such that for any finite covering $\t{N}$ of $N$ the induced map $\t{f}\co\t{M}\to\t{N}$ is a ${\Z}$-homology equivalence. Then $f$ is homotopic to a homeomorphism if and only if  the   induced map $f|W_M\co W_M\to N^3$ is $\pi_1$-injective.
\end{thh}
{\bf Third step: A factorization theorem}\qua  It follows from Theorem \ref{toresinjectifs} that to show our homeomorphism criterion it is sufficient to see that the canonical tori do not degenerate under the map $f$. So in the following we will suppose the contrary. The purpose of this step is to understand the behavior (up to homotopy) of the map $f$ in the case of degenerate tori. To do this we recall the definition of {\it degenerate maps} of Jaco-Shalen. 
\begin{df} Let $S$ be a Seifert fibered space and let $N$ be a closed Haken manifold. A map $f\co S\to N$ is said to be degenerate if either:

(1)\qua the group Im$(f_{\ast}\co\pi_1(S)\to\pi_1(N))=\{1\}$, or

(2)\qua the group Im$(f_{\ast}\co\pi_1(S)\to\pi_1(N))$ is cyclic, or

(3)\qua the map $f|\gamma$ is homotopic in $N$ to a constant map for some fiber $\gamma$ of $S$.
\end{df}
So we first state  the following result which explains how certain submanifolds of $M^3$ can degenerate.
\begin{thh}\label{toresdegeneres} Let $f\co M\to N$ be a map between two closed Haken manifolds satisfying hypothesis of Theorem \ref{MT} and suppose that $N$ satisfies the conclusion of Proposition \ref{rev}. Let $T$ be a canonical torus  in $M$ which degenerates under the map $f$.  Then $T$ separates $M$ in two submanifolds $A$, $B$, one and only one (say $A$) satisfies the followings properties:

{\rm(i)}\qua $\H{A}={\Z}$ and each Seifert component of $A\setminus W_M$ admits a Seifert fibration whose orbit space is a surface of genus 0,

{\rm(ii)}\qua each Seifert component of $A\setminus W_M$ degenerates under the map $f$, $A$ is a graph manifold and the group $f_{\ast}(\pi_1(A))$ is either trivial or infinite cyclic.
\end{thh} 
With this result we may write the following definitions.
\begin{df}\label{cdm} Let $M^3$ and $N^3$ be two closed, connected, Haken manifolds and let $f\co M^3 \to N^3$ be a map satisfying  hypothesis of Theorem \ref{MT}. We say that a codimesion 0 submanifold $A$ of $M$ is a  maximal end of $M$ if $A$ satisfies the following three properties:  

(i)\qua $\partial A $ is a single incompressible torus,  $H_1(A,{\bf Z}) = {\bf Z}$ and  $f_{\ast}(\pi_1(A)) = {\bf Z}$,

(ii)\qua if $p\co\tilde{M}\to M$ is any finite covering induced by $f$ from some finite covering $\tilde{N}$ of $N$ then each component of $p^{-1}(A)$ satisfies (i), 

(iii)\qua if $C$ is a submanifold of $M$ which contains $A$ and  satisfying  (i) and (ii) then $A = C$.
\end{df}  
To describe precisely the behavior of the map $f$ (up to homotopy) we still need  the following definition:
\begin{df}\label{eff} Let $M$ be a closed, connected, compact 3-manifold and let $A$ be a compact, connected codimension 0 submanifold of $M$ whose boundary is a torus in $M$. We say that $M$  collapses along $A$ if there exists a homeomorphism $\varphi\co\partial(D^2\times S^1)\to\partial A=\partial(\overline{M\setminus A})$ and a map $\pi\co M\to (\overline{M\setminus A})\cup_{\varphi}D^2\times S^1$ such that $\pi|\overline{M\setminus A} = id$ and $\pi(A) = D^2\times S^1$.
\end{df}
So using Theorem \ref{toresinjectifs} and Theorem \ref{toresdegeneres} we obtain the following {\it factorization} Theorem  which will be used to get a good decription of the behavior of the map $f$.
\begin{thh}\label{facto} Let $M^3$ and $N^3$ be two closed, connected, Haken manifolds  satisfying hypothesis of Theorem \ref{MT} and assume that $N$ satisfies the conclusion of Proposition \ref{rev}. Then there exists a finite family $\{A_1,...,A_{n_M}\}$ (eventually empty) of disjoint maximal ends of $M$, a  Haken manifold $M_1$ obtained from $M$ by collapsing $M$ along the family $\{A_1,...,A_{n_M}\}$ and a homeomorphism $f_1\co M_1 \to N$ such that  $f$ is homotopic to the map $f_1\circ\pi$, where $\pi$ denotes the collapsing map $\pi\co M \to M_1$.
\end{thh}
Note that Theorems \ref{toresdegeneres} and \ref{facto} remain true if we simply assume that  the given manifolds $M^3$, $N^3$ and the  map $f\co M^3\to N^3$ satisfies hypothesis of Theorem \ref{MT}. But it is more convenient for our purpose to suppose that $N^3$ satisfies the conclusion of Proposition \ref{rev}.

{\bf Fourth step}\qua The purpose of this step is to show that the hypothesis which says that certain canonical tori degenerate is finally absurd. To do this, we will show that if $A$ is a maximal end of $M$ then we can construct a finite covering $p\co\t{M}\to M$ induced by $f$ from some finite covering of $N$, such that the connected components of $p^{-1}(A)$ are not maximal ends, which contradicts Definition \ref{cdm}. But to construct such a covering, it is first necessary to have good informations about the behavior of the induced map $f|A\co A\to N$ up to homotopy. To do this we state the following result whose proof depends crucially on Theorem \ref{facto} (see Section \ref{consequencedefacto}):
\begin{prop}\label{trivialise} Let $f\co M\to N$ be a map between two closed Haken manifolds with the same Gromov Invariant satisfying the hypothesis of Theorem \ref{MT} and assume that $N$ satisfies the conclusion of Proposition \ref{rev}. If $A$ denotes a maximal end of $M$ then there exists a Seifert piece $S$ in $A$, whose orbit space is a disk such that $f_{\ast}(\pi_1(S))\not=\{1\}$,  a Seifert piece $B=F\times{\bf S}^1$ in $N$ such that $f(S)\subset f(A)\subset B$ and $f_{\ast}(\pi_1(S))\subset\l t\r$, where $t$ denotes the homotopy class of the fiber in $B$.
\end{prop} 
The aim of this  result is to replace  the {\it Mapping Theorem} (see \cite[Chapter III]{Ja-S}) which says that if a map between a Seifert fibered space and a Haken manifold satisfies certains {\it good} properties of non-degeneration then it can be changed by a homotopy in such a way that its whole image is contained in a Seifert fibered space. But when such a map degenerates (which is the case for $f|A$) its behavior can be very complicated a priori. 

The above result shows that the map $f|A$ is homotopically very simple.
 We next construct a finite covering $p\co\t{M}\to M$, induced by $f$ from some finite covering of $N$ such that the component of $p^{-1}(S)$ admits a Seifert fibration whose orbit space is a surface of genus $>0$. Then using \cite[Lemma 3.2]{Pe-S} we show that the components of $p^{-1}(A)$ are not  maximal ends which gives a contradiction. The construction of our finite covering depends crucially on the following result  
  which completes the proof of the fourth step  and  whose proof is based on the Thurston {\it Deformation Theory} of complete finite volume hyperbolic structures  and will be proved in Section 6.3.
\begin{prop}\label{def} Let $N^3$ be a closed Haken manifold with non-trivial Gromov simplicial volume. Then there exists a finite covering $\t{N}$ of $N$ satisfying the following property: for every integer $n_0>0$ there exists an integer $\alpha>0$ and a finite covering $p\co\hat{N}\to\t{N}$ such that for each Seifert piece $\t{S}$ of $\t{N}\setminus W_{\t{N}}$ and for each component $\hat{S}$ of $p^{-1}(\t{S})$ the map $p|\hat{S}\co\hat{S}\to\t{S}$ is fiber preserving and induces the $\alpha n_0$-index covering on the fibers of $\t{S}$.
\end{prop}
Note that this result plays a Key Role in the proof of Theorem \ref{MT}. Indeed, this Proposition \ref{def} allows us to avoid the main difficulty stated in paragraph 1.3.

\section{\label{Srev} Preliminary results on Haken manifolds}
In this section we state some general results on Haken manifolds and their finite coverings which will be useful in the following of this article. On the other hand  we will always suppose in the following that the given manifold $N$ has non trivial Gromov simplicial volume which implies in particular that $N$ has no finite cover which is fibered over the circle by tori.

\subsection{Outline of proof of Proposition \ref{rev}}
 In this section we outline the proof of Proposition \ref{rev} which  extends in the Haken manifolds case the result of \cite{L-W} which concerns graph manifolds. For a complete proof of this result see \cite[Proposition 1.2.1]{De}.

First note that since $N$ is a non geometric Haken manifold then $N$ is not a Seifert fibered space (in particular $N$ has a non empty torus decomposition) and has no finite cover that fibers as a torus bundle over the circle. By \cite[Theorem 2.6]{Lu} we may assume, after passing possibly to a finite cover, that each component of $N\setminus W_N$ either has hyperbolic interior or is Seifert fibered over an orientable surface whose base 2-orbifold has strictly negative Euler characteristic.

By applying either \cite[Theorem 2.4]{Lu} or \cite[Theorem 3.2]{Lu} to each piece $Q$ of $N\setminus W_N$ (according to whether the piece is Seifert fibered or hyperbolic, resp.) there is a prime $q$, such that for every $Q$ in $N\setminus W_N$ there is a finite, connected, regular cover $p_Q\co\t{Q}\to Q$ where, if $T$ is a component of $\b Q$, then $(p_Q)^{-1}(T)$ consists of more than one component; furthermore, if $\t{T}$ is a component of $(p_Q)^{-1}(T)$, then $p_Q|\t{Q}\co\t{Q}\to Q$ is the $q\times q$-characteristic covering. This allows us to glue the covers of the pieces of $N\setminus W_N$ together to get a covering $\t{N}$ of $N$ in which each piece of $\t{N}\setminus W_{\t{N}}$ has at least two boundary components. By repeating this process, we may assume, after passing to a finite cover, that each component of $N\setminus W_N$ has at least three boundary components. 

Let $S$ be a Seifert piece of $N$ and let $F$ be the orbit space of $S$. Let $T_1,...,T_p$ ($p\geq 3$) be the components of $\b S$, $D_1,...,D_p$ those of $\b F$ and set $d_i=[D_i]\in\pi_1(F)$ (for a choice of base point). With these notations we have: $\pi_1(T_i)=\l d_i,h\r$ where $h$ denotes the regular fiber in $S$. Since $S$ has at least three boundary components then using the presentation of $\pi_1(S)$ one can show that for all but finitely many primes $q$ there exists an epimorphism $\varphi\co\pi_1(S)\to\z{q}\times\z{q}$ such that: 

(i)\qua $\varphi(d_j)\not\in\l\varphi(h)\r$ for $j=1,...,p$,

(ii)\qua $\ker(\varphi|\pi_1(T_j))$ is the $q\times q$-characteristic subgroup of $\pi_1(T_j)$ for $j=1,...,p$.

Let $\pi\co\t{S}\to S$ be the finite covering of $S$ corresponding to $\varphi$ and let $\pi_F\co\t{F}\to F$ be the finite (branched) covering induced by $\pi$ between the orbit spaces of $S$ and $\t{S}$. Then using (i) and (ii) combined with the Riemann-Hurwitz formula \cite[pp.\ 133]{Pr-S} one can show that $\t{g}> g$ where $g$ (resp.\ $\t{g}$) denotes the genus of $F$ (resp.\ of $\t{F}$). Thus, by applying this result combined with \cite[Theorem 3.2]{Lu} to each piece $Q$ of $N\setminus W_N$ (according to whether the piece is Seifert fibered or hyperbolic, resp.) there is a prime $q$, such that for every $Q$ in $N\setminus W_N$ there is a finite, connected, regular cover $p_Q\co\t{Q}\to Q$ where, if $T$ is a component of $\b Q$ and if $\t{T}$ is a component of $(p_Q)^{-1}(T)$, then $p_Q|\t{T}\co\t{T}\to T$ is the $q\times q$-characteristic covering. This allows us to glue the covers of the pieces of $N\setminus W_N$ together to get a covering $\t{N}$ of $N$. Furthermore, if $Q$ is a Seifert piece of $N$ whose orbit space is a surface of genus $g$ then $\t{Q}$ is a Seifert piece of $\t{N}$ whose orbit space is a surface of genus $\t{g}>g$. 

 It remains to see that $N$ is finitely covered by a Haken manifold in which each Seifert piece is a trivial circle bundle. Since the Euler characteristic of the orbit space of the Seifert pieces of $N$ is non-positive then by Selberg Lemma each orbit space is finitely covered by an orientable surface. This covering induces a finite covering (trivial when restricted on the boundary) of the Seifert piece by a circle bundle over an orientable surface, which is trivial because the boundary is not empty. Now we can (trivially) glue these coverings together to get the desired covering of $N$.       
 
\subsection{A technical result for Haken manifolds}
\begin{prop}\label{technique} Let $N^3$ be a closed Haken manifold satisfying the conclusion of Proposition \ref{rev} and let $B$ be a Seifert piece of $N$. Let $g$ and $h$ be elements of $\pi_1({B})\subset\pi_1({N})$ such that either $[g,h]\not=1$ or the group $\langle g,h\rangle$ is the free abelian group of rank two. Then there exists a finite group $H$ and a homomorphism $\phi\co\pi_1({N})\to H$ such that $\phi(g)\not\in\langle\phi(h)\rangle$.
\end{prop}
The proof of this result depends on the following lemma which allows to extend to the whole manifolds $N$  certain ``good" coverings of a given Seifert piece in $N$. 
\begin{lem}\label{etendre}  Let $N$ be a closed Haken manifold such that each Seifert piece is a product and has more than one boundary component and let $B_0$ be a Seifert piece in $N$. Then there exists a prime $q_0$ satisfying the following property: for every finite covering $\t{B}_0$ of $B_0$ which induces the $q^r\times q^r$-characteristic covering on the bounbary components of $B_0$ with $q\geq q_0$ prime and $r\in{\Z}$, there exists a finite covering $\pi\co\t{N}\to N$ such that

{\rm(i)}\qua the covering $\t{N}$ induces the $q^r\times q^r$-characteristic covering on each of the canonical tori of $N$,

{\rm(ii)}\qua each component of the covering of $B_0$ induces by $\t{N}$ is equivalent to $\t{B}_0$.
\end{lem}
The proof of this result depends of the following Lemma which is a slight generalization of Hempel's Lemma, \cite[Lemma 4.2]{He2} and whose proof may be found in \cite[Lemma 1.2.3]{De}.
\begin{lem}\label{Luecke}  Let $G$ be a finitely generated group and let $\tau\co G\to SL(2,{\C})$ be a discret and faithful representation of $G$. Let $\lambda_1,...,\lambda_n$ be elements of $G$ such that $\lambda_i\not=1_G$ and $tr(\tau(\lambda_i))=\pm 2$. Then for all but finitely many  primes $q$ and for all integers $r$ there exists a finite ring ${\A}_{q^r}$ over $\z{q^r}$ and a representation $\tau_q\co G\to SL(2,{\A}_{q^r})$ such that for each element $g\in G$ satisfying $tr(\tau(g))=\pm 2$ the element $\tau_q(g)$ is of order $q^{r_g}$, with $r_g\leq r$ in $SL(2,{\A}_{q^r})$ and the elements $\tau_q(\lambda_i)$ are of order $q^r$ in  $SL(2,{\A}_{q^r})$.
\end{lem}
\begin{proof}[Outline of proof of Lemma \ref{etendre}] We show that if $B$ denotes a component of $N\setminus W_N$ such that $B\not= B_0$ then for each $r\in{\Z}$ and for all but finitely many  primes $q$ there exists a connected regular finite covering $\t{B}$ of $B$ which induces the $q^r\times q^r$-characteristic covering on each of the boundary component of $B$. Next we use similar arguments as in \cite{Lu} using Lemma \ref{Luecke} (see \cite[Lemma 1.2.2]{De}).
\end{proof}
\begin{proof}[Proof of Proposition \ref{technique}] Recall that $B$ can be identified to a product $F\times S^1$, where $F$ is an orientable surface of genus $\geq 1$ with at least two boundary components. Let $D_1,...,D_n$ denote the components of $\b F$ and set $d_i=[D_i]$, for $i=1,...,n$ (for a choice of base point).   

\medskip
{\bf Case 1}\qua If $[g,h]\not=1$, then since $\pi_1(N)$ is a residually finite group (see \cite[Theorem 1.1]{He1}) there is a finite group $H$ and an epimorphism $\varphi\co\pi_1(N)\to H$ such that $\varphi([g,h])\not=1$ and so $\varphi(g)\not\in\langle\varphi(h)\rangle$.

\medskip
{\bf Case 2}\qua If  $[g,h]=1$ then we may write $g=(u^{\beta},t^{\alpha})$ and  $h=(u^{\beta'},t^{\alpha'})$ with $u\in\pi_1(F)$ and where $t$ is a generator of $\pi_1({\bf S^1})={\Z}$. Since $\langle g,h\rangle$ is the free abelian group of rank 2 then $\beta\alpha'-\beta'\alpha=\gamma\not=0$ and $u\not=1$. We first show the following assertion:

{\sl For all but finitely many  primes $p$ there exists an integer $r_0$ such that for each integer $r\geq r_0$ there is a finite group $K$ and a homomorphism $\psi\co\pi_1(B)\to K$ inducing the $p^r\times p^r$-characteristic homomorphism on $\pi_1(\b B)$ and such that $\psi(g)\not\in\langle\psi(h)\rangle$.}

To prove this assertion we consider two cases.

{\bf Case 2.1}\qua Assume first that $\alpha'\not=0$. Choose a prime $p$ such that $(p,\alpha')=1$ and $(p,\gamma)=1$. Then using Bezout's Lemma we may find an integer $n_0$ such that $\beta-n_0\beta'\not\in p{\Z}$. Then using the Key Lemma  on surfaces of B. Perron and P. Shalen, \cite[Key Lemma 6.2]{Pe-S}, by taking $g=u$ we get a homomorphism $\rho\co\pi_1(F)\to H_F$, where $H_F$ is a  $p$-group and satisfying $\rho(u)\not=1$ and $\rho(d_i)$ has order $p^r$ in $H_F$. Let $\lambda\co{\Z}\to\z{p^r}$ denote the canonical epimorphism and consider the following homomorphism:
$$\psi=\rho\times\lambda\co\pi_1(F)\times{\Z}\to H_F\times\z{p^r}$$   
It follows now easily from the above construction that $\psi(g)\not\in\langle\psi(h)\rangle$ and $\ker(\psi|\l d_i,t\r)=\l d_i^{p^r},h^{p^r}\r$.

\medskip
{\bf Case 2.2\qua} We now suppose that $\alpha'=0$. Thus we have  $g=(u^{\beta},t^{\alpha})$ and  $h=(u^{\beta'},1)$ with $\beta'\alpha\not=0$. Recall that $\pi_1(F)=\l d_1\r\ast...\ast\l d_{n-1}\r\ast L_q$ with $d_i=[D_i]$, where $D_1,...,D_n$ denote the components of $\b F$ and where $L_q$ is a free group. Let $\rho_2\co\pi_1(F)\to{\Z}$ be an epimorphism such that $\rho_2(d_1)=...\rho_2(d_{n-1})=1$ and $\rho_2(L_q)=0$. This implies that $\rho_2(d_n)=-(n-1)$. Choose a prime $p$ satisfying $(p,\alpha)=1$, $(p,n-1)=1$ and let $\varepsilon\co{\Z}\to\z{p^r}$ be the canonical epimorphism. So consider the following homomorphism.
$$\psi=(\varepsilon\circ\rho_2)\times\varepsilon\co\pi_1(B)=\pi_1(F)\times{\Z}\to\z{p^r}\times\z{p^r}$$ We now check easily that $\psi(g)\not\in\langle\psi(h)\rangle$ and $\ker(\psi|\l d_i,t\r)=\l d_i^{p^r},h^{p^r}\r$ which completes the proof of the above assertion.

Let $\hat{\pi}\co\hat{B}\to B$ be the covering corresponding to the above homomorphism $\psi$. Since this covering induces the $p^r\times p^r$-covering on each component of $\b B$ then using Lemma  \ref{etendre} there is a finite covering $\pi:\hat{N}\to N$ of $N$ such that each component of $\pi^{-1}(B)$ is equivalent to $\hat{\pi}$. We identify $\pi_1(\hat{N})$ as a subgroup of finite index of $\pi_1(N)$.  Let $\Lambda$ be a subgroup of $\pi_1(\hat{N})$ such that $\Lambda$ is a finite index regular subgroup of $\pi_1(N)$. Then the canonical epimorphism $\varphi:\pi_1(N)\to\pi_1(N)/\Lambda$ satisfies the conclusion of Proposition \ref{technique}.
\end{proof}    
\section{\label{Storesinjectifs} Proof of Theorem \ref{toresinjectifs}}
In this section we always assume that the manifold $N^3$ has non-trivial Gromov Invariant and satisfies the conclusion of Proposition \ref{rev}. 
\subsection{Main ideas of the proof of Theorem \ref{toresinjectifs}}
It follows from the Rigidity Theorem of Soma (see Theorem \ref{soma}) that $f$ is properly homotopic to a map, still denoted by $f$, such that $f|(H_M,\b H_M):(H_M,\b H_M)\to(H_N,\b H_N)$ is a homeomorphism. Then we will prove (see Lemma \ref{Rong}) that we may arrange $f$ by a homotopy fixing $f|H_M$ such that $f(S_M)\subset S_N$: this is the Mapping Theorem of W. Jaco and P. Shalen with some care. So our main purpose here is to find a finite covering $\t{N}$ of $N$ such that for each component $\t{B}$ of $S_{\t{N}}$ there exists exactly one Seifert piece $\t{A}$ of $S_{\t{M}}$ such that $f(\t{A},\b\t{A})\subset(\t{B},\b\t{B})$. We next prove that the induced map  $f|(\t{A},\b\t{A})$ is homotopic to a homeomorphism. To do this the key step consists, for technical reasons, in finding a covering  $\tilde{M}_0$ of $M$, induced by $f$, such that for each Seifert piece $A_i$ of $S_{\tilde{M}_0}$ the induced covering $\t{A}_i$ over $A_i$ is a Seifert fibered space whose orbit space is  a surface of genus $\geq 3$.  This step depends on Proposition \ref{technique}. Indeed the construction of $\t{M}_0$ will be splitted in two steps:
   
\medskip
{\bf First step}\qua The first step is to prove that there exists a finite covering $\tilde{M}_0$ of $M$ induced by $f$ from some finite covering $\tilde{N}_0$ of $N$ in which each Seifert piece is either based on a surface of genus $\geq 3$ ({\it type I})  or based on an annulus ({\it type II}) (see Lemma \ref{1etape}). More precisely the result of Lemma \ref{1etape} is the ``best" that we may obtain using Proposition \ref{technique}.

\medskip
{\bf Second step}\qua The main purpose of this step is to prove, using specific arguments,  that  $\tilde{M}_0$ contains no Seifert piece of type II.  More precisely, if $A_i$ denotes a Seifert piece of type II in $\tilde{M}_0$ then using  \cite[Characteristic Pair Theorem]{Ja-S} we know that there is a Seifert piece $B_j$ in $N$ such that $f(A_i)\subset$ int$(B_j)$ (up to homotopy). Then we construct a vertical torus $U$ in $B_j$ such that if $T$ is a component of $\b A_i$ then $f$ may be changed by a homotopy fixing $\tilde{M}_0\setminus A_i$ so that $f|T\co T\to U$ is a homeomorphism. We next use the structure of $\pi_1(B_j)$ to show that this implies that $A_i$ has no exceptional fiber (i.e.\ $A_i=S^1\times S^1\times I$) which contradicts the minimality of the Torus Decomposition of $\tilde{M}_0$. 
\medskip

Finally we show  that the results obtained in the above steps allows us to use arguments similar to those of \cite[paragraphs 4.3.15 and 4.3.16]{Pe-S} to complete the proof (see paragraph 3.5).
\subsection{Proof of the first step}
This section is devoted to the outline of proof of the following result (for a complete proof see \cite[Lemma 3.2.1]{De}).
\begin{lem}\label{1etape} There exists a finite covering $\t{M}_0$ of $M$ induced by  $f$ from some finite covering $\t{N}_0$ of $N$ in which each Seifert piece $\t{A}$ is either based on a surface of genus   $\geq 3$ (Type I) or satisfies the following properties (Type II):

{\rm(i)}\qua the orbit space of  $\tilde{A}$ is an annulus,

{\rm(ii)}\qua the group $f_{\ast}(\pi_1(\tilde{A}))$ is isomorphic to $\Z$$\oplus$$\Z$,

 {\rm(iii)}\qua for each finite covering  $\pi$: $\hat{M}\to\tilde{M}_0$ induced by $f$ from some finite covering of $\t{N}_0$ then each component of  $\pi^{-1}(\tilde{A})$ satisfies  points  (i) and (ii).
 \end{lem}
 The proof of this result depends on the following lemma.
\begin{lem}\label{augmenter} Let  S be a Seifert piece in $M$ whose orbit space is a surface of genus 0. Suppose that  $S$ contains at least three non-degenerate boundary components.  Then there exists a finite covering $\tilde{S}$ of S satisfying the two followings properties:

{\rm(i)}\qua  $\tilde{S}$ admits a Seifert fibration whose orbit space  is a  surface of genus $\geq$ 1,

{\rm(ii)}\qua $\tilde{S}$ is equivalent to a component of the covering induced from some finite covering of $N$ by $f$.
\end{lem}
\begin{proof} Denote by $F$ the orbit space of the Seifert piece $S$. Let $T_1,...,T_j$, $j\geq 3$, be the non-degenerate tori in $\b S$ and $\pi_1(T_l)=\l d_l,h\r$, $1\leq l\leq j$, the corresponding fundamental groups. Since Rk($\l f_\ast(d_l), f_\ast(h)\r$) = 2 for $l\leq j$, it follows from Proposition 1.4, that there exists a finite group $H$ and a homomorphism $\varphi\co\pi_1(N)\to H$ such that $\varphi\circ f_{\ast}(d_l)\not\in\l\varphi f_{\ast}(h)\r$ for $1\leq l\leq j$.

Let $K$ be the group $\varphi f_{\ast}(\pi_1(S))$ and denote by $\pi\co\t{S}\to S$ the finite covering corresponding to $\varphi\circ f_{\ast}\co\pi_1(S)\to K$. Then $\t{S}$ inherits a Seifert fibration with some base $\t{F}$. We denote by $\tau$ the order of $K$, by $t$ the order of $\varphi f_{\ast}(h)$ in $K$ and by $\beta_i$ the order of $\varphi f_{\ast}(c_i)$ where $c_1,...,c_r$ denote the exceptional fibers of $S$ with index $\mu_1,...,\mu_r$. The map $\pi$ induces a covering $\pi_F\co\t{F}\to F$ on the orbit spaces of $S$ and $\t{S}$ with degree $\sigma=\tau/t$, ramified at the points $\o{c}_i\in F$ corresponding to the exceptional fiber $c_i$ of $S$. Let $\delta_l$, $l=1,...,p$, denote the boundary components of $F$ corresponding to $d_l$ and let $\t{\delta}^1_l,...,\t{\delta}^{r_l}_l$ denote the components of $\pi_F^{-1}(\delta_l)$. Then we have $r_ln_l=\sigma$ for each $l$, where $n_l$ is the index of the subgroup generated by $\varphi f_{\ast}(d_l)$ and $\varphi f_{\ast}(h)$ in $K$. Then by the Riemann-Hurwitz formula (\cite[pp.\ 133]{Pr-S}, see also \cite[Section 4.2.12]{Pe-S}) we get:  
$$2\tilde{g} = 2 + \sigma \left(2g+p + r - 2 -  \sum_{l=1}^{l=p}\frac{1}{n_l} - \sum_{i=1}^{i=r}\frac{1}{(\mu_i,\beta_i)}\right)$$   where $\tilde{g}$ (resp.\ $g$) denotes the genus of $\tilde{F}$ (resp.\ of $F$, here $g=0$), $p$ denotes the number of boundary components of $F$ and $(\mu_i,\beta_i)$ denotes the greatest common divisor of $\mu_i$ and $\beta_i$. Remark that $n_l\geq 2$ for $l\leq j$. Indeed if $n_l=1$ then $|K:\l\varphi f_{\ast}(h),\varphi f_{\ast}(d_l)\r|=r_l=\sigma=\tau/t=|K,\l\varphi f_{\ast}(h)\r|$. Hence $\l\varphi f_{\ast}(h),\varphi f_{\ast}(d_l)\r=\l\varphi f_{\ast}(h)\r$ which is impossible since $\varphi f_{\ast}(d_l)\not\in\l\varphi f_{\ast}(h)\r$. In particular we have $\sigma\geq 2$.

\medskip
{\bf Case 1}\qua If $j\geq 4$ then $n_l\geq 2$ for $l = 1,...,j$ and so $2\tilde{g} \geq 2 + \sigma(p - p + 4 - 2 - 2) =2$. Thus $\tilde{g}$ $\geq$ 1.

\medskip
{\bf Case 2}\qua If $j=3$ we have 2$\tilde{g}\geq 2 + \sigma(1-\frac{1}{n_1} - \frac{1}{n_2} - \frac{1}{n_3})$ with $n_1, n_2, n_3\geq 2$.

 If $\sigma = 2$ then $n_1=n_2=n_3=2$ and thus $\tilde{g}\geq 1$.

If $\sigma> 2$ then either  $n_l> 2$ for $l= 1,...3$, and thus $\tilde{g}\geq 1$ or there is an element $l$ in $\{1,...3\}$ such that $n_l = 2$. Since $\sigma = n_lr_l$ we have $r_l\geq 2$ and thus $\tilde{S}$ contains at least four boundary components which are non-degenerate and we have a reduction to  {\it Case 1}.
This proves the Lemma.
\end{proof}
\begin{proof}[Outline of proof of Lemma \ref{1etape}] Let $A$ be a Seifert piece of $M$ whose orbit space is a surface of genus $g=2$ (resp.\ $g=1$). We prove here that such a Seifert piece is neccessarily of type I. It follows from the hypothesis of Theorem \ref{toresinjectifs} that 
$f|A : A \to N$ is a non-degenerate map thus using \cite[Mapping Theorem]{Ja-S}  we can change $f$ in such a way that $f(A)$ is contained in a (product) Seifert piece $B$ of $N$. Then combining the fact that $f|\b A$ is non-degenerate and Proposition \ref{technique} we may easily construct a finite (regular) covering of $M$ induced by $f$ from a finite covering of $N$ in which each component of the pre-image of $A$ is a Seifert piece whose orbit space is a surface of genus $g\geq 3$ (resp.\ $g\geq 2$).

Suppose now that the orbit space $F$ of $A$ is a surface of genus 0. It is easily checked that $F$ has at least two boundary components. 
If $A$ has at least three boundary components then it follows easly from Lemma \ref{augmenter} that there is a finite covering of $M$ induces by $f$ from  a finite covering of $N$ in which the lifting of $A$ is a Seifert piece of type I. Thus we may assume that $A$ has exactly two boundary components  (and then $F\simeq S^1\times I$).

If $f_{\ast}(\pi_1(A))$ is non-abelian then we check that   $A$ has at least three boundary components and thus we have a reduction to the ``Type I" case. So suppose now that  $f_{\ast}(\pi_1(A))$ is abelian. Since $f$ is a non-degenerate map and since $f_{\ast}(\pi_1(A))$ is a subgroup of a torsion free three-manifold group it is a free abelian group of rank 2 or 3 (see \cite[Theorem V.I and paragraph V.III]{Ja-S}).  If $f_{\ast}(\pi_1(A)) = {\Z}\times{\Z}\times{\Z}$ then $A$ has at least three boundary components and we have a reduction to the type I case. 
So we may assume that $f_{\ast}(\pi_1(A)) = {\Z}\times{\Z}$. If there is a finite covering $p$ of $M$ induced by $f$ from some finite covering of $N$ such that some component of $p^{-1}(A)$ does not satisfy (i) or (ii) of Lemma 3.1 then using the above argument we show that $A$ is a component of type I, up to finite covering.  If (i) and (ii) of Lemma 3.1 are always checked for any finite covering then $A$ is a component of type II.
\end{proof}
\subsection{Preliminaries for the proof of the second step}
\subsubsection{\label{arrangement}Introduction}
In the following we set $\{ A_i,\  i = 1,...,s(M) \}$ (resp.\ $\{ B_{\alpha},\  \alpha = 1,...,s(N) \}$) the Seifert pieces of a minimal torus decomposition   of  $M$ (resp.\ $N$). On the other hand we will denote by $W_M^S$ (resp.\ $W_N^S$) the canonical tori of $M$ (resp.\ of $N$) which are adjacent on both sides to Seifert pieces of $M$ (resp.\ of $N$).  
We set $A'_i$ = $A_i \setminus W_{M} \times [-1,1]$, for $i = 1,...,s(M)$. Using hypothesis of Theorem \ref{toresinjectifs} and applying the Characteristic Pair Theorem  of \cite{Ja-S} we may assume that for each $i$ there is an $\alpha_i$ such that $f(A'_i) \subset$ int$( B_{\alpha_i})$. Thus if $\Sigma_M$ (resp.\ $\Sigma_N$) denotes the union of the components of $S_M$ (resp.\ $S_N$) with the components $T_i\times[-1,1]$ of $W_M\times[-1,1]$ (resp.\ $W_N\times[-1,1]$) such that $T_i\times\{\pm 1\}\subset\b H_M$ (resp.\ $T_i\times\{\pm 1\}\subset\b H_N$) then $f(\Sigma_M)\subset$ int$(\Sigma_N)$.  Moreover, by identifying a regular neighborhood of $W_M^S$ with  $W_M^S\times I$ we may suppose, up to homotopy, that $f^{-1}(W_N^S)$ is a collection of incompressible tori in $W_M^S\times I$. Indeed since for each $i=1,...,s(M)$ we have $f(A'_i) \subset$ int$( B_{\alpha_i})$ then using standard cut and paste arguments we may suppose, after modifying $f$ by a homotopy which is constant on $\cup A_i'\cup H_M$ that $f^{-1}(W_N^S)$ is a collection of incompressible surfaces in $W_M^S\times I$. Since each component $T_j$ of $W_M^S$ is an incompressible torus then $f^{-1}(W_N^S)$ is a collection of tori parallel to the $T_j$. In the following the main purpose (in the second step) is to prove the following key result.
\begin{lem}\label{etape2} Let $\{A_i, i = 1,...,s(\tilde{M_0})\}$ (resp.\  $\{B_{\alpha}, \alpha = 1,...,s(\tilde{N_0})\}$)  be the Seifert pieces of $S_{\tilde{M_0}}$ (resp.\ of $S_{\tilde{N_0}}$).  Then  $f$ is homotopic to a map $g$ such that:

{\rm(i)}\qua  $g|(H_{\tilde{M_0}},\partial H_{\tilde{M_0}}) : (H_{\tilde{M_0}},\partial H_{\tilde{M_0}}) \to (H_{\tilde{N_0}},\partial H_{\tilde{N_0}})$ is a homeomorphism,

{\rm(ii)}\qua for each $\alpha\in\{1,...,s(\tilde{N_0})\}$ there is a single $i\in\{1,...,s(\tilde{M_0})\}$ such that $f(A_i,\partial A_i) \subset (B_{\alpha},\partial B_{\alpha})$. Moreover the induced maps $f_i=f|A_i\co A_i\to B_{\alpha}$ are ${\Z}$-homology equivalences and  $f|\partial A_i\co\partial A_i\to\partial B_{\alpha}$ is a homeomorphism.
\end{lem}
The proof of this result will be given in paragraph 3.4 using Lemma 3.4 below. In the remainder of Section 3  we will always assume that $(M,N,f)$ is equal to $(\tilde{M}_0,\tilde{N}_0,\tilde{f}_0)$ given by Lemma \ref{1etape}. The goal of this paragraph 3.3 is to prove the following Lemma which simplifies by a homotopy the given map $f$.
\begin{lem}\label{etape2bis} There is a subfamily of canonical tori  $\{T_j, j\in J\}$ in M which cuts $M\setminus H_M$  into graph manifolds $\{V_i,\ i = 1,...,t(M) \leq s(M) \}$ such that: 

{\rm(i)}\qua for each $\alpha_i \in \{ 1,...,s(N) \}$ there is a single $i \in \{ 1,...,t(M) \}$ such that $f$ is homotopic to a map $g$ with $g(V_i,\partial V_i) \subset (B_{\alpha_i},\partial B_{\alpha_i})$. Moreover we have:

{\rm(ii)}\qua   $(V_i,\partial V_i)$ contains at least one Seifert piece of type I,

{\rm(iii)}\qua   $g|\partial V_i : \partial V_i \to \partial B_{\alpha_i}$ is a homeomorphism,

{\rm(iv)}\qua  $g_i=g|(V_i,\partial V_i) : (V_i,\partial V_i) \to (B_{\alpha_i},\partial B_{\alpha_i})$ is a ${\Z}$-homology equivalence.
\end{lem}
 \subsubsection{Some useful lemmas} The proof of Lemma \ref{etape2bis} depends on the following results. In particular  Lemma \ref{description} describes precisely the {\it subfamily of canonical tori}  $\{T_j, j\in J\}$. Here hypothesis and notations are the same as in the above paragraph. The following result is a  consequence of \cite[Main Theorem]{So} and \cite[Lemma 2.11]{Ro}.  
\begin{lem}\label{Rong} There is a homotopy   $(f_t)_{0\leq t\leq 1}$ such that $f_0=f\co M\to N$,  $f_t|\Sigma_M = f|\Sigma_M$ and such that $f_1|(H_M,\partial H_M)\co(H_M,\partial H_M)\to(H_N,\partial H_N)$ is a homeomorphism.
\end{lem}
\begin{proof} Let $T$ be a component of $\b H_N$. We first prove that, up to homotopy fixing $f|\Sigma_M$, we may assume that each component of $f^{-1}(T)$ is a torus which is parallel to a component of $W_M$. Indeed since $f(\Sigma_M)\subset$ int$(\Sigma_N)$ then $f^{-1}(T)\cap\Sigma_M=\emptyset$. On the other hand since $\b\Sigma_M$ is incompressile, then using standard cut and paste arguments (see \cite{Wa}) we may suppose that, up to homotopy fixing $f|\Sigma_M$,  $f$ is transversal to $T$ and that $f^{-1}(T)$ is a collection of incompressible surfaces in $M\setminus\Sigma_M$. The hypothesis of Theorem \ref{toresinjectifs} together with  Theorem \ref{soma}   imply that $H_M\simeq H_N$. Hence we may use similar arguments as those of \cite[Proof of Lemma 2.11]{Ro} to show that each component of $f^{-1}(T)$ is a torus. Thus $f^{-1}(T)$ is a collection of incompresible tori in $H_M\bigcup\cup_j(T_j\times [-1,1])$. Since each incompressible torus in $H_M$ is $\b$-parallel then we may change $f$ by a homotopy  on a regular neighborhood of $H_M$ to push these tori in $\b H_M$. Finally $f^{-1}(\b H_N)$ is made of tori parallel to some components of $W_M$. Each component $E_i$ of $f^{-1}(H_N)$ is a component of $M$ cutted along $f^{-1}(\b H_N)$. Since $f(\Sigma_M)\cap H_N=\emptyset$ then $f^{-1}(H_N)\subset M\setminus\Sigma_M=H_M\cup\bigcup T_i\times I$. Then each component $E_i$ is either a component of $H_M$ or a component $T_j\times [-1,1]$. For each component $H$ of $H_N$ we have deg$\{f^{-1}(H)\stackrel{f|}{\to}H\}=$ deg$(f)=1$. Since a map $T\times [-1,1]\to H$ has degree zero then $f^{-1}(H)$ must contain a component of $H_M$. Since $H_M\simeq H_N$ then $f^{-1}(H)$ contains exactly one component $H'$ of $H_M$ which is sent by $f$ with degree equal to 1. So it follows from \cite[Lemma 1.6]{Th1} that after modifying $f$ by a homotopy  on a regular neighborhood of $H'$ then $f$ sends $H'$ homeomorphically on $H$. We do this for each component of $H_N$. This finishs  the proof of  Lemma 3.5.
\end{proof}
We next prove the following result.
\begin{lem}\label{typeII} Let $A$ be a Type II Seifert piece in  $M$ given by Lemma 3.1 (recall that we have replaced $\t{M}_0$ by $M$). Then we have the following properties:

{\rm(i)}\qua $A$ is not adjacent to a hyperbolic piece in $M$, 

{\rm(ii)}\qua let $S$ be a Seifert piece adjacent to $A$ and let $B$ be the Seifert piece in $N$ such that  $f(S') \subset$  int$(B)$  then necessarily $f(A) \subset$ int$(B)$.
\end{lem}
The proof of this lemma depends on the following  result whose proof is straightfoward.
\begin{lem}\label{trivial}  Let $A$ be a codimension 0 graph  submanifold of $M$ whose boundary is made of a single canonical torus  $T\subset M$ and such that Rk($H_1(A,{\Z}))=1$. If each canonical torus in $A$ separates $M$ then $A$ contains a component which admits a Seifert fibration whose orbit space is the disk  ${D}^2$.
\end{lem}
\begin{proof}[Proof of Lemma \ref{typeII}] We first prove (i).  Let $T_1$ and $T_2$ be the boundary components of $A$. Suppose that there is a hyperbolic piece   $H$ in $M$ which is  adjacent  to  $A$ along $T_1$.
Up to homotopy we know that $f(A') \subset$ int$(B)$ where $B$ is a Seifert piece in $N$, $f(H,\partial H) \subset (H_i,\partial H_i)$ where $H_i$ is a hyperbolic piece in $N$ and that $f|(H,\partial H)\co(H,\partial H)\to(H_i,\partial H_i)$ is a  homeomorphism.
Denote by $W(T_1)$ a regular neighborhood of $T_1$ in $M$. Then  $f(W(T_1))$  contains necessarily one component of $\partial B \cap \partial H_i$ and so $f$ induces a map $f_1\co(A,T_1)\to(B,\partial B)$. Since  $f|(H,\partial H)\co(H,\partial H)\to(H_i,\partial H_i)$ is a  homeomorphism we have found a canonical torus $U$ in $\partial B$ such that $f|T_1\co T_1\to U$ is a homeomorphism.
Recall that $\pi_1(A)$  has a presentation:
 $$\l d_1, d_2, q_1,...q_r, h : [h,q_i] = [h,d_j] = 1,\ \ q_i^{\mu_i} = h^{\gamma_i}, \ \ d_1d_2q_1...q_r = h^b\r$$
and $\pi_1(B)$:  
$$\left\l a_1, b_1,...a_g , b_g, \delta_1, ...\delta_p, t :[t,\delta_k] = [t,a_i] = [t,b_j] = 1,\prod_{i=1}^{i=g}[a_i,b_i]\delta_1...\delta_p = 1\right\r$$ 
with $\pi_1(U) = \l\delta_1,t\r$.
So we get $f_{\ast}(h) = (\delta_1^{\alpha},t^{\beta})$, with $(\alpha,\beta) = 1$.
Let $c_i$ be the homotopy class of an exceptional fiber in $A$ which exists, otherwise $A$ would be homeomorphic to $S^1\times S^1\times I$, which is excluded. So $c_i^{\mu_i}=h$ for some $\mu_i>1$. Since $f_{\ast}(\pi_1(A))$ is isomorphic to ${\Z}\times{\Z}$, we get: $f_{\ast}(c_i) =  (\delta_1^{\alpha_i},t^{\beta_i})$. So we have $\mu_i | (\alpha,\beta)$. This is a contradiction which proves (i).
\end{proof}
Before  continuing the proof of Lemma \ref{typeII} we state the following result.  
\begin{lem}\label{homologique} Let $M$, $N$ be two Haken manifolds and let  $f : M\to N$ be a ${\Z}$-homology equivalence. Moreover we assume that  $M$ and $N$ satisfy the conclusions of Lemma  \ref{1etape}.
If $T$ is a separating canonical torus which is a boundary of a type II Seifert piece 
in $M$ then there exists a finite covering $p$ of $M$ induced by $f$ from a finite covering of $N$ such that some component of  $p^{-1}(T)$ is non-separating.
\end{lem}
\begin{proof} Let $T$  be a separating torus in $M$ and let $X_1$ and $X_2$ be the components of $M \setminus T$. 
We first prove that $H_1(X_1,{\Z}) \not\simeq {\Z}$ and $H_1(X_2,{\Z}) \not\simeq {\Z}$. Suppose the contrary. Thus we may assume that $H_1(X_1,{\Z}) \simeq {\Z}$. It follows from  (i) of Lemma \ref{typeII}, from Lemma 3.1 and from    \cite[Lemma 3.2]{Pe-S} that  $X_1$ is made of Seifert pieces of  Type II. Since  $T = \partial X_1$ is a separating torus in  $M$ then   each canonical torus in $X_1$ separates $M$.
Indeed to see this it is sufficient to prove that if $A$ is a Seifert piece of $X_1$ (of type II) whose a boundary component, say $T_1$ is separating in $M$ then so is the second component of $\b A$, say $T_2$. This fact follows easily from the homological exact sequence of the pair $(A,\partial A)$.
Thus we may apply Lemma \ref{trivial} to $X_1$ which gives a contradiction with the fact that $M$ contains no Seifert piece whose orbit space is a disk.
Hence we get $H_1(X_1,{\Z}) \not\simeq {\Z}$. The same argument shows that $H_1(X_2,{\Z}) \not\simeq {\Z}$.
So to complete the proof it is sufficient to apply arguments of \cite{Pe-S} in paragraph 4.1.4.
\end{proof}
\begin{proof}[End of proof of Lemma \ref{typeII}] We now prove (ii) of Lemma \ref{typeII}.
Let $S$ be a Seifert piece adjacent to $A$ along $T_1$.
 Let $B_S$ and $B_A$ be the Seifert pieces in $N$ such that $f(A') \subset$ int$(B_A)$, $f(S') \subset$ int$(B_S)$ and let $T_1$, $T_2$ be the $\partial$-components of $A$.
If $B_A \not= B_S$, then by identifying a regular neighborhood  $W(T_1)$ of $T_1$ with $T_1 \times [-1,1]$ in such a way that $f(T_1 \times \{-1\}) \subset$ int$(B_A)$ and  $f(T_1 \times \{+1\}) \subset$ int$(B_S)$ we see, using paragraph \ref{arrangement}, that $f(W(T_1))$ must contain a component $U$ of $\partial B_A$. Thus, modifying  $f$ by a  homotopy supported on a regular neighborhood of  $T_1$, we may assume that  $f$ induces a map $f : (A,T_1) \to (B_A,U)$.

\medskip
{\bf Case 1}\qua Suppose first that $T_1$ is  non-separating in $M$. We may choose a simple closed curve  $\gamma$ in $M$ such that $\gamma$ cuts $T_1$ in a single point. Since  $f$ is a  ${\Z}$-homology equivalence  it must preserve intersection number and then we get:   $$[T_1].[\gamma] = \mbox{deg}(f|T_1 : T_1 \to U)\times [U].[f_{\ast}(\gamma)] = 1$$
Hence  deg$(f|T_1 : T_1 \to U) = 1$ and then $f|T_1 : T_1 \to U$ induces an isomorphism $f_{\ast}|\pi_1(T_1) : \pi_1(T_1) \to \pi_1(U)$. Thus we get a contradiction as in the proof of (i) using the fact that   $f_{\ast}(\pi_1(A))$ is abelian.

\medskip
{\bf Case 2}\qua Suppose now that $T_1$ separates $M$ and denote by $X_S$ the component of  $M\setminus T_1$ which contains  $S$ and by  $X_A$ the component of $M\setminus T_1$ which contains  $A$.
Let $p:\t{M}\to M$ be the finite covering of $M$ given by Lemma \ref{homologique} with $T_1$. There is a component   $\tilde{T}$ of $p^{-1}(T_1)$ which is non-separating in $\t{M}$. Let  $\tilde{A}$, $\tilde{S}$ be the Seifert components of $\tilde{M}$ adjacent to  both sides of $\tilde{T}$. Recall that $\tilde{A}$ is necessarily a Seifert piece of type II such that $f_{\ast}(\pi_1(\tilde{A}))$ is abelian (see Lemma 3.1). Let $B_{\tilde{A}}$ (resp.\ $B_{\tilde{S}}$) be the Seifert pieces of $\tilde{N}$ such that $$f(\tilde{A}') \subset \mbox{int}(B_{\tilde{A}}) \ \ \ \ \ \ \ \  f(\tilde{S}') \subset \mbox{int}(B_{\tilde{S}}).$$ Since $B_A \not= B_S$ then  $B_{\tilde{A}} \not= B_{\tilde{S}}$, and thus there is a component  $\tilde{U}$ in $\partial B_{\tilde{A}}$ such that  $\tilde{f}$ induces a map $\tilde{f} : (\tilde{A}, \tilde{T}) \to (B_{\tilde{A}},\tilde{U})$. Since $\tilde{T}$ is non-separating we have a  reduction to case 1. This proves  Lemma 3.6.
\end{proof}
\begin{lem}\label{description} There is a homotopy $(f_t)_{0\leq t\leq 1}$ with $f_0=f$ and $f_t|(H_M,\partial H_M) = f|(H_M,\partial H_M)$ and a collection of canonical tori $\{T_j, j \in J\}\subset W_M^S$ such that:

{\rm(i)}\qua $f_1$ is transversal to  $W_N^S$,

{\rm(ii)}\qua $f_1^{-1}(W_N^S) = \bigcup_{j\in J}T_j$,

{\rm(iii)}\qua the family  $\{T_j, j \in J\}$  corresponds exactly to tori of $W_M^S$ which are adjacent on both sides to Seifert pieces of type I.     \end{lem}
\begin{proof} The proof of (i) and (ii) are similar to paragraphs 4.3.3 and 4.3.6 of \cite{Pe-S}. Thus we only prove (iii). Let $T$ be a component of $W_M^S$ which is adjacent to Seifert pieces of type I denoted by $A_i$, $A_{i'}$ in $M$. Using the same arguments as in paragraph 4.3.7 of \cite{Pe-S} we prove that $T\times [-1,1]$ contains exactly one component of  $f^{-1}(W_N^S)$.

On the other hand if $T$ is the boundary component of a Seifert piece of type II denoted by $A_i$ we denote by    $A_j$ the other Seifert piece adjacent to  $T$. It follows from Lemma \ref{typeII} that  $B_{\alpha_i} = B_{\alpha_j}$. Thus we get $f(A'_i\cup (T\times [-1,1]) \cup A'_j) \subset$ int$(B_{\alpha_i})$, and hence $T\times [-1,1]$ contains no component of  $f^{-1}(W_N^S)$. This completes the proof of Lemma 3.9.
\end{proof}
\subsubsection{End of proof of lemma \ref{etape2bis}}
Let $V_1,...,V_{t(M)}$ be the components of $(M\setminus H_M)\setminus (\cup_{j\in J}T_j)$ where $\{T_j,j \in J\}$ is the family of canonical tori given by Lemma \ref{description}. It follows from Lemma \ref{description} that $f$ induces a map $f_i : (V_i,\partial V_i) \to (B_{\alpha_i},\partial B_{\alpha_i})$. Since deg$(f) = 1$, then the correspondance: $\{1,...,t(M)\}\ni i\mapsto\alpha_i\in\{1,...,s(N)\}$ is surjective.

$(a)$ The fact that the graph manifolds $V_1,...,V_{t(M)}$ contain  some Seifert piece of type I comes from the construction of the  $V_i$ and from Lemma \ref{typeII}. 
Remark that the construction implies that if $A$ is a Seifert piece in  
 $V_i$ such that $\partial V_i\cap\partial A \not=\emptyset$ then $A$ is of Type I (necessarily).

$(b)$ We next show that the correspondence $i\mapsto\alpha_i$ is bijective.
Since $f$ is a degree one map  then to see this it is sufficient to prove that this map is injective. Suppose the contrary. Hence we may choose two  pieces  $V_1$ and $V_2$ which are sent in the same Seifert piece  $B_{\alpha}$ in $N$. If $V_1$ and $V_2$ are adjacent  we denote by $T$ a common boundary component and by  $A_1 \subset V_1$ and $A_2 \subset V_2$ the Seifert pieces (necessarily of type I) adajacent to $T$. Thus by  \cite[Lemma 4.3.4]{Pe-S} we have a contradiction. Thus we may assume that $V_1$ and $V_2$ are non-adjacent. Since deg($f)=1$ we may assume, after re-indexing, that   $f_1 : (V_1,\partial V_1) \to (B_{\alpha},\partial B_{\alpha})$ has non-zero degree and that $f_2 : (V_2,\partial V_2) \to (B_{\alpha},\partial B_{\alpha})$ with $V_1$ and $V_2$ non-adjacent. Moreover, if $A_i^{\star}$ (resp.\ $V_i^{\star}$) denotes the space obtained from $A_i$ (resp.\ $V_i$) by identifying each component of $\b A_i$ (resp.\ $\b V_i$) to a point, we have: Rk$(H_1(A_i^{\star},{\Q})) \leq$ Rk$(H_1(V_i^{\star},{\Q}))$. Since $A_i$ is of Type I, using \cite[Lemma 3.2]{Pe-S}, we get Rk$(H_1(A_i^{\star},{\Q})) \geq 4$ and thus Rk$(H_1(V_i^{\star},{\Q})) \geq 4$.
Thus to obtain a contradiction we  apply the same arguments as in the proof of Lemma 4.3.9 of \cite{Pe-S} to $V_1$ and $V_2$. This proves point (i) of Lemma \ref{etape2bis}.

We now show that we can arrange  $f$ so that $f_i|\b V_i : \b V_i\to\b B_{\alpha_i}$ is a homeomorphism for all  $i$. The above paragraph implies that $f$ induces maps  $f_i : (V_i,\b V_i)\to (B_{\alpha_i},\b B_{\alpha_i})$ such that deg($f_i) = $ deg$(f) = 1$ for all $i$. Thus we need only to show that $f_i$ induces  a one-to-one map from the set of components of $\b V_i$ to the set of components of $\b B_{\alpha_i}$. To see this we apply arguments of paragraph 4.3.12 of \cite{Pe-S}.

Since $f$ is a ${\Z}$-homology equivalence and since $f_i$ is a degree one map and restricts to a homeomorphism on the boundary, by a Mayer-Vietoris argument we see that $f_i$ is a ${\Z}$-homology equivalence for every $i$. This achieves the proof of Lemma \ref{etape2bis}.
\subsection{Proof of Lemma \ref{etape2}}
It follows from Lemma \ref{etape2bis} that to prove Lemma \ref{etape2} it is sufficient to show that any graph manifold $V=V_i$ of $\{V_1,...,V_{t(M)}\}$  contains exactly one Seifert piece  (necessarily of Type I). In fact it is sufficient to prove that $V$ does not contains type II components. Indeed, in this case, if there were two adjacent pieces of type I, they could not be sent into the same Seifert piece in $N$, by an argument made in paragraph 3.3.3. So we suppose that $V$ contains pieces of type II. Then we can find a finite chain $(A_1,...,A_n)$ of Seifert pieces of type II in $V$ such that:

(i)\qua $A_i\subset $ int$(V)$  for $i \in\{1,...,n\}$,

(ii)\qua $A_1$ is adjacent in $V$ to a Seifert piece of type I, denoted by  $S_1$, along a canonical torus $T_1$ of $W_M$ and $A_n$ is adjacent to a Seifert piece of  type I, denoted by $S_n$ in $V$ along a canonical torus $T_n$,

(iii)\qua for each $i\in\{1,...,n-1\}$ the space $A_i$ is adjacent to $A_{i+1}$ along a single canonical torus in $M$.

This means that each Seifert piece of type II in $M$ can be included in a maximal chain of Seifert pieces of type II.
In the following we will denote by $X$ the connected space  $\bigcup_{1\leq i\leq n}A_i$ corresponding to a maximal chain of Seifert pieces of type II in  $V$ and by $B = F\times S^1$ the Seifert piece of $N$ such that $f(V,\b V) = (B,\b B)$.

\begin{rem}\label{bete} In the following we can always assume, using Lemma \ref{homologique},  up to finite covering, that  $M\setminus X$ is connected (i.e.\  $T_1$ is non-separating in $M$).
\end{rem}
In the proof of lemma \ref{etape2} it will be convenient to separate the two following (exclusive) situations:

{\bf Case 1}\qua We assume that  $T_1$ is a non-separating torus in  $V$ (i.e.\ $V\setminus X$ is connected),

{\bf Case 2}\qua We assume that  $T_1$ is a separating torus in $V$ (i.e.\ $V\setminus X$ is disconnected).

We first prove that Case 1 is impossible (see section \ref{non-separant}). We next show (see section \ref{separant}) that in Case 2 there is a finite covering   $p : \t{M}\to M$ induced by $f$ from some finite covering of $N$ such that for each component  $\t{V}$ of $p^{-1}(V)$ the component  $\t{X}$ of $p^{-1}(X)$ which is included in  $\t{V}$ is non-separating in $\t{V}$, which gives a reduction to Case  1.
This will imply that the family $\cal X$ of components of type II in $V$ is empty and then the proof of Lemma \ref{etape2} will be complete.
Before the beginning of the proof we state the following result (notations and hypothesis are the same as in the above paragraph).

\begin{lem}\label{homologique1} Let $V$ be a graph piece in $M$ correponding to the decomposition given by Lemma \ref{etape2bis} and let $X$ be a maximal chain of Seifert pieces of type II in $V$. Then the homomorphism $(i_X)_{\ast} : H_1(\b X,{\Z})\to H_1(X,{\Z})$, induced by the inclusion $\b X\hookrightarrow X$ is surjective.
\end{lem}
\begin{proof} Let $G$ be the space $M\setminus X$ (connected by  Remark \ref{bete}). Since $G$ contains at least one Seifert piece of type I, then using \cite[Lemma 3.2]{Pe-S}, we get Rk$(H_1(G,\h)\geq 6$. Thus the homomorphism $(i_G)_{\ast} : H_1(\b G,\h\to H_1(G,\h$ induced by the inclusion  $i_G : \b G\to G$ is not surjective. Thus there exists a non-trivial torsion group $L_G$ and a  surjective homomorphism:
$$\rho_G : H_1(G,{\Z})\to L_G$$
such that  $(\rho_G)_{\ast}\circ (i_G)_{\ast} = 0$.
On the other hand if we assume that  $(i_X)_{\ast} : H_1(\b X,{\Z})\to H_1(X,{\Z})$ is not surjective, then there is a non-trivial torsion group   $L_X$ and a surjective homomorhism: 
$$\rho_X : H_1(X,{\Z})\to L_X$$
such that $(\rho_X)_{\ast}\circ (i_X)_{\ast} = 0$, where $i_X$ is the inclusion $\b X\to X$. 
Thus using the  Mayer-Vietoris exact sequence of the   decomposition $M = X\cup G$, we get: $\H{M}=\H{G}\oplus\H{X}\oplus{\Z}$ which allows us to construct a surjective homomorphism $$\rho : H_1(M,\h\to L_X\oplus L_G$$ such that $\rho(H_1(G,\h)\not= 0, \ \rho(H_1(X,\h)\not= 0$ and $\rho(H_1(\b X,\h)= 0$.
Let $p : \t{M}\to M$ be the finite covering corresponding to $\rho$. Then $p^{-1}(\b X)$ has  $|L_X\oplus L_G|$ components and each component of  $p^{-1}(G)$ (resp.\ of $p^{-1}(X)$) contains $2|L_X|>2$ (resp.\ $2|L_G|>2$) boundary components. This implies that for each component of $p^{-1}(X)$ the number of boundary components over $T_1$ is $|L_G|>1$, which implies the each component of  $p^{-1}(X)$ contains some Seifert piece which are not of type II. Moreover since $p$ is an abelian covering and since $f$ is a ${\Z}$-homology equivalence then $\t{M}$ is induced by $f$ from a finite covering of $N$. Since $X$ is made of Seifert pieces of type II this contradicts  Lemma \ref{1etape} and proves  Lemma 3.10.
\end{proof}
\subsubsection{\label{non-separant}The ``non-separating" case}
In this section we prove that if  $V\setminus X$ is connected then we get a contradiction. This result depends on the following Lemma:
\begin{lem}\label{Hempel} Let $W(T_1)$ be a regular neighborhood of $T_1$. Then there exists an incompressible vertical torus  $U=\Gamma\times S^1$ in $B\simeq F\times S^1$ where $\Gamma\subset F$ is a simple closed curve and a homotopy $(f_t)_{0\leq t\leq 1}$ such that:

{\rm(i)}\qua $f_0 = f$, the homotopy $(f_t)_{0\leq t\leq 1}$ is equal to $f$ when restricted to $M\setminus W(T_1)$ and $f_1(T_1)=U$,

{\rm(ii)}\qua $\pi_1(U,x) = \l u,t_B\r$ with $x \in f_1(T_1)$,  $u$ is represented by the  curve $\Gamma$ in $F$ and $t_B$ is represented by the fiber of $\pi_1(B,x)$. 
\end{lem}
{\bf Proof}\ \ Denote by $X_1$ the space $f(T_1)$. Since $T_1$ is a non-separating torus in $V$ we can choose a simple closed curve $\gamma$ in int$(V)$ such that:

(i)\qua  $\gamma$ cuts each component of $\partial A_i$, $i = 1,...,n$ transversally in a single point and the other canonical tori of  int$(V)$ transversally,

 (ii)\qua $\gamma$ representes a generator of $H_1(M,{\Z})/{T}(M)$ where ${T}(M)$ is the torsion submodule of $H_1(M,{\Z})$.

\begin{figure}[ht!]
\centerline{
\relabelbox\small
\epsfbox{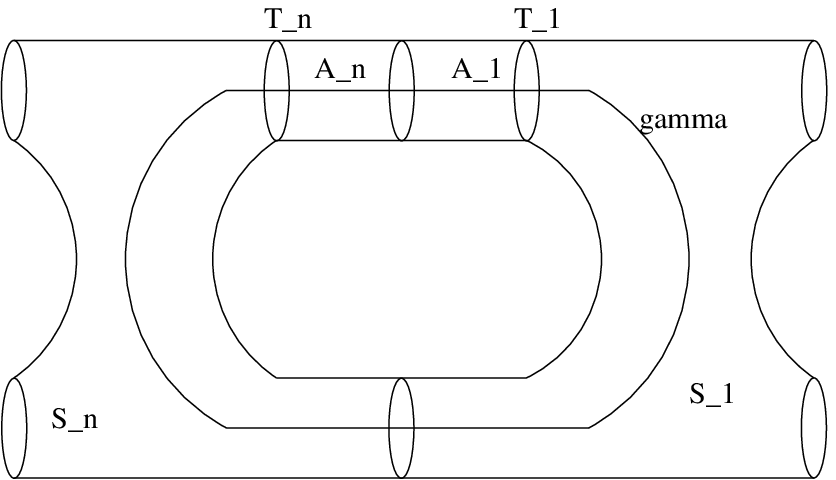}
\relabel {gamma}{$\gamma$}
\relabel {T_1}{$T_1$}
\relabel {T_n}{$T_n$}
\relabel {S_1}{$S_1$}
\relabel {S_n}{$S_n$}
\relabel {A_1}{$A_1$}
\relabel {A_n}{$A_n$}
\endrelabelbox}
\caption{}
\end{figure}

Let $\star$ be a  base point in $T_1$ such that $\gamma \cap T_1 = \{\star\}$ and set $x = f(\star)$. Let $h$ be the homotopy class of the regular fiber of $A_1$ and let $d_1$ be an element in $\pi_1(A_1,\star)$ such that $\l d_1,h\r = \pi_1(T_1,\star)$. We now choose a basis of $\H{M}/{T}(M)$ of type $\{[\gamma],e_2,...,e_n\}$. Since $f$ is a ${\Z}$-homology equivalence then the family $\{f_{\ast}([\gamma]),f_{\ast}(e_2),...,f_{\ast}(e_n)\}$ is a basis of $\H{N}/{T}(N)$. We want to construct an epimorphism $p_1:\H{N}\to{\Z}$ such that $p_1(f_{\ast}([\gamma]))$ is a generator of ${\Z}$ and such that $p_1(f_{\ast}(\l [d_1],[h]\r))=0$.

To do this we choose a basis $\{[\gamma],e_2,...,e_n\}$ of $\H{M}/{T}(M)$ so that $[T_1]\cdot e_i=0$ for $i=2,...,n$. Denote by $i$ the inclusion $T_1\hookrightarrow M$. Since $[T_1]\cdot i_{\ast}(h)=[T_1]\cdot i_{\ast}(d_1)=0$ then $i_{\ast}(h)$ and $i_{\ast}(d_1)$ are in the subspace $K$ of $\H{M}/{T}(M)$ generated by $\{e_2,...,e_n\}$. So it is sufficient to choose $p_1$ equal to the projection of $\H{N}$ on ${\Z}f_{\ast}([\gamma])$ with respect to $f_{\ast}(K)$.       
Denote by $\epsilon$ the following homomorphism: $$\pi_1(N,x) \stackrel{Ab}{\to} H_1(N,{\Z})  \stackrel{p_1}{\to} {\Z}$$
Thus we get an epimorphism  $\epsilon : \pi_1(N,x) \to {\Z}$ such that $\epsilon([f(\gamma )]) = z^{\pm 1}$   where $z$ is a generator of ${\Z}$ and $x = f(\star)$. Since $\pi_1(B,x)$ is a subgroup of $\pi_1(N,x)$ and since  $[f(\gamma )]$ is represented by  $f(\gamma )$ in $B$ then $\epsilon$ induces an epimorphism $\rho_{\ast} = \epsilon|\pi_1(B,x) : \pi_1(B,x) \to {\Z} = \pi_1(S^{1})$ with $\rho_{\ast}([f(\gamma )]) = z^{\pm 1}$ and $\rho_{\ast}(\pi_1(X_1,x)) = 0$ in ${\Z}$. Since $B$ and $S^1$ are both $K(\pi,1)$, it follows from Obstruction theory (see \cite{He1}) that there is a continuous map $\rho : (B,x) \to (S^1,y)$ which induces the above homomorphism and such that $y = \rho(x)$.

The end of proof of Lemma \ref{Hempel} depends on the following result. Notations and hypothesis are the same as in the above paragraph.
\begin{lem} There is a homotopy $(\rho_t)_{0\leq t\leq 1}$ with $\rho_0 = \rho$ such that:

{\rm(i)}\qua $\rho_1(X_1) = \rho_1(f(T_1)) = y$,

{\rm(ii)}\qua $\rho_1^{-1}(y)$ is a  collection of incompressible surfaces in $B$.
\end{lem}
\begin{proof} Since $\rho_{\ast}(\pi_1(X_1,x)) = 0$ in $\pi_1(S^1,y)$ then  the homomorphism $(\rho|X_1)_{\ast} : \pi_1(X_1,x) \to \pi_1(S^1,y)$ factors through $\pi_1(z)$ where $z$ is a 0-simplexe. Then there exist two maps   $\alpha_{\ast} : \pi_1(X_1,x) \to \pi_1(z)$ and $\beta_{\ast} : \pi_1(z) \to\pi_1(S^1,y)$ such that $(\rho|X_1)_{\ast} = \beta_{\ast}\circ \alpha_{\ast}$. Since $z$ and $S^1$ are both  $K(\pi,1)$ then the homomorphisms on $\pi_1$ are induced by maps $\alpha : (X_1,x) \to z$,  $\beta : z \to (S^1,y)$  and $\rho|X_1$ is homotopic to  $\beta\circ \alpha$. Thus we extend this homotopy to  $B$ and we denote by $\rho'$ the resulting map. Then the map $\rho' : (B,x) \to (S^1,y)$ is homotopic to   $\rho$ and  $\rho'(X_1) = y$. This proves point  (i) of the Lemma.

Using  \cite[Lemma 6.4]{He1}, we may suppose that each component of  ${\rho'}^{-1}(y)$ is a surface in $B$.
To complete the proof of the lemma it is  sufficient to show that after changing $\rho'$ by a homotopy fixing  $\rho'|X_1$, then each component of ${\rho'}^{-1}(y)$ is incompressible in $B$. In \cite[pp.\ 60-61]{He1}, J. Hempel proves this point using chirurgical arguments on the map $\rho'$ to get a simplical map $\rho_1$ homotopic to $\rho'$ such that $\rho_1$ is ``simpler" than $\rho'$, (this means that $c(\rho_1) < c(\rho')$ where $c(\rho)$ is the complexity of $\rho$) and inducts on the complexity of  $\rho'$. But these chirurgical arguments can a priori modify the behavior of $\rho'|X_1$. So we will use some other arguments. 
Let $U$ be the component of ${\rho'}^{-1}(y)$ which contains  $f(T_1) = X_1$. Then since $f|T_1 : T_1 \to N$ is non-degenerate the  map $f : (T_1,\star) \to (U,x)$ induces an injective homomorphism $(f|T_1)_{\ast}\co\pi_1(T_1,\star)\to\pi_1(U,x)$. Since $\pi_1(U,x)$ is a surface group then  $\pi_1(U,x)$ has one of the following forms:

(i)\qua a free abelian group of rank $\leq 2$ or,

(ii)\qua a non-abelian free group (when $\partial U \not= \emptyset$)  or,

(iii)\qua a free product with amalgamation of two non-abelian free groups.

Since $\pi_1(U,x)$ contains a subgroup isomorphic to ${\Z}\oplus{\Z}$ then $\pi_1(U,x)\simeq{\Z}\oplus{\Z}$
 and hence $U$ is an incompressible torus in $B$. Note that we necessarily have $f(T_1)=U$. Indeed if there were a point $\star\in U$ such that $f(T_1)\subset U-\{\star\}$ then the two generators free group $\pi_1(U-\{\star\})$ would contain the group $f_{\ast}(\pi_1(T_1))={\Z}\times{\Z}$, which is impossible.
\end{proof}
\begin{proof}[End of proof of Lemma \ref{Hempel}] We show here that $U$ satisfies the  conclusion of Lemma \ref{Hempel}. Since $(\rho')_{\ast}(f_{\ast}(\gamma)) = z^{\pm 1}$ then the intersection number (counted with sign) of $f(\gamma)$ with $U$ is an odd number and then $U$ is a non-separating incompressible torus in $B$. 
Let $t_{S_1}$ be an element of $\pi_1(S_1,\star)$ represented by a regular fiber in $S_1$ and let $t_B \in \pi_1(B,x)$ be represented by the fiber in $B$. Since $S_1$ is a Seifert piece of   Type I, we get $f_{\ast}(t_{S_1}) = t_B^{\alpha}$.
 Indeed, the image of $t_{S_1}$ in $\pi_1(S_1,\star)$ is central, hence the centralizer of $f_{\ast}(t_{S_1})$ in $\pi_1(B,x)$ contains $(f|S_1)_{\ast}(\pi_1(S_1,\star))$ and since $S_1$ is of type I, by the second assertion of \cite[Lemma 4.2.1]{Pe-S} the latter group is non abelian, which implies, using \cite[addendum to Theorem VI.1.6]{Ja-S} that $f_{\ast}(t_{S_1})\in\l t_B\r$.
 Thus $\pi_1(U,x) \supset \l t_B^{\alpha}\r$ i.e.\ $\pi_1(U,x)$ contains an infinite  subgroup which is   central in $\pi_1(B,x)$ and $\pi_1(U,x) \supset {\Z}\oplus {\Z} = f_{\ast}(\pi_1(T_1,\star)).$
Then using   \cite[Theorem  VI.3.4]{Ja} we know that $U$ is a satured torus in $B$, then $\pi_1(U,x) = \l u,t_B\r$ where $u$ is represented by a simple closed curve in $F$. This ends the proof of  Lemma 3.11.
\end{proof}   
\begin{proof}[End of proof of  case 1] It follows from the above paragraph that $[T_1].[\gamma] = [U].[f_{\ast}(\gamma)] = 1$ and thus $f|T_1 : T_1 \to U$ is a degree one map. So $f : (A_1,T_1,\star) \to (B,U,x)$  induces an  isomorphism $f_{\ast} \co \pi_1(T_1,\star)   \to \pi_1(U,x).$
Recall that $\pi_1(A_1,\star)$ has a presentation:
 $$\l d_1, d_2, q_1,...q_r, h: [h,d_i] = [h,q_j] = 1, q_j^{\mu_j} = h^{\gamma_j}, d_1d_2 = q_1...q_rh^b\r$$ where $d_1$ is chosen in such a way that $\pi_1(T_1,\star) = \l d_1,h\r$.
Hence there are two integers  $\alpha$ and $\beta$ such that $f_{\ast}(h) = (u^{\alpha},t_B^{\beta})$ and $(\alpha,\beta) = 1$. Since $f_{\ast}(\pi_1(A_1,\star))$ is an abelian group we have  $f_{\ast}(c_i) = (u^{\alpha_i},t_B^{\beta_i})$ where $c_i$ denotes the homotopy class  of an exceptional fiber in $A_1$. Since $c_i^{\mu_i} = h$ then $\mu_i | (\alpha,\beta).$ This is a contradiction.
\end{proof} 
\subsubsection{\label{separant} The ``separating" case}
We suppose here that $T_1$ is a separating torus in $V$.
We set $X = \bigcup_{1\leq i\leq n}A_i$. Moreover it follows from Remark \ref{bete}, that the space $M\setminus X$ is connected. Let $G$ denote the space  $M\setminus X$ and let $T_1$, $T_n$ be the canonical tori of $M$ such that $T_1\coprod T_n = \b X = \b G$.
Consider the following commutative diagram:
$$\xymatrix{
\H{\b G} \ar[r]^{i_{\ast}} \ar[d]_{\Vert} & \H{S_1\cup S_n} \ar[r]^{j_{\ast}} \ar[d]_{\Vert}  & \H{G} \ar[d]_{\Vert}\\
\H{T_1}\oplus\H{T_n} \ar[r]  & \H{S_1}\oplus\H{S_n} \ar[r] & \H{G}}
$$
Since $S_1$ and $S_n$ are Seifert pieces of type I then Rk$(\H{S_1}\to\H{G})\geq 6$ and  Rk$(\H{S_n}\to\H{G})\geq 6$ (see \cite[Lemma 3.2]{Pe-S}).

\begin{figure}[ht!]
\centerline{
\relabelbox\small
\epsfbox{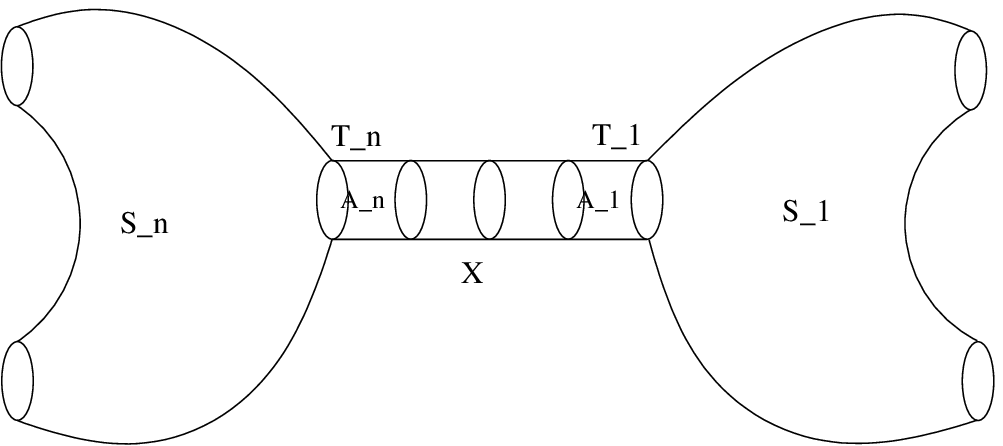}
\relabel {X}{$X$}
\relabel {T_1}{$T_1$}
\relabel {T_n}{$T_n$}
\relabel {S_1}{$S_1$}
\relabel {S_n}{$S_n$}
\relabel {A_1}{$A_1$}
\relabel {A_n}{$A_n$}
\endrelabelbox}
\caption{}
\end{figure}

 So there exists a non-trivial torsion group $L_G$ and an epimorphism:
$$\rho_G : \H{G}\to L_G$$ such that $\rho_G\circ(i_G)\et = 0$, Rk$(\rho_G(\H{S_1})\not= 0$ and Rk$(\rho_G(\H{S_1})\not= 0$, where $i_G$ denotes the inclusion of $\b G$ in $G$ ($(i_G)\et = (j)\et\circ (i)\et$).
It follows from Lemma \ref{homologique1} that the homomorphism
 $(i_X)\et : \H{\b X}\to\H{X}$ is surjective. Then by the Mayer-Vietoris exact sequence of $M = X\cup G$ we get an epimorphism:
$$\rho : \H{M}\to L_G.$$ 
such that $\rho\circ I\et = 0$ and $\rho\circ(i_X)\et = 0$ where $I : \b X\hookrightarrow M$ and $i_X : X\hookrightarrow M$ denote the inclusion and Rk$(\rho_G(\H{S_1}))\not= 0$,  Rk$(\rho_G(\H{S_1}))\not= 0$.

Let $p : \t{M}\to M$ be the finite covering induced by  $\rho$. Since it is an abelian covering and since $f$ is a homology equivalence this covering is induced from a finite covering  $\t{N}$ of $N$. Moreover it follows from the above contruction that  $p^{-1}(X)$ (resp.\ $p^{-1}(G)$) has $|L_G|>1$ (resp.\ 1) components and if  $\t{S}_1$ (resp.\ $\t{S}_n$) denotes a component of  $p^{-1}(S_1)$ (resp.\ of $p^{-1}(S_n)$) then $\b\t{S}_1$ (resp.\ $\b\t{S}_n$) contains at least two components of $p^{-1}(T_1)$ (resp.\ of $p^{-1}(T_n)$).
Let 
 $\t{V}$ be a component of $p^{-1}(V)$ in $\t{M}$ and let $\t{S}_1^1,...,\t{S}_1^{p_1}$ (resp.\ $\t{S}_n^1,...,\t{S}_n^{p_n}$) denote the components of $p^{-1}(S_1)$ (resp.\ $p^{-1}(S_n)$) which are in $\t{V}$.

It follows from the construction of $p$ that each component of  $\t{S}_i^j$ (for $i = 1,n$ and $j\in\{1,...,p_i\}$) has at least two boundary components and the components $\t{X}_1,...,\t{X}_r$ of $p^{-1}(X)\cap\t{V}$ are all homeomorphic to $X$ (i.e.\ the covering is trivial over $X$ because of the  surjectivity of  $\H{\b X}\to\H{X}$).
Let  $\cal A$ denote the submanifold  $\t{V }$ equal to $(\cup_j\t{S}_1^j)\cup(\cup_i\t{X}_i)\cup(\cup_j\t{S}_n^j)$ where we have glued the boundary components of the $\b\t{X}_i$ with the boundary components of the correponding spaces $\t{S}_i^j$.

\begin{figure}[ht!]
\centerline{
\relabelbox\small
\epsfbox{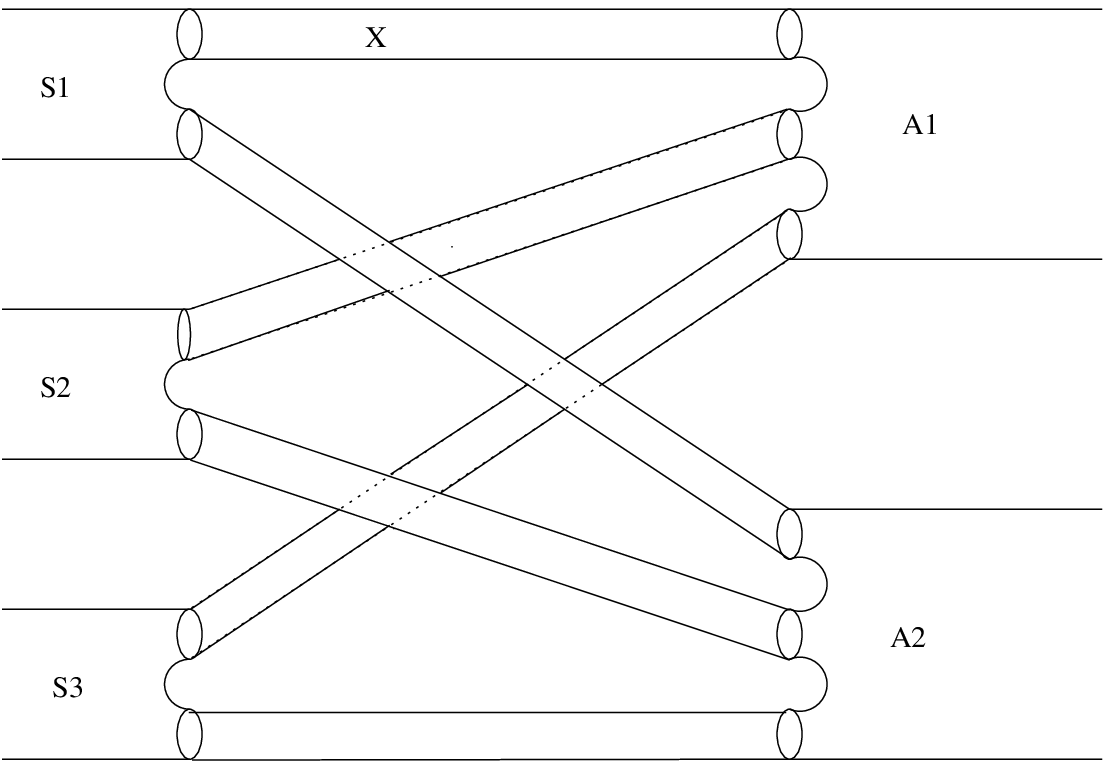}
\relabel {X}{$\t{X}_1$}
\relabel {S1}{$\t{S}^1_n$}
\relabel {S2}{$\t{S}^2_n$}
\relabel {S3}{$\t{S}^3_n$}
\relabel {A1}{$\t{S}^1_1$}
\relabel {A2}{$\t{S}^2_n$}
\endrelabelbox}
\caption{}
\end{figure}

Hence it follows from the construction that there is a submanifold $\t{X}_i$ with a boundary component, say $T_i$, which is non-separating in $\cal A$ (and thus in  $\t{V}$). Let $\t{B}$ be the Seifert piece of  $\t{N}$ such that $\t{f}(\t{V})\subset\t{B}$.
So we can choose a simple closed curve $\gamma$ in $\cal A$ such that $\gamma$ cuts transversally the canonical tori of ${\cal A}$    in at most one point, such that $\t{f}(\gamma)\subset\t{B}$. Thus we have a reduction to the non-separating case. This completes the proof of Lemma 3.3.
\subsection{Proof of the third step} 
We complete here the proof of Theorem \ref{toresinjectifs}. Let $B_{\alpha_i}$ be a Seifert piece of the decomposition of $N$ given by Lemma \ref{etape2} and let $A_i$ be the Seifert piece in $M$ such that $f(A_i,\b A_i)\subset(B_{\alpha_i},\b B_{\alpha_i})$. On the other hand, it follows from Lemma \ref{etape2} that the induced map $f_i=f|(A_i,\b A_i):(A_i,\b A_i)\to (B_{\alpha_i},\b B_{\alpha_i})$ is a ${\Z}$-homology equivalence and the map $f_i|\b A_i:\b A_i\to\b B_{\alpha_i}$ is a homeomorphism. So to complete the proof it is sufficient to show that we can change $f_i$ by a homotopy (rel. $\b A_i$) to a homeomorphism. To see this we first prove that $f_i$ induces an isomorphism  on fundamental groups and we next use \cite[Corollary 6.5]{Wa} to conclude. To prove that maps $f_i$ induce an isomorphism $(f_i)_{\ast}\co\pi_1(A_i)\to\pi_1(B_{\alpha_i})$ we apply arguments of \cite[Paragraphs 4.3.15 and 4.3.16]{Pe-S}. This completes the proof of Theorem \ref{toresinjectifs}. 
\section{Study of the degenerate canonical tori}
 This section is devoted to the proof of Theorem \ref{toresdegeneres}. Recall that the Haken manifold $N^3$ has large first Betti number ($\beta_1(N^3)\geq 3$) and that each Seifert piece in $N^3$ is homeomorphic to a product $F\times{\bf S^1}$ where $F$ is an orientable surface with at least two boundary components.
\subsection{A key lemma for Theorem \ref{toresdegeneres}}
This section is devoted to the proof of the following result.  
\begin{lem}\label{key}  Let $f:M\to N$ be a map satisfying hypothesis of Theorem \ref{toresdegeneres}. If $T$ denotes a degenerate canonical torus in $M$ then $T$ separates $M$ into two submanifolds and there is a component (and only one), say $A$, of $M\setminus T$, such that:

{\rm(i)}\qua $\H{A}={\Z}$,

{\rm(ii)}\qua for any finite covering $p$ of $M$ induced by $f$ from some finite covering of $N$ the components of  $p^{-1}(A)$ have connected boundary,

{\rm(iii)}\qua for any finite covering $p$ of $M$ induced by $f$ from some finite covering of $N$ then each component $\t{A}$ of $p^{-1}(A)$ satisfies $\H{\t{A}}={\Z}$.
\end{lem}
\begin{proof} It follows from \cite[paragraph 4.1.3]{Pe-S} that if $T$ is a degenerate canonical torus in $M$ then $T$ separates $M$ into two submanifolds $A$ and $B$ such that $\H{A}$ or $\H{B}$ is isomorphic to ${\Z}$. Fix notations in such a way that $\H{A}={\Z}$. Note that since $\beta_1(N^3)\geq 3$ then it follows from the Mayer-Vietoris exact sequence of the decomposition $M=A\cup_{T}B$ that $\beta_1(B)\geq 3$. So to complete the proof of Lemma \ref{key} it is sufficient to prove (ii) and (iii). 

We first prove (ii) for regular coverings. Let $\tilde{N}$ be a regular finite covering of $N$ and denote by $\tilde{M}$ the induced finite covering over $M$. Since $p : \tilde{M} \to M$ is regular we can denote by   $k$ (resp.\ $k'$) the number of connected components of $p^{-1}(A)$ (resp.\ $p^{-1}(B)$) and by $p$ (resp.\ $p$') the number of boundary components of  each component of $p^{-1}(A)$ (resp.\ of $p^{-1}(B)$).

\begin{figure}[ht!]
\centerline{
\relabelbox\small
\epsfbox{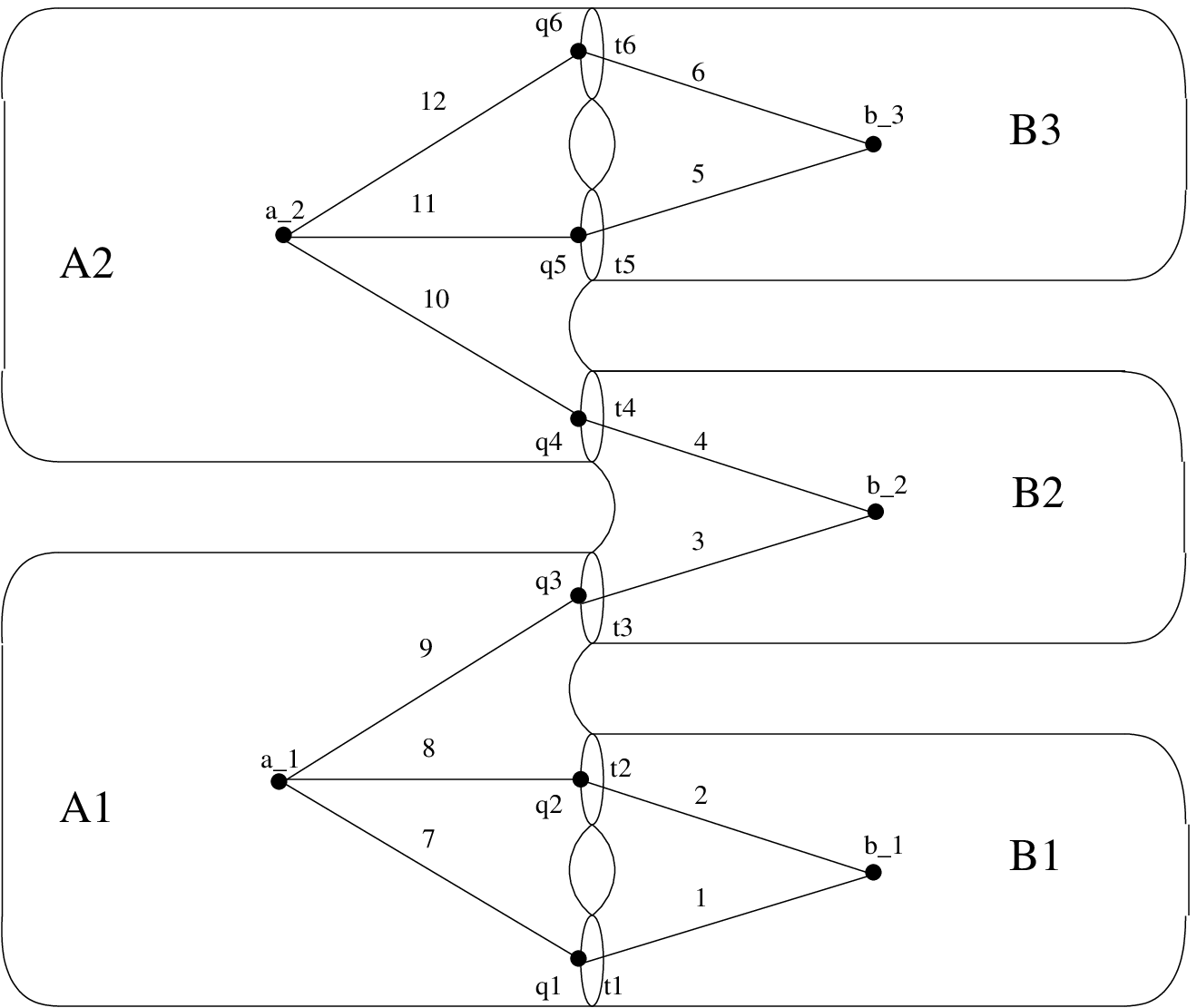}
\relabel {a_1}{$a_1$}
\relabel {a_2}{$a_2$}
\relabel {b_1}{$b_1$}
\relabel {b_2}{$b_2$}
\relabel {b_3}{$b_3$}
\relabel {1}{$\beta_1^1$}
\relabel {2}{$\beta_1^2$}
\relabel {3}{$\beta_2^3$}
\relabel {4}{$\beta_2^4$}
\relabel {5}{$\beta_3^5$}
\relabel {6}{$\beta_3^6$}
\relabel {7}{$\alpha^1_1$}
\relabel {8}{$\alpha^2_1$}
\relabel {9}{$\alpha^3_1$}
\relabel {10}{$\alpha^4_2$}
\relabel {11}{$\alpha^5_2$}
\relabel {12}{$\alpha^6_2$}
\relabel {q1}{$Q_1$}
\relabel {q2}{$Q_2$}
\relabel {q3}{$Q_3$}
\relabel {q4}{$Q_4$}
\relabel {q5}{$Q_5$}
\relabel {q6}{$Q_6$}
\relabel {t1}{$T_1$}
\relabel {t2}{$T_2$}
\relabel {t3}{$T_3$}
\relabel {t4}{$T_4$}
\relabel {t5}{$T_5$}
\relabel {t6}{$T_6$}
\relabel {B1}{$\t{B}_1$}
\relabel {B2}{$\t{B}_2$}
\relabel {B3}{$\t{B}_3$}
\relabel {A1}{$\t{A}_1$}
\relabel {A2}{$\t{A}_2$}
\endrelabelbox}
\caption{}
\end{figure}

Let $\tilde{A}_1,...,\tilde{A}_k$ (resp.\ $\tilde{B}_1,...,\tilde{B}_{k'}$) denote the components of $p^{-1}(A)$ (resp.\break $p^{-1}(B)$). For each $ i = 1,...,k$ ($j = 1,...,k'$) choose a base point $a_i$ (resp.\ $b_j$)  in the interior of each space $\tilde{A}_i$ (resp.\ $\tilde{B}_j$) and choose a base point $Q_l$ in each component of $p^{-1}(T)$ (for $l=1,..., $Card$(p^{-1}(T))$). For each $\tilde{A}_i$ (resp.\ $\tilde{B}_j$) and each component $\tilde{T}_l \subset \partial\tilde{A}_i$ (resp.\ $\tilde{T}_l \subset \partial\tilde{B}_j$) we choose an embedded path $\alpha_i^l$ in $\tilde{A}_i$ joining $a_i$ to $Q_l$ (resp.\ a path $\beta_j^m$ in $\tilde{B}_j$ joining $b_j$ to $Q_m$); we choose these path in such a way that  they don't meet in their interior. Their union is a connected graph denoted by  $\Gamma$.

Then the fundamental group  $\pi_1(\Gamma)$ is a free group with 1 - $\chi$($\Gamma$) generators. In particular $H_1(\Gamma,\Z)$ is the free abelian group of rank 1 - $\chi$($\Gamma$) where $\chi$($\Gamma$) denotes the Euler characteristic of $\Gamma$. Thus we have : $$\chi(\Gamma) = pk + k + k' \ \  \mbox{with} \ \  pk=p'k'$$
So suppose that $p$ and $p'$ $\geq$ 2.
Then we get : $$\chi(\Gamma) \leq  k - k'\ \  \mbox{and} \ \  \chi(\Gamma) \leq k'- k.$$ 
Thus we get $\chi(\Gamma) \leq 0$ and Rk($H_1(\Gamma,\Z)$) $\geq$ 1. Then there exists at least one 1-cycle in $\Gamma$, and thus we can find a component  of  $p^{-1}(T)$ which is a non-separating torus in $\tilde{M}$. So it follows from   \cite[paragraph 4.1.3]{Pe-S} that there exists a canonical torus  $\tilde{T}$ in $p^{-1}(T)$ such that $\tilde{f}|\tilde{T} : \tilde{T} \to \tilde{N}$ is a non-degenerate map. Since $f|T : T \to N$ is a degenerate map, we have a contradiction.

So we can suppose that  $p$ or $p'$  is equal to 1.  So suppose that $p>1$.   Hence we have $p'=1$,   $p^{-1}(A)$ is connected with $p$ boundary components and $p^{-1}(B)$ has $p$ components $\t{B}_1,...,\t{B}_p$ and each of them have connected boundary. Note that since $\beta_1(B)\geq 3$ then it follows from \cite[Lemma 3.4]{Pe-S} that $\beta_1(\t{B}_i)\geq 3$ for $i=1,...,p$. Set $\t{T}_i=\b\t{B}_i$ and $\t{A}_{p-1}=p^{-1}(A)\cup_{\t{T}_1}\t{B}_1\cup_{\t{T}_2}...\cup_{\t{T}_{p-1}}\t{B}_{p-1}$.  It follows easily by a Mayer-Vietoris argument that $\beta_1(\t{A}_{p-1})\geq 2$. So we get a contradiction with the first step of the lemma since $\t{M}=\t{A}_{p-1}\cup\t{B}_p$ and since $\beta_1(\t{A}_{p-1})\geq 2$ and $\beta_1(\t{B}_p)\geq 3$. This proves that $p=1$.

To complete the proof of  (ii) it is sufficient to consider the case of a finite covering (not necessarily regular)  $q \co \tilde{N} \to N$. Then there exists a finite covering $q_1 \co \hat{N} \to \tilde{N}$ such that $\pi_N = q\circ q_1 \co \hat{N} \to N$ is regular. Denote by $p$ (resp.\ $p_1$, resp.\ $\pi_M$) the  covering induced by  $f$ which comes from $q$ (resp.\ $q_1$, resp.\ $\pi_N$). It follows from the above paragraph that each component of $\pi_M^{-1}(A) = q_1^{-1}(q^{-1}(A))$ has connected boundary. So each component of $p^{-1}(A)$ has connected boundary too, which completes the proof of (ii).

We now prove (iii).  So suppose that there is  a finite covering $p:\t{M}\to M$, induced by $f$ from some finite covering of $N$ such that a component $\t{A}$ of $p^{-1}(A)$ satisfies $\H{\t{A}}\not\simeq{\Z}$. Then as in   \cite[paragraph 4.1.4]{Pe-S} we can construct a finite abelian covering $q:\hat{M}\to\t{M}$ in such a way that the components of $q^{-1}(\t{A})$ have at least two boundary components which contradicts (ii). This completes the proof of Lemma \ref{key}. 
\end{proof}
\subsection{Proof of Theorem \ref{toresdegeneres}}
 In the following we denote by  $T$  a canonical torus in $M$ which degenerates under the map  $f \co M \to N$, by $A$ the component of $M\setminus T$ (given by Lemma \ref{key}) satisfying $H_1(A,{\Z}) = {\Z}$ and we set $B = M\setminus A$ with $\beta_1(B)\geq 2$. We will show that the piece $A$ satisfies the conclusion of Theorem \ref{toresdegeneres}.
\subsubsection{Characterization of the non-degenerate components of A} 
 To prove Theorem \ref{toresdegeneres} we will show that each Seifert piece in $A$ degenerates under the map $f$. Suppose the contrary. The purpose of this section is to prove the following result which describes the (eventually) non-degenerate Seifert pieces of $A$.
\begin{lem}\label{describ} Let $T$ be a degenerate canonical torus in  $M$ and let $A$ be the component of $M\setminus T$ such that $H_1(A,{\Z})$ is isomophic to $\Z$. Let $S$ be a Seifert piece in $A$ (using  \cite[Lemma 3.2]{Pe-S} we know that $S$ admits a base of genus 0) such that $f|S : S \to N$ is non-degenerate. Then we get the following properties:

{\rm(i)}\qua there exist exactly two components $T_1, T_2$ of $\partial S$ such that the map $f|T_i:T_i\to N$ is non-degenerate,

{\rm(ii)}\qua $f_{\ast}(\pi_1(S)) = \Z \oplus \Z$,

{\rm(iii)}\qua if p : $\tilde{M} \to M$ denotes a finite covering of $M$ induced by $f$ from some finite covering of  $N$ then each component of $p^{-1}(S)$ satisfies (i) and (ii).
\end{lem}
This result will be used in the paragraph \ref{demo}   to get a contradiction. The proof of Lemma \ref{describ}  depends on the  following  result. 
\begin{lem}\label{aide} Let  $S$ be a Seifert piece in   $M$ whose orbit space is surface of genus 0. Suppose that  $f|S : S \to N$ is a non-degenerate map. Then there exist at least two components $T_1$ and $T_2$ in $\partial S$ such that $f|T_i : T_i \to N$ is non-degenerate.
\end{lem} 
\begin{proof} Let us recall that the group $\pi_1(S)$ has a presentation $(a)$:$$\l d_1,.....,d_p, h, q_1,....., q_r: [h;q_i]=[h;d_j]=1,
 q_i^{\mu_i}=h^{\gamma_i}, d_1...d_pq_1...q_r=h^b\r$$
Since $f|S : S \to N$ is a non-degenerate map, then using    \cite[Mapping Theorem]{Ja-S} we may suppose, after modifying $f$ by a homotopy, that  $f(S)$ is contained in a Seifert piece  $B \simeq F\times S^1$ in $N$.

1. We first show that if the map $f|S:S\to N$ is non-degenerate then $S$ contains at least one boundary component which is non-degenerate under  $f$. To see this, we suppose the contrary:  we will show that if each boundary component  of  $S$ degenerates under $f$ then: $f_{\ast}(\pi_1(S))\simeq\Z$ which gives a contradiction with the definition of  non-degenerate maps  (see \cite{Ja-S}).

Since $f|S :S\to N$ is non-degenerate, we have $f_{\ast}(h) \not= 1$ and then $f_{\ast}(\l d_i,h\r)\simeq{\Z}$. Thus there exist two integers $\alpha_i$ and $\beta_i$ such that 
$$f_{\ast}^{\alpha_i}(d_i) = f_{\ast}^{\beta_i}(h)\ \mbox{and} \ f_{\ast}^{\mu_i}(q_i) = f_{\ast}^{\gamma_i}(h)\eqno{(\star)}$$ 

{\bf Case 1.1}\qua We suppose that the group  $f_{\ast}(\pi_1(S))$ is abelian (remember that the group $f_{\ast}(\pi_1(S))$ is torsion free). Thus it follows from equalities $(\star)$ and from the presentation $(a)$ above  that   $f_{\ast}(\pi_1(S))$ is necessarily isomorphic to the free abelian group of rank 1.

\medskip{\bf Case 1.2}\qua We suppose that the group $f_{\ast}(\pi_1(S))$ is non-abelian. Since  $h$ is central in $\pi_1(S)$, the centralizer $(f|S)_{\ast}(h)$ in $\pi_1(B)$ contains $f_{\ast}(\pi_1(S))$. Since the latter group is non-abelian, it follows from  \cite[addendum to Theorem VI 1.6]{Ja-S}  that $f_{\ast}(h) \in\l t\r$ where $t$ denotes the homotopy class of the regular fiber in $B$. Then  equality $(\star)$ implies that $f_{\ast}(d_i)$ and $f_{\ast}(q_i)$ are in $\l t\r$. Thus using the presentation $(a)$ we get $f_{\ast}(\pi_1(S)) \simeq {\Z}$ which is a contradiction.

2. We show now that if $f|S:S\to N$ is a non-degenerate map then $S$ contains at least two boundary components which are non-degenerate under  $f$. To do this we suppose the contrary. This means that we can assume that
 $f_{\ast}|\l d_1,h\r$ is an injective map and that   $f_{\ast}|\l d_2,h\r$,..., $f_{\ast}|\l d_p,h\r$ are degenerate.  
 
\medskip
{\bf Case 2.1}\qua We suppose that the group    $f_{\ast}(\pi_1(S))$ is abelian. Thus since $$d_1...d_pq_1...q_r=h^b\eqno{(\star \star)}$$ we get Rk($\l f_{\ast}(d_1),f_{\ast}(h)\r$ = 1. This is a contradiction.
 
\medskip
{\bf Case 2.2}\qua  We suppose that the group  $f_{\ast}(\pi_1(S))$ is non-abelian. Since  $h$ is central in $\pi_1(S)$, then by the same argument as in Case 1.2 we get $f_{\ast}(h)$ $\in$ $\l t\r$  where  $t$ the homotopy class of the regular fiber in $B$.
Thus $f_{\ast}(q_i) \in \l t\r$ for $i = 1,...,r$ and  $f_{\ast}(d_j)\in\l t\r$ for $j = 2,...,p$. Then using $(\star \star)$ we get Rk($\l f_{\ast}(d_1),f_{\ast}(h)\r$) = 1. This is  contradiction. This completes the proof of Lemma \ref{aide}.
\end{proof}
\begin{proof}[Proof of Lemma \ref{describ}] Since $S$ is non-degenerate, we denote by $B\simeq F\times{\bf S}^1$ the Seifert piece of $N$ such that $f(S)\subset B$ and by $t$ the (regular) fiber in $B$. Suppose that $S$ contains at least three injective tori in $\partial S$. Denote by $\tilde{N}$ the finite covering of $N$ given by  Lemma \ref{augmenter}. $\tilde{N}$ admits a finite covering $(\hat{N},p)$ which is regular over $N$.
Then each component of the  covering over  $S$ induced from $\hat{N}$  by $f$ admits a Seifert fibration whose orbit space is a surface of genus $\geq$ 1 and then, by regularity,  each component of $p^{-1}(A)$ contains a Seifert piece whose orbit space is a surface of genus $\geq$ 1.

Let $A_1,...A_p$ be the components of $p^{-1}(A)$ and set $\hat{B} = 
p^{-1}(B)$. It follows from Lemma \ref{key} that  $\hat{B}$ is connected and each component 
 $A_i$, $i = 1,...,p$ has a connected boundary.
 Since $\beta_1(A_i) \geq 2$ using \cite[Lemma 3.2]{Pe-S} and $\beta_1(\hat{B})\geq\beta_1(B) \geq 3$,  we get a  contradiction with Lemma \ref{key}. This proves  (i).

Suppose now that the group  $f_{\ast}(\pi_1(S))$ is non-abelian. Since $S$
admits a Seifert fibration over a surface  of genus 0  then $\pi_1(S)$
admits a presentation as in $(a)$ (see the proof of Lemma \ref{aide}). Using   (i) of Lemma 4.2  we may assume that $\l d_1,h\r, \l d_2,h\r$ are injective tori and that $\l d_i,h\r$, $i = 3,...,p$ are degenerate. Then we know that the elements $f_{\ast}(d_i)$ and $f_{\ast}(q_j)$ are in $\l t\r$, (for $i \geq 3$ and $j = 1,...,r$), and then it follows from  $(\star \star)$ that:  
$$f_{\ast}(d_1)f_{\ast}(d_2) \in\l t\r. \eqno{(1)}$$
Since $B$ is a product, we may write : $f_{\ast}(d_1) = (u,t^{\alpha_1})$ and $f_{\ast}(d_2) = (v,t^{\alpha_2}).$
Thus it follows from $(1)$ that  $v = u^{-1}$, and then $f_{\ast}(\pi_1(S))$ is an abelian group. This is a  contradiction.
So $f_{\ast}(\pi_1(S))$ is abelian. Since $f|S\co S\to N$
is a non-degenerate map and since $\pi_1(N)$ is a torsion free group,
$f_{\ast}(\pi_1(S))$ is a finitely generated abelian   free subgroup
of $\pi_1(N)$. Using  \cite[Theorem V.6]{Ja} we know that there exists a compact 3-manifold 
 $V$ and an  immersion $g : V \to N$
such that $g_{\ast} : \pi_1(V) \to \pi_1(N)$ is an isomorphism
onto $f_{\ast}(\pi_1(S))$. Finally $f_{\ast}(\pi_1(S))$ is a free abelian group of rank at least two which is the fundamental group of a 3-manifold.  Then using   \cite[exemple V.8]{Ja} we get that
$f_{\ast}(\pi_1(S))$ is a free abelian group of rank 2 or 3.

Then we prove here that we necessarily have  $f_{\ast}(\pi_1(S)) \simeq {\Z}\oplus{\Z}$. We know that Rk$(\l f_{\ast}(q_j),f_{\ast}(h)\r) = 1$ for $j = 1,...,r$ and by (i) 
Rk$(\l f_{\ast}(d_i),f_{\ast}(h)\r) = 1$ for $i = 3,...,p$ and Rk$(\l f_{\ast}(d_1),f_{\ast}(h)\r) =$ Rk$(\l f_{\ast}(d_2),f_{\ast}(h)\r) = 2$.
Then using $(\star \star)$ and the fact that $f_{\ast}(\pi_1(S))$ is an abelian group, we can find two integers $\alpha, \beta$ such that $f_{\ast}(d_1)^{\alpha}f_{\ast}(d_2)^{\alpha} = f_{\ast}(h)^{\beta}$. This implies that $f_{\ast}(d_2)^{\alpha} \in \l f_{\ast}(d_1),f_{\ast}(h)\r$ and then $f_{\ast}(\pi_1(S)) \simeq {\Z}\oplus {\Z}$. This proves (ii). The proof of (iii) is a direct consequence of (i) and (ii).
\end{proof}
\subsubsection{\label{demo} End of proof of Theorem \ref{toresdegeneres}}
To complete the proof of Theorem \ref{toresdegeneres} it is sufficient to prove (ii). So we first prove that each Seifert piece of  $A$ degenerates and that $A$ is a graph manifold. 
Denote by $S_0$ the component of $A$ which is adjacent to $T=\b A$. It follows from    \cite[Lemma 2]{So} that  $S_0$ is necessarily a Seifert piece of $A$. We prove that  $f|S_0 : S_0 \to N $ is a degenerate map. Suppose the contrary.  Thus $S_0$ satisfies the conclusion of Lemma  \ref{describ}. Let $T_1, T_2$ be the non-degenerate components of  $\b S_0$ and $\pi_1(T_1)=\l d_1,h\r, \pi_1(T_2)=\l d_2,h\r$ the corresponding fundamental groups.
Let $\varphi:\pi_1(N)\to H$ be the correponding epimorphism given by Proposition \ref{technique}, where $H$ is a finite group such that  $\varphi f_\ast(d_1), \varphi f_\ast(d_2)\not\in\l \varphi f_\ast(h)\r$. Denote by $\t{N}$ the (finite) covering given by $\varphi$, $\t{M}$ (resp.\ $\t{S}_0$) the covering of $M$ (resp.\ of $S_0$) induced by $f$. Then formula of paragraph 3.2 applied to $S_0$ and $\t{S}_0$ becomes:
$$2\tilde{g} + \tilde{p} = 2 + \sigma \left(p + r - \sum_{i=1}^{i=r}\frac{1}{(\mu_i,\beta_i)} -2 \right)\eqno{(1)}$$ 
where $\t{p}=\sum_{j=1}^{p}r_j=\sigma\left(\sum_{j=1}^{p}\frac{1}{n_j}\right)$ (resp.\ $p$) is the number of boundary components of the finite covering $\t{S}_0$ of $S_0$ (resp.\ of $S_0$) and where $\tilde{g}$ denotes the genus of the orbit space of  $\t{S}_0$. We can write: $p = 2 + p_1$, where $p_1$ denotes the number of degenerate boundary components of $S_0$ and $\tilde{p} = 2 + \tilde{p_1}$ (where $\tilde{p}_1$ denotes the number of degenerate boundary components of $\tilde{S}_0$). It follows from Lemma  \ref{describ} that we may assume that $\tilde{g}=0.$ Thus using (1), we get: $$\tilde{p_1} 
 = \sigma\left(p_1 + r - \sum_{i=1}^{i=r}\frac{1}{(\mu_i,\beta_i)}\right)$$
Since $(\mu_i,\beta_i)\geq 1$, we have $\tilde{p_1}\geq\sigma p_1$ and then $\tilde{p_1}=\sigma p_1$.
This implies that for each degenerate torus $U$ in  $\b S_0$ there are at least two (degenerate) tori  in $\t{S}_0$ which project onto $U$.
Let us denote by $P\co\tilde{M}\to M$ the finite regular covering of $M$ correponding to $\varphi\circ f_{\ast}$. Then each component of  $P^{-1}(A)$ contains at least two components in its boundary. This contradicts Lemma \ref{key} and so $f|S_0\co S_0\to N$ is a degenerate map. This proves, using 
\cite[Lemma 2]{So} that each component of $A$ adjacent to $S_0$ is a Seifert manifold which allows to apply the above arguments to each of them and prove that they degenerate.   Then we apply these arguments successively to each Seifert piece of  $A$, which proves that $A$ is a graph manifold whose all  Seifert pieces degenerate. 

We now prove that the group $f_{\ast}(\pi_1(A))$ is either trivial or infinite cyclic by induction on the number of Seifert components $c(A)$ of $A$.  
 If $c(A) = 1$ then $A$ admits a Seifert fibration over the disk ${D}^2$. Then the group $\pi_1(A)$ has a presentation: $$\l d_1, h, q_1,...,q_r : [h,d_1] =  [h,q_j] = 1, q_j^{\mu_j} = h^{\gamma_j}, d_1 = q_1...q_rh^b\r$$
We know that $f|A : A\to N$ is a degenerate map. Thus either $f_{\ast}(h) = 1$ or $f_{\ast}(\pi_1(A))$ is isomorphic to
$\{1\}$ or ${\Z}$. So it is sufficient to consider the case  $f_{\ast}(h) = 1$. Since $\pi_1(N)$ is a torsion free group then $f_{\ast}(q_1) = ...=f_{\ast}(q_r) = 1$ and thus
$f_{\ast}(d_1) = f_{\ast}(q_1)...f_{\ast}(q_r)f_{\ast}(h)^b = 1$.
So we have   $f_{\ast}(\pi_1(A)) = \{1\}$.

Let us suppose now that $c(A)>1$. Denote by $S_0$ the Seifert piece  adjacent to
$T$ in $A$ and by $T_1,...,T_k$ its boundary components in int$(A)$. It follows from Lemma \ref{key} that $A\setminus S_0$ is composed of $k$ submanifolds $A_1,...,A_k$ such that $\partial A_i = T_i$ for $i = 1,...,k$. Furthermore, again by Lemma 4.1, $H_1(A_1,{\Z})\simeq ...\simeq H_1(A_k,{\Z})\simeq{\Z}$. Thus the induction hypothesis applies and  implies that $f_{\ast}(\pi_1(A_i)) = \{1\}$ or   $f_{\ast}(\pi_1(A_i)) = {\Z}$ for $i = 1,...,k$.
Let $h_0$ denote the homotopy class of the regular fiber of  $S_0$.

\medskip
{\bf Case 1}\qua Suppose first that $f_{\ast}(h_0) \not = 0$. Since the map $f|S_0 : S_0 \to N$ is degenerate, it follows from the definition that the group $f_{\ast}(\pi_1(S_0))$ is abelian.  Denote by $x_1,...,x_k$   base points in $T_1,...,T_k$. Since $f_{\ast}(\pi_1(A_i))$  is an abelian group, we get the following commutative diagram:
$$\xymatrix{
\pi_1(\partial A_i,x_i) \ar[r]^{i_{\ast}} \ar[d] &   \pi_1(A_i,x_i) \ar[r]^{(f|A_i)_{\ast}} \ar[d] & \pi_1(N,y_i) \ar[d]^{Id}\\
H_1(\partial A_i,{\Z}) \ar[r]^{i_{\ast}} &  H_1(A_i,{\Z})\simeq {\Z} \ar[r] &    \pi_1(N,y_i)
}$$
Since $H_1(A_i,{\Z})\simeq {\Z}$ and since $\partial A_i = T_i$ is connected, then   \cite[Lemma 3.3.$(b)$]{Pe-S} implies that the homomorphism $H_1(\partial A_i,{\Z})\to H_1(A_i,{\Z})$ is surjective and then $$f_{\ast}(\pi_1(A_i,x_i)) = f_{\ast}(\pi_1(T_i,x_i)) \eqno{(\bullet)} $$
Let $(\lambda_i,\mu_i)$ be a base of $\pi_1(T_i,x_i)\subset\pi_1(A_i,x_i)$. Recall that the group $\pi_1(S_0,x_i)$ has a presentation: $$\l d_1,...d_k, d, h_0, q_1,...q_r : [h_0,q_j] = [h_0,d_i] = [h_0,d] = 1,$$
$$ q_i^{\mu_i} = h_0^{\gamma_i}, d_1...d_kd = q_1...q_rh_0^b\r$$
where the element  $d_i$ is chosen in such a way that $\pi_1(T_i,x_i) = \l d_i,h_0\r\subset\pi_1(S_0,x_i)$  for $i = 1,...,k$. Set $A^1 = S_0\cup_{T_1}A_1$ and $A^j = A^{j-1}\cup_{T_j}A_j$ for $j = 2,...,k$ (with this notation we have $A^k = A$). Applying the  the Van-Kampen Theorem     to these decompositions we get:
$$\pi_1(A^1,x_1) = \pi_1(S_0,x_1) \ast_{\pi_1(T_1,x_1)} \pi_1(A_1,x_1)$$
so we get
$$f_{\ast}(\pi_1(A^1,x_1)) = f_{\ast}(\pi_1(S_0,x_1))\ast_{f_{\ast}(\pi_1(T_1,x_1))} f_{\ast}(\pi_1(A_1,x_1))$$
On the other hand it follows from $(\bullet)$ that the injection $f_{\ast}(\pi_1(T_1,x_1))\hookrightarrow f_{\ast}(\pi_1(A_1,x_1))$ is an epimorphism, which implies that the canonical injection $f_{\ast}(\pi_1(S_0))\!\hookrightarrow\! f_{\ast}(\pi_1(S_0,x_1))\ast_{f_{\ast}(\pi_1(T_1,x_1))} f_{\ast}(\pi_1(A_1,x_1))$ is an epimorphism. Thus $f_{\ast}(\pi_1(A^1,x_1))$ is a quotient of the free abelian group of rank 1 $f_{\ast}(\pi_1(S_0,x_1))$ which implies that $f_{\ast}(\pi_1(A^1,x_1))=\{1\}$ or ${\Z}$.  
Applying the same argument with the spaces $A^1$, $A_2$ with  base point $x_2$ we obtain that $f_{\ast}(\pi_1(A^2,x_2))$ is a quotient of   $f_{\ast}(\pi_1(A^1,x_2))$, which implies that  $f_{\ast}(\pi_1(A^2))=\{1\}$ or ${\Z}$. By repeating this method a finite number of times we get:  $f_{\ast}(\pi_1(A))=\{1\}$ or  ${\Z}$.

\medskip
{\bf Case 2}\qua We suppose that $f_{\ast}(h_0)  = 0$. Since $c_i^{\mu_i} = h_0$ (where $c_i$ is any exceptional fiber of $S_0$)  and since
$\pi_1(N)$ is a torion free group, we conclude that $f_{\ast}(\gamma)
= 1$ for every fibers $\gamma$ of $S_0$. Let $F_0$ denote the orbit space 
(of genus 0) of the Seifert fibered manifold $S_0$. Then the map
$f_{\ast}\co\pi_1(S_0)\to\pi_1(N)$ factors through 
$\pi_1(S_0)/\l\mbox{all fibers}\r\simeq\pi_1(F_0)$. Let
$D_1,...,D_n$ denote the boundary components of $F_0$ in such a way that
$[D_i] = d_i \in \pi_1(F_0)$. Then there exist two homomorphisms
$\alpha_{\ast}:\pi_1(S_0)\to\pi_1(F_0)$ and $\beta_ {\ast}:\pi_1(F_0)\to\pi_1(N)$ such that $(f|S_0)_{\ast}=\beta_ {\ast}\circ\alpha_{\ast}$.

We may suppose, after re-indexing, that there exists an integer $n_0
\in\{1,...,k\}$ such that $f_{\ast}(d_1) = ... = f_{\ast}(d_{n_0})
= 1$ and $f_{\ast}(d_j) \not= 1$ for $j = n_0+1,...,k$. If $n_0 = k$ then $f_{\ast}(\pi_1(S_0)) = \{1\}$ and we have a reduction to   Case 1. Thus we may assume that $n_0<n$. Let $\hat{F}_0$ be the 2-manifold obtained from  $F_0$ by gluing a disk ${D}_i^{2}$ along $D_i$ for $i = 1,...,n_0$. The homomorphism $\beta_{\ast} : \pi_1(F_0) \to \pi_1(N)$ factors through the group $\pi_1(\hat{F_0})$. Finally we get two homomorphisms $\hat{\alpha}_{\ast} : \pi_1(S_0) \to \pi_1(\hat{F}_0)$ and $\hat{\beta }_{\ast} : \pi_1(\hat{F}_0) \to \pi_1(N)$ satisfying $(f|S_0)_{\ast}=\hat{\beta}_{\ast}\circ\hat{\alpha}_{\ast}$ where   $\hat{\alpha}_{\ast} : \pi_1(S_0) \to \pi_1(\hat{F}_0)$  is an epimorphism.
It follows from $(\bullet)$ that  $f_{\ast}(\pi_1(A_i)) = \{1\}$ for $ i = 1,...,n_0$ and $f_{\ast}(\pi_1(A_j)) = {\Z}$ for $j = n_0+1,...,k$. Thus the homomorphism $(f|A_i)_{\ast} : \pi_1(A_i) \to \pi_1(N)$ factors through $\pi_1({D}_i^2)$, where ${D}_i^2$ denotes a disk, for  $i = 1,...,n_0$  and the homomorphism $(f|A_j)_{\ast} : \pi_1(A_j) \to \pi_1(N)$ factors through $\pi_1({S}_j^1)$, where ${S}_j^1$ denotes the circle, for  $j = n_0+1,...,k$. 
So we can find two homomorphisms $\pi : \pi_1(A) \to \pi_1(\hat{F}_0)$ and $g : \pi_1(\hat{F}_0) \to \pi_1(N)$ such that $(f|A)_{\ast}=g\circ\pi$ where  $\pi : \pi_1(A) \to \pi_1(\hat{F}_0)$  is an epimorphism. Then consider the following commutative diagram: 
$$\xymatrix{
\pi_1(A) \ar[r]^{\pi} \ar[d] &   \pi_1(\hat{F}_0) \ar[d]\\
H_1(A,{\Z}) \ar[r]^{\hat{\pi}} &  H_1(\hat{F}_0,{\Z})}$$
Since  ${\pi} : \pi_1(A) \to
\pi_1(\hat{F}_0)$ is an epimorphism, then so is
$H_1(A,{\Z}) \to H_1(\hat{F}_0,{\Z})$. Moreover we know that $H_1(A,{\Z})\simeq {\Z}$. Thus we get:
$H_1(\hat{F}_0,{\Z}) \simeq H_1(A,{\Z})\simeq {\Z}$.
Recall that $\pi_1(\hat{F}_0) = \l d_{n_0+1}\r\ast...\ast\l d_{k-1}\r$. Thus
$H_1(\hat{F}_0,{\Z})$ is an abelian free group of rank  $k-1-n_0$ and thus we have: $n_0 = n-2$.
Finally we have proved that $\pi_1(\hat{F}_0) \simeq\l d_{k-1}\r \simeq {\Z}$ which implies that $g_{\ast}(\pi_1(\hat{F}_0))$ is isomorphic to ${\Z}$ and thus $f_{\ast}(\pi_1(A))\simeq {\Z}$. The proof of Theorem \ref{toresdegeneres} is now complete.
\section{Proof of the Factorization Theorem and some consequences} 
This section splits in two parts. The first one (paragraph \ref{Preuve}) is devoted to the proof of Theorem \ref{facto} and the second one gives a consequence of this result (see Proposition \ref{trivialise})  which will be useful in the remainder of this paper.
\subsection{\label{Preuve} Proof of Theorem \ref{facto}} 
The first step is to prove that there exists a finite collection
 $\{T_1,...,T_{n_M}\}$ of degenerate canonical tori satisfying $f_{\ast}(\pi_1(T_i)) =
{\Z}$ in $M$ which define a finite family $\cal{A}$ =
$\{A_1,...,A_{n_M}\}$ of maximal ends  of $M$
such that $\partial A_i = T_i$ and $f|(M\setminus\cup A_i)$ is a non-degenerate map.   We next show  that the map $f \co
M^3 \to N^3$ factors through $M_1$, where $M_1$ is a collapse of $M$ along $A_1,...,A_{n_M}$ and we will see that the map  $f_1 \co M_1^3 \to N^3$, induced by $f$,
satisfies the hypothesis of Theorem \ref{toresinjectifs}. Then the conclusion of Theorem \ref{toresinjectifs} will complete the proof of Theorem \ref{facto}.

\subsubsection{First step}  Let $\{T^0_1,...,T^0_{n_0}\} = W^0_M \subset W_M$ be the canonical tori in $M$ which degenerate under
$f \co M \to N$. If $W^0_M = \emptyset$, by setting $A_i =
\emptyset$, $\pi = f = f_1$ and $M = M_1$ then   Theorem \ref{facto} is obvious by Theorem \ref{toresinjectifs}. So we may assume that $W^0_M \not= \emptyset$.  It follows from  
\cite[Lemma 2.1.2]{So}  that for each component  $T$ of $\partial H_M$, the induced map $f|T \co T \to N$ is $\pi_1$-injective and thus  $W^0_M
\not= W_M$. Then we can choose a degenerate canonical torus 
$T_1$ such that $T_1$ is a boundary component of a Seifert piece  $C_1$ in $M$ which does not degenerate under
$f$. It follows from Theorem \ref{toresdegeneres} that $T_1$ is a separating torus in
$M$. Using Theorem \ref{toresdegeneres} there is a component   $A_1$ of $M\setminus T_1$ such that:

(a)\qua $A_1$ is a graph manifold, $\H{A_1}={\Z}$ and the group  $f_{\ast}(\pi_1(A_1))$ is either trivial or infinite cyclic,

(b)\qua each Seifert piece of $A_1$ degenerates under the map $f$,

(c)\qua $A_1$ satisfies the hypothesis of a maximal end of $M$ (see Definition 1.8).

\medskip
This implies that int$(A_1)\cap$ int$(C_1)=\emptyset$ and $f_{\ast}(\pi_1(A_1))={\Z}$ (if $f_{\ast}(\pi_1(A_1))=\{1\}$, $C_1$ would degenerate under $f$).
Set $B_1 = M\setminus A_1$. If $W^0_{B_1} = \{T^1_1,...,T^1_{n_1}\}
\subset W^0_M$ denotes the family of degenerate canonical tori in int$(B_1)$
then   $n_1 < n_0$. If $n_1 =
0$ we take $\cal{A}$ = $\{A_1\}$. So suppose that $n_1 \geq 1$; we may choose a canonical torus $T_2$ in $W^0_{B_1}$ in the same way as above. Let  $C_2$ denote the non-degenerate Seifert piece  in $M$ such that $T_2 \subset \partial C_2$ and let  $A_2$
be the component of $M\setminus T_2$ which does not meet   int($C_2$). It follows from Theorem \ref{toresdegeneres} that:

$(1)$\qua $A_1 \cap A_2 = \emptyset$,

$(2)$\qua $A_2$ satisfies the above properties  (a), (b) and (c).

Thus by repeating these arguments a finite number of times we get a finite collection  $\{A_1,...,A_{n_M}\}$  of pairwise disjoint maximal ends of  $M$ such that each canonical torus of $M\setminus \bigcup_{1\leq i\leq n_M}A_i$ is non-degenerate.
\subsubsection{Second step}

We next show that the map $f\co M\to N$ factors through a manifold $M_1$ which is obtained from
$M$ by collapsing $M$ along $A_1,...,A_{n_M}$ (see Definition 1.9). To see this it is sufficient to consider the case of a   single maximal end   (i.e.\ $\cal{A}$ = $\{A_1\}$).
Let
$T_1$ be the canonical torus $\partial A_1$ and let $C_1$ be the (non-degenerate) Seifert piece in $M$ adjacent to $A_1$ along $T_1$. Since $f_{\ast}(\pi_1(A_1)) = {\Z}$,  the homomorphism
$f_{\ast} \co \pi_1(A_1) \to \pi_1(N)$ factors through
${\Z}$. Then there are two homomorphisms $(\pi_0)_{\ast}:\pi_1(A_1)\to \pi_1(V_1)$, $(f_0)_{\ast}:\pi_1(V_1)\to\pi_1(N)$
such that $(f|A_1)_{\ast}=(f_0)_{\ast}\circ(\pi_0)_{\ast}$
(where $V_1$ denotes a solid torus) and where $(\pi_0)_{\ast} \co \pi_1(A_1) \to \pi_1(V_1)$ is an epimorphism. Since $V_1$ and $N$ are $K(\pi,1)$,
it follows from Obstruction Theory \cite{He1} that these homomorphisms on $\pi_1$ are induced by  two maps
$\pi_0 \co A_1 \to V_1$ and $f_0 \co V_1\to N$. Moreover we can assume that $f_0$ is an embedding.
We show that we can choose  $\pi_0$ in its homotopy class in such a way that its behavior is sufficiently ``nice". This means that we want that  $\pi_0$ satisfies the  following two conditions:

(i)\qua $\pi_0\co(A_1,\partial A_1)\to(V_1,\partial V_1)$,

(ii)\qua $\pi_0$ induces a homeomorphism $\pi_0|\b A_1\co\b A_1\to\b V_1$.

Indeed since $f_{\ast}(\pi_1(T_1)) = {\Z}$, then there is a basis
$(\lambda,\mu)$ of $\pi_1(T_1)$ such that $(\pi_0)_{\ast}(\lambda) = 1$ in
$\pi_1(V_1)$ and $\l (\pi_0)_{\ast}(\mu)\r= \pi_1(V_1)$. So we may suppose that $\pi_0(\mu) = l_{V_1}$ (resp.\ $\pi_0(\lambda)=m$) where $l_{V_1}$ is a parallel (resp.\ $m$ is a meridian) of  $V_1$. So we have defined a map $\pi_0 \co \partial A_1 \to \partial V_1$ which induces an isomorphism $(\pi_0|\partial A_1)_{\ast} : \pi_1(\partial
A_1)\to\pi_1(\partial V_1)$. So we may assume that   condition (ii)
is checked.
Thus it is sufficient to show that the map $\pi_0|\partial A_0$ can be extended to a map $\pi_0 \co A_1 \to V_1$.
For this consider a handle presentation of  $A_1$ from
$T_1$:$$T_1\cup(e_1^1\cup...e^1_i...\cup e^1_{n_1})\cup(e_i^2\cup...e^2_j...\cup e^2_{n_2})\cup(e_1^3\cup...e^3_k...\cup e_{n_3}^3)$$ where $\{e^k_i\}$ are $k$-cells ($k=1,2,3$).
Since $(\pi_0)_{\ast}(\pi_1(A_1))=(\pi_0)_{\ast}(\pi_1(\b A_1))$, we can extend the map $\pi_0$ defined on $\b A_1$ to the 1-skeletton. Since $V_1$ is a $K(\pi,1)$ space, we can extend $\pi_0$ to $A_1$. Thus, up to homotopy, we can suppose that the map $f\co M\to N$ is such that $f|A_1=f_0\circ\pi_0$, where $\pi_0|\b A_1$ is a homeomorphism.

Set $B_1 = M\setminus A_1$.   Attach a solid torus $V_1$ to $B_1$ along $T_1$ in such a way that the meridian of  $V_1$ is identified to $\lambda$ and the parallel
$l_{V_1}$ of $V_1$ is identified to $\mu$. Let  $\varphi$ denote
the corresponding gluing homeomorphism $\varphi : \partial V_1\to \partial B_1$ and denote by   $\hat{B_1}$ the resulting manifold.
  Let $\pi_0' \co B_1 \to
\hat{B_1}\setminus V_1$ be the  identity map. We define a map $\pi_1  : M = A_1\cup B_1 \to M_1$ such that
$\pi_1|A_1 = \pi_0$ and $\pi_1|B_1 = \pi_0'$ and $M_1 = \hat{B_1}$. Thus it follows from the above construction that $\pi_1 : M
\to M_1$ is a well defined continuous map. Since the map
$\pi_1|B_1\setminus T_1 : B_1\setminus T_1 \to
\hat{B_1}\setminus V_1$ is equal to the  identity, we can define the map
$f_1|\hat{B_1}\setminus V_1$ by setting
$f_1|\hat{B_1}\setminus V_1=f\circ(\pi_1)^{-1}|\hat{B_1}\setminus V_1$ and $f_1 : V_1 \to N$
as $f_0$. Thus we get a map $f_1 : M_1
\to N$ such that $f=f_1\circ\pi_1$.

We now check that $M_1$ is still a Haken manifold of finite volume. Let $\hat{C}_1$ be the space
$C_1\cup_{\varphi}V_1$. Since $M\setminus (A_1\cup C_1)$ is a Haken manifold,   it is sufficient to prove that   $\hat{C}_1$ admits a Seifert fibration. Since $f|C_1 : C_1 \to N$ is a  non-degenerate map,
then $f_{\ast}(h_1) \not= 1$, where $h_1$ denotes the homotopy class of the regular fiber in $C_1$. Therefore the curve $\lambda$
is not a fiber in $C_1$. Thus the Seifert fibration of $C_1$ extends to a Seifert fibration in $\hat{C}_1$.
On the other hand, since $f$ is homotopic to $f_1\circ\pi_1$, we have deg($f_1) =
$ deg($\pi_1$) = deg($f$) = 1 and since $\|N\| =\|M\|$ then $\|N\|
=\|M\| = \|M_1\|$.

In the following, if $A$ denotes a ${\Z}$-module, let ${T}(A)$ (resp.\
${\F}(A)$) be the torsion submodule (resp.\ the free submodule)  of $A$.
To complete the proof of the second step we show that $f_1$ satisfies the homological hypothesis of Theorem \ref{MT}. Let
$q : \tilde{N} \to N$ be a finite cover of $N$, $p :
\tilde{M} \to M$ the  finite covering induced from $\tilde{N}$ by
$f$ and  $p : \tilde{M}_1 \to M_1$ the finite covering induced from
$\tilde{N}$ by $f_1$. Denote by $\tilde{f} : \tilde{M} \to
\tilde{N}$ and $\tilde{f}_1 : \tilde{M}_1 \to \tilde{N}$ the induced maps. Fix base points: $x\in M$, $\tilde{x}\in p^{-1}(x)$, $x_1 =
\pi_1(x)$, $y = f(x)$, $\tilde{y} = \tilde{f}(\tilde{x})$, and $\tilde{x}_1$
such that $\tilde{f}_1(\tilde{x}_1) = \tilde{y}$. In the following  diagram
we first show that there is a map $\tilde{\pi}_1 :
(\tilde{M},\tilde{x}) \to (\tilde{M}_1,\tilde{x}_1)$ such that diagrams  (I) and (II) are consistent.
$$\centerline{
\relabelbox\small
\epsfbox{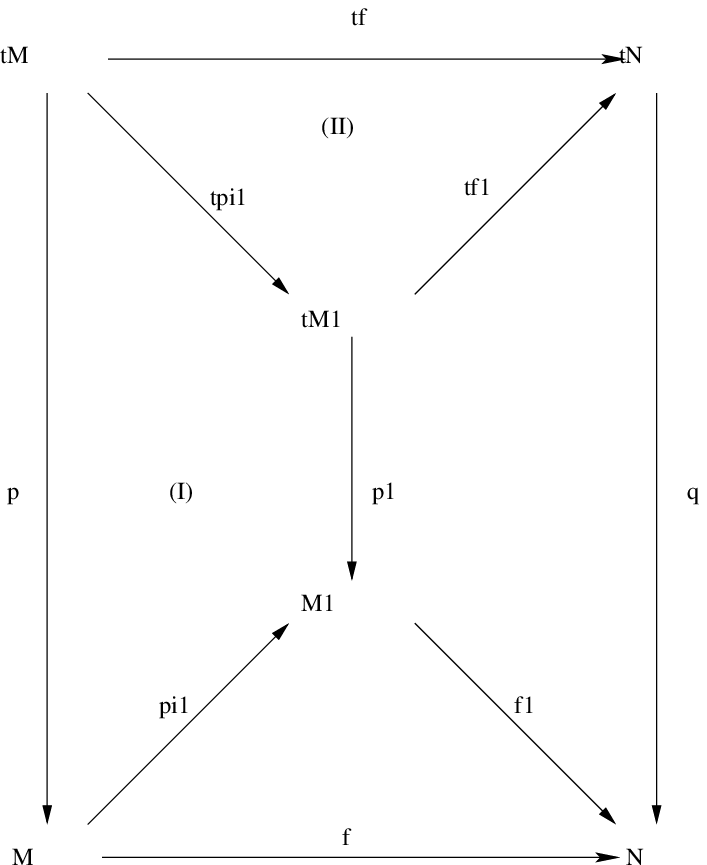}
\adjustrelabel <-3pt,0pt> {M}{$(M,x)$}
\relabel {N}{$(N,y)$}
\adjustrelabel  <-2pt,0pt> {M1}{$(M_1,x_1)$}
\adjustrelabel  <-2pt,0pt> {tM1}{$(\t{M}_1,\t{x}_1)$}
\relabel {tM}{$(\t{M},\t{x})$}
\relabel {tN}{$(\t{N},\t{y})$}
\relabel {f}{$f$}
\relabel {pi1}{$\pi_1$}
\relabel {f1}{$f_1$}
\relabel {p}{$p$}
\relabel {q}{$q$}
\relabel {p1}{$p_1$}
\relabel {tpi1}{$\t{\pi}_1$}
\relabel {tf1}{$\t{f}_1$}
\adjustrelabel <-0pt,-4pt> {tf}{$\t{f}$}
\endrelabelbox}
$$
We know that: $$(p_1)_{\ast}(\pi_1(\tilde{M}_1,\tilde{x}_1)) =
(f_1)_{\ast}^{-1}(q_{\ast}(\pi_1(\tilde{N},\tilde{y}))$$ and 
$$p_{\ast}(\pi_1(\tilde{M},\tilde{x})) =
(f)_{\ast}^{-1}(q_{\ast}(\pi_1(\tilde{N},\tilde{y}))$$ So we get:
$$(\pi_1)_{\ast}(p_{\ast}(\pi_1(\tilde{M},\tilde{x}))) =
(\pi_1)_{\ast}(f)_{\ast}^{-1}q_{\ast}(\pi_1(\tilde{N},\tilde{y})) =
(\pi_1)_{\ast}(\pi_1)_{\ast}^{-1}(f_1)_{\ast}^{-1}q_{\ast}(\pi_1(\tilde{N},\tilde{y}))$$
and thus finally:
$$(\pi_1)_{\ast}(p_{\ast}(\pi_1(\tilde{M},\tilde{x}))) =
(f_1)_{\ast}^{-1}q_{\ast}(\pi_1(\tilde{N},\tilde{y})) =
(p_1)_{\ast}(\pi_1(\tilde{M}_1,\tilde{x}_1))$$ Thus it follows from the lifting criterion, that there is a map
$\tilde{\pi}_1$ such that the  diagram (I) is consistent. Denote by $\hat{f}$
the map $\tilde{f}_1\circ\tilde{\pi}_1$. We easily check  that
$q\circ\hat{f} = f\circ p$ and thus we have  $\tilde{f} = \hat{f}$.
We next show that the maps  $\pi_1$, $f_1$, $\tilde{\pi}_1$ and
$\tilde{f}_1$ induce isomorphisms on $H_1$ (with
coefficients ${\Z}$). Since $f$ (resp.\ $\tilde{f}$) is a 
${\Z}$-homology equivalence  then $(\pi_1)_{\ast} : H_q(M,{\Z})
\to H_q(M_1,{\Z})$ (resp.\ $(\tilde{\pi}_1)_{\ast} :
H_q(\tilde{M},{\Z}) \to H_q(\tilde{M}_1,{\Z})$) is injective
and $(f_1)_{\ast} : H_q(M_1,{\Z}) \to
H_q(N,{\Z})$ (resp.\ $(\tilde{f}_1)_{\ast} : H_q(\tilde{M}_1,{\Z})
\to
H_q(\tilde{N},{\Z})$) is surjective (for $q = 0,...,3$). Since deg$(f)
=$ deg$(\tilde{f}) = 1$ then deg$(\tilde{\pi}_1) =$ deg$(\tilde{f}_1) =$
deg$(\pi_1) =$ deg$(f_1) = 1$. Thus the homomorphism $(\pi_1)_{\ast} :
H_1(M,{\Z}) \to H_1(M_1,{\Z})$ (resp.\ $(\tilde{\pi}_1)_{\ast} :
H_1(\tilde{M},{\Z}) \to H_1(\tilde{M}_1,{\Z})$) is surjective
and therefore is an isomorphism, which implies that $\tilde{f}_1$ and $f_1$
induce  isomorphisms on $H_1$.

We now check that the maps $\tilde{f}_1$ and $\tilde{\pi}_1$
are ${\Z}$-homology equivalences. Recall that $M =
A_1\cup_{T_1}B_1$ and $M_1 = V_1\cup_{T_1}B_1$ where $V_1$ is a 
solid torus and where $\pi_1|(B_1,\partial B_1):(B_1,\partial B_1) \to
(B_1,\partial B_1)$ is  the identity map. On the other hand we see directly that the map $\pi_1|(A_1,\partial A_1) : (A_1,\partial
A_1) \to (V_1,\partial V_1)$ is a  ${\Z}$-homology equivalence
 and deg$(p)$ = deg($p_1$) = deg($q)$.
Set $\tilde{B_1} = p^{-1}(B_1)$ and $\tilde{B}_{1,1} = (p_1)^{-1}(B_1)$.
Since $V_1$ is a  solid torus, it follows from Lemma \ref{key} that:

(i)\qua $\tilde{B_1}$ and $\tilde{B}_{1,1}$ are connected and have the same number  $k_1$ of boundary components,

(ii)\qua $p^{-1}(A_1)$ is composed of $k_1$ connected components 
$\tilde{A_1}^1$,...,$\tilde{A_1}^{k_1}$; $\partial\tilde{A_1}^j$ is
connected; $H_1(\tilde{A_1}^j,{\Z}) = {\Z}$ and
$(p_1)^{-1}(V_1)$ is composed of $k_1$ connected components 
$\tilde{V_1}^1$,...,$\tilde{V_1}^{k_1}$ where the $\tilde{V_1}^j$
are solid tori,

(iii)\qua the map $\tilde{\pi}_1$ induces a map
$\tilde{\pi}_1^j : (\tilde{A_1}^j,\partial\tilde{A_1}^j) \to
(\tilde{V_1}^j,\partial\tilde{V_1}^j)$.\\
Thus we get the two following commutative diagram:
$$\xymatrix{
(\tilde{B}_1,\partial\tilde{B}_1) \ar[r]^{\tilde{\pi}_1|\tilde{B_1}} \ar[d]_{p} & (\tilde{B}_{1,1},\partial\tilde{B}_{1,1}) \ar[d]^{p_1}\\
(B_1,\partial B_1) \ar[r]^{Id} &  (B_1,\partial
B_1)}$$
Since deg$(p|\tilde{B}_1)$ = deg$(p_1|\tilde{B}_{1,1})$ then
deg$(\tilde{\pi}_1|\tilde{B_1}) = 1$ and so the map
$\tilde{\pi}_1|\tilde{B_1}$ is homotopic to a  homeomorphism and is a ${\Z}$-homology equivalence. Consider the following commutative diagram:
$$\xymatrix{
(\tilde{A_1}^j,\partial\tilde{A_1}^j) \ar[r]^{\tilde{\pi}_1} \ar[d]_{p} & 
(\tilde{V_1}^j,\partial\tilde{V_1}^j) \ar[d]^{p_1}\\
(A_1,\partial A_1) \ar[r]^{\pi_1} &  (V_1,\partial V_1)}$$ 
Then we show that we have the following properties:
$$H_1(\tilde{A_1}^j,{\Z}) = {\Z} \ \  \mbox{and} \ \  H_q(\tilde{A_1}^j,{\Z}) = 0 \
\   \mbox{for} \ \   q\geq 2$$
The first identity comes directly from Lemma \ref{key}. On the other hand since  $\partial \tilde{A_1}^j \not=\emptyset$, and since
$\tilde{A_1}^j$ is a 3-manifold, the homology exact sequence of the pair
 $(\tilde{A_1}^j,\partial \tilde{A_1}^j)$ implies that
 $H_3(\tilde{A_1}^j,{\Z}) = 0$. Using  \cite[Corollary 4, p. 244]{Sp} and combining this with Poincar\'e duality, we get: $H_2(\tilde{A_1}^j,{\Z})
\simeq H^1(\tilde{A_1}^j,\partial\tilde{A_1}^j,{\Z})$ and thus ${T}(H_2(\tilde{A_1}^j,{\Z})) = {T}(H_0(\tilde{A_1}^j,\partial\tilde{A_1}^j,{\Z})) = 0$. Moreover , ${F}(H^1(\tilde{A_1}^j,\partial\tilde{A_1}^j,{\Z})) =  {F}(H_1(\tilde{A_1}^j,\partial\tilde{A_1}^j,{\Z}))$ and since
$\beta_1(\tilde{A_1}^j,\partial\tilde{A_1}^j) + 1 = \beta_1(\tilde{A_1}^j) =
1$, we have: $H_2(\tilde{A_1}^j,{\Z}) = 0$.
So the map $\tilde{\pi}_1$ induces an isomorphism on
$H_q(\tilde{A_1}^j,{\Z})$ for $q = 0,1,2,3$. Thus using the Mayer-Vietoris exact sequence of the decompositions $\tilde{M} = (\bigcup_{1\leq
i\leq k_1}\tilde{A}_1^i)\cup(\tilde{B}_1)$ and  $\tilde{M_1} =
(\bigcup_{1\leq i\leq k_1}\tilde{V}_1^i)\cup(\tilde{B}_{1,1})$ we check that the map $\tilde{\pi}_1$ and then $\tilde{f}_1$ are
${\Z}$-homology equivalences.  This proves that $f_1$ satisfies hypothesis of Theorem \ref{MT}. Then using Theorem \ref{toresinjectifs} the proof of  Theorem \ref{facto} is now complete.
\subsection{\label{consequencedefacto} Some consequences of the Factorization Theorem}
We assume here that the manifold $M^3$ contains some canonical tori which degenerate under the map $f$. Then we fix a maximal end $A$ of $M$, whose existence is given by Theorem \ref{facto}. We state here a result which shows that the induced map $f|A$ can be homotoped to a very nice map. More precisely we prove here  Proposition \ref{trivialise}. The proof of this result splits in two lemmas. 
\begin{lem}\label{simpli} If $A$ denotes a maximal end of $M$ then the space   $A\setminus W_M$ contains at least one component, denoted by  $S$,  which admits a Seifert fibration whose orbit space is a disk  ${\bf D}^2$ in such a way that  $f_{\ast}(\pi_1(S))\not=\{1\}$.
\end{lem}
\begin{proof} The fact that the maximal end $A$ contains at least one Seifert piece whose orbit space is a disk (called {\it an extremal component} of $A$) comes directly from Lemma \ref{trivial}  since $A$ is a graph submanifold of $M$   whose Seifert pieces are based on a surface of genus zero and whose canonical tori are separating in   $M$. To prove the second part of Lemma \ref{simpli} we  suppose the contrary. This means that we suppose, for each extremal component  $S$ of $A$, that the induced map  $f|S$ is homotopic in $N$ to a  constant map. Then we show, arguing inductively on the number of connected components of  $A\setminus W_M$, denoted by  $k_A$, that this hypothesis implies that $f_{\ast}(\pi_1(A))=\{1\}$ which gives a contradiction with  Definition \ref{cdm}.

If $k_A=1$, this result is obvious since the component $A$ is a Seifert space whose orbit space is a 2-disk. Then we now suppose that $k_A>1$. The induction hypothesis is the following: 

{\sl If $\hat{A}$ is a  degenerate graph submanifold of $M$   made of $j<k_A$ Seifert pieces and if each Seifert piece $\hat{S}$ of $\hat{A}$ based on a   disk satisfies $f_{\ast}(\pi_1(\hat{S}))=\{1\}$ then the group  $f_{\ast}(\pi_1(\hat{A}))$ is trivial.} 

 Denote by $S_0$ the Seifert piece of $A$ which contains $\b A$, $T_1,...,T_k$ the components of $\b S_0\setminus\b A$ and $A_1,...,A_k$ the connected components of  $A\setminus$ int$(S_0)$ such that $\b A_i=T_i$ for $i=1,...,k$. So we may apply the induction hypothesis to the spaces  $A_1,...,A_k$ which implies that the groups  $f_{\ast}(\pi_1(A_1)),...,f_{\ast}(\pi_1(A_k))$ are trivial. Recall that the group  $\pi_1(S_0)$ has a  presentation:
$$\langle d_1,...,d_k,d,h,q_1,...,q_r: [h,d_i]=[h,q_j]=1,   q_j^{\mu_j}=h^{\gamma_j}, d_1...d_kdq_1...q_r=h^b\rangle$$
where the group $\langle d_i,h\rangle $ is conjugated to $\pi_1(T_i)$ for $i=1,...,k$ and where $\langle d,h\rangle $ is conjugated to $\pi_1(T)$, where $T=\b A$. Since $h$ admits a representative in each component of $\b S_0$ and since $f_{\ast}(\pi_1(T_i))=1$ then $f_{\ast}(h)=1$ and $f_{\ast}(d_1)=...f_{\ast}(d_{k})=1$. This implies that $f_{\ast}(q_1)=...=f_{\ast}(q_r)=1$ and since  $d_1...d_kdq_1...q_r=h^b$ we get $f_{\ast}(d)=1$, which proves that $f_{\ast}(\pi_1(S_0))=1$. Since $A=S_0\cup A_1\cup...\cup A_k$, then applying the   Van Kampen Theorem to this decomposition of $A$, we get $f_{\ast}(\pi_1(A))=\{1\}$ which completes the proof of Lemma \ref{simpli}.         
\end{proof}
\begin{lem} Let $A$ be a maximal end of $M^3$. Let $S$ be a submanifold of $A$ which admits a Seifert fibration whose orbit space is a disk  such that $f_{\ast}(\pi_1(S))\not=\{1\}$. Then there exists a Seifert piece  $B$ of $N$  such that $f_{\ast}(\pi_1(S))\subset \langle t\rangle $, where $t$ denotes the homotopy class of the fiber in $B$.
\end{lem}
\begin{proof} Applying Theorem \ref{facto} to the map $f:M\to N$, we know that  $f$ is homotopic to the comopsition  $f_1\circ\pi$ where $\pi:M\to M_1$  denotes the collapsing map of $M^3$ along its maximal ends and where $f_1:M_1\to N$ is a homeomorphism. More precisely, if $C$ denotes the Seifert piece of $M^3$ adjacent to $A$ along $\b A$  then we know, by the proof of Theorem \ref{facto} that there is a solid torus  $V$ in $M_1$ and a homeomorphism $\phi:\b V\to\b A$ such that:
\begin{itemize}
{\item[\rm(i)] the space $C_1=C\cup_{\varphi}V$ is a Seifert piece in $M_1$,}
{\item[\rm(ii)] $\pi(A,\b A)=(V,\b V)\subset$ int$(B)$ and the map $\pi|\o{M\setminus A}:\o{M\setminus A}\to\o{M_1\setminus V}$ is  the identity.}
\end{itemize}
Since the map $f_1$ is a homeomorphism from $M_1$ to $N$, then by the proof of Theorem \ref{toresinjectifs}, we know that there exists a Seifert fibered space  of $N$, denoted by $B_N$, such that $f_1$ sends $(C_1,\b C_1)$ to $(B_N,\b B_N)$ homeomorphically. Hence the map $f$ is homotopic to the map $f_1\circ\pi$ still denoted  by $f$,  such that $f(A)\subset$ int$(B_N)$ where $B_N$ is a Seifert piece in $N\setminus W_N$. In particular, we have $f(S)\subset$ int$(B_N)$. On the other hand, since $\H{A}={\Z}$ then it follows from    \cite[lemma 5.3.1(b)]{Pe-S}, that the map $\H{\b A}\to\H{A}$, induced by  inclusion, is surjective and since $f_{\ast}(\pi_1(A))$ is an abelian group (in fact isomorphic to ${\Z}$) we get $f_{\ast}(\pi_1(A))=f_{\ast}(\pi_1(\b A))$. Since $f=f_1\circ\pi$, if $h_1$ denotes the homotopy class of the fiber in  $C$ represented in $\b A$, then $f_{\ast}(h_1)=t^{\pm 1}$ where $t$ denotes the homotopy class of the fiber in $B_N$. Moreover, since $f_{\ast}|\pi_1(\b A)$ is a   homomorphism of rank 1 and since $B_N$ is homeomophic to a product  $F_n\times \bf{S}^1$, then we get $f_{\ast}(\pi_1(A))=f_{\ast}(\pi_1(\b A))=\langle t\rangle \subset\pi_1(B_N)\simeq\pi_1(F_n)\times\langle t\rangle $. Finally, since  $\pi_1(S)$ is a subgroup of $\pi_1(A)$ we get $f_{\ast}(\pi_1(S))\subset\langle t\rangle $ which completes the proof of  Lemma 5.2. The proof of Proposition \ref{trivialise} is now complete.
\end{proof}
\section{Proof of Theorem \ref{MT}}
\subsection{Preliminary}
\subsubsection{Reduction of the general problem}
It follows from the form of the hypothesis of Theorem \ref{MT} that to prove this result it is sufficient to find a finite cover $\t{N}$ of $N$ such that the lifting $\t{f}:\t{M}\to\t{N}$ of $f$ is homotopic  to a homeomorphism. So we may always assume without loss of generality that the manifold $N$ satisfies the conclusions of Proposition \ref{rev}. It follows from Theorem  \ref{toresinjectifs} that to prove Theorem \ref{MT} it is sufficient to show that the canonical tori in $M$ do not degenerate under $f$. Thus suppose the contrary: using Theorem \ref{facto} this means that there is a finite collection ${\cal A}=\{A_1,...,A_n\}$ of codimension-0 submanifolds of $M$ which degenerate under $f$ (the maximal ends). We denote by $M_1$ the Haken manifold obtained from $M$ by collapsing along the components of ${\cal A}$, by $\pi:M\to M_1$ the collapsing projection  and by $f_1:M_1\to N$ the homeomorphism such that $f\simeq f_1\circ\pi$.         
Let $A=A_1$ be a maximal end in ${\cal A}$ and let $S$ be a Seifert piece of $A$ whose orbit space is a disk, given by Proposition 1.11. Then the proof of Theorem \ref{MT} will depend on the following result:
\begin{lem}\label{contradiction}There exists a finite covering $p:\t{M}\to M$ induced by $f$ from some finite covering of $N$ such that each component of $p^{-1}(S)$ admits a Seifert fibration whose orbit space is a surface of genus $\geq 1$.
\end{lem}
This result implies that the components of $p^{-1}(A)$ are not maximal ends. Indeed since each component of $p^{-1}(A)$ contains at least one Seifert piece whose orbit space is a surface of genus $\geq 1$ then it follows from   \cite[Lemma 3.2]{Pe-S} that their first homology group is an abelian group of rank $\geq 2$ which contradicts  Definition 1.8. This result gives the desired contradiction.
\subsubsection{Proposition 1.12 implies Lemma \ref{contradiction}.}
In this paragraph we show that to prove Lemma \ref{contradiction}  it is sufficient to prove Proposition \ref{def}.

Let $f:M\to N$ be a map between two Haken manifolds satisfying hypothesis of Theorem \ref{MT}. Let $A$ be a maximal end of $M$ and let $S$ be the extremal Seifert piece of $A$ given by Proposition 1.11 and we denote by $B_N$ the Seifert piece of $N$ such that $f(A)\subset B_N$. Let $h$ (resp.\ $t$) denote the homotopy class of the fiber in $S$ (resp.\ in $B_N$). Then Proposition 1.11 implies that  $f_{\ast}(\pi_1(S))\subset\l t\r$. Recall that the group $\pi_1(S)$ has a presentation: 
$$\l d_1,q_1,...,q_r,h : [h,d_1]=[h,q_i]=1\ \ q_i^{\mu_i}=h^{\gamma_i}\ \ d_1q_1...q_r=h^b\r$$
 Let us denote by $\{\alpha_1,...,\alpha_r\}$ the integers such that $f_{\ast}(c_1)=t^{\alpha_1},...,f_{\ast}(c_r)=t^{\alpha_r}$ where $c_1,...,c_r$ denote the homotopy class of the exceptional fibers in $S$ (i.e.\ $c_i^{\mu_i}=h$). In particular we have $f_{\ast}(h)=t^{\mu_i\alpha_i}$ for $i=1,...,r$. Since the canonical tori in $M$ are incompressible,   the manifold $S$ contains at least two exceptional fibers  $c_1$ and $c_2$ (otherwise $S=D^2\times {\bf S}^1$ which is impossible). Set $n_0 = \alpha_1\alpha_2\mu_1\mu_2$, where $\mu_i$ denotes the index of the exceptional fiber $c_i$. Then, we apply Proposition 1.12 to the manifold $N^3$ with the integer $n_0$ defined as above. Let $\t{B}_N$  be a component of $p^{-1}(B_N)$ in $\t{N}$, where $p$ is the finite covering given by Proposition 1.12. Thus there exists an integer $m$ such that the fiber preserving map $p|\t{B}_N:\t{B}_N\to B_N$ induces the $mn_0$-index covering on the fibers  $\t{t}$  of $\t{B}_N$. Let  $\pi$ denote the homomorphism correponding to the covering induced on the fibers. Thus the covering induces, via  $f$, a regular finite covering $q$ over $S$ which corresponds to the following homomorphism $\theta$: $$\pi_1(S)\stackrel{(f|S)_{\ast}}{\to}{\Z}\simeq\l t\r\stackrel{\pi}{\to}\frac{\Z}{mn_0{\Z}}=\frac{\l t\r}{\l t^{mn_0}\r}$$ 
 Let $\t{S}$ be a component of the covering of $S$ corresponding to $\theta$. Our goal here is to comput the genus of the orbit space, denoted by  $\t{F}$ of $\t{S}$. For each $i\in\{1,...,r\}$, we denote by $\beta_i$ the order of the element $\theta(c_i)=\o{\alpha_i}$  in ${\Z}/{mn_0{\Z}}$. Thus we get the following equalities:
$$\beta_1=m\mu_1\mu_2\alpha_2\ \ \ \beta_2=m\mu_1\mu_2\alpha_1\ \ \mbox{and}\ \ (\beta_1,\mu_1)=\mu_1\ \ \ (\beta_2,\mu_2)=\mu_2$$
Let $\pi_F:\t{F}\to F$ denote the (branched) covering induced by $q$ on the orbit spaces of $\t{S}$ and $S$ and denote by  $\sigma$ the degree of the map $\pi_F$. It follows from Lemma 4.1 (applied to $S$) that each component of   $q^{-1}(S)$ has connected boundary. Using  paragraph 3.2   we know that the genus  $\t{g}$ of $\t{F}$ is given by the following formula:
$$2\t{g}=2+\sigma\left(r-1-\frac{1}{\sigma}-\sum_{i=1}^{i=r}\frac{1}{(\mu_i,\beta_i)}\right)$$
Since $\b\t{S}$ is connected, then  using the above equalities, the last one implies that:
$$2\t{g}\geq 1+\sigma\left(1-\frac{1}{\mu_1}-\frac{1}{\mu_2}\right)\geq 1$$
 which proves that Proposition 1.12 implies Lemma \ref{contradiction}. Hence the remainder of this section will be devoted to the proof of Proposition \ref{def}.
\subsection{Preliminaries for the proof of Proposition \ref{def}}
We assume that $N^3$ satisfies the conclusion of Proposition 1.4. In this section we begin by constructing a class of finite coverings for hyperbolic manifolds. This is the heart of the proof of Theorem \ref{MT}: we use deep results of W. P. Thurston on the theory of deformation of hyperbolic structure. Next (in subsection 6.2.4) we construct  special finite coverings of Seifert pieces, that can be glued to the previous coverings over the hyperbolic pieces, to get a covering of $N^3$ having the desired properties.
\subsubsection{A finite covering lemma for hyperbolic manifolds}
In this paragraph we construct a special class of finite coverings for hyperbolic manifolds (see Lemma 6.2). To state this result precisely we need some notations. Throughout this paragraph we assume that the manifold $N^3$ satisfies the conclusion of Proposition 1.4. 

In this section we deal with a class, denoted by ${\cal H}$ of three-manifolds with non-empty boundary made of pairwise disjoint tori whose interior is endowed with a complete, finite volume hyperbolic structure. Let $H$ be an element of ${\cal H}$ and let $T_1,...,T_h$ be the components of $\b H$. We consider $H$ as a submanifold of the Haken manifold $N$ and we cut $\b H$ in two parts: the first one is made of the tori, denoted  $T_1,...,T_l$,  which are adjacent to Seifert pieces in $N$ and the second one is made of tori, denoted  $U_{1},...,U_r$,  which are adjacent to hyperbolic manifolds along their two sides. For each  $T_i$ ($i\in\{1,...,l\}$)  in $\b H$ we fix generators   $(m_i,l_i)$ of $\pi_1(T_i)\simeq{\Z}\oplus{\Z}$ and we assume these generators are represented by simple smooth closed curves (denoted by $l_i$ and $m_i$ too) meeting transversally at one point and such that  $T_i\setminus (l_i\cup m_i)$ is diffeomorphic to the open disk. The curves $(m_i,l_i)$ will be abusively called  system of   ``longitude-meridian'' (we use notation `` '' as $T_i$ is not the standard torus but a subset of  $N^3$). On the other hand we denote by ${\cal P}^{\ast}$ the set of all prime numbers in ${\N}^{\ast}$ and for each integer $n_0$, we denote by ${\cal P}_{n_0}$ the set:
$${\cal P}_{n_0}=\{n\in{\cal P}^{\ast}\ \ \mbox{such that there is an} \  m\in{\N} \ \mbox{with} \ \ n=mn_0+1\}$$
It follows from  the Dirichlet Theorem (see \cite[Theorem 1, Chapter 16]{I-R}) that for each integer $n_0$ the set ${\cal P}_{n_0}$ is infinite. The goal if this paragraph is to prove the following result:
\begin{lem}\label{hyperbolique1} For each integer $n_0$ and for all but finitely many  primes $q$ of the form   $mn_0+1$, there exists a finite group  $K$, a cyclic subgroup $G_n\simeq{\Z}/n{\Z}$ of $K$,   an element $\o{c}\in G_n$ of multiplicative order $mn_0$,  elements $\o{\gamma}^1,...,\o{\gamma}^l$ in $G_n$ and a homomorphism  $\varphi:\pi_1(H)\to K$ satisfying the following properties:

{\rm(i)}\qua for each $i\in\{1,...,l\}$ there exists an element $g_i\in K$ such that $\varphi(\pi_1(T_i))\subset g_iG_ng_i^{-1}=G_n^i\simeq{\Z}/n{\Z}$ with the following equalities: $\varphi(m_i)=g_i\o{c}g_i^{-1}$ and $\varphi(l_i)=g_i\o{\gamma}_ig_i^{-1}$,

{\rm(ii)}\qua for each $j\in\{1,...,r\}$ the group $\varphi(\pi_1(U_j))$ is either isomorphic  to ${\Z}/q{\Z}$ or to ${\Z}/q{\Z}\times{\Z}/q{\Z}$.
\end{lem}
\subsubsection{Preparation of the proof of Lemma \ref{hyperbolique1}.}
We first recall some results on {\it deformation} of hyperbolic structures for three-manifolds. These results come from chapter 5 of \cite{Th2}. Let $Q$ be a 3-manifold whose interior admits a complete finite volume  hyperbolic structure and whose boundary is made of tori $T_1,...,T_k$. This means that $Q$ is obtained as the orbit space of the action  of a discret, torsion free subgroup $\Gamma$  of ${\cal I}^{+}({\bf H}^{3,+})\simeq PSL(2,{\C})$ on ${\bf H}^{3,+}$  (where ${\bf H}^{3,+}$ denotes the Poincar\'e half space) denoted by $\Gamma/{\bf H}^{3,+}$. Hence we may associate to the complete hyperbolic structure of $Q$ a discret and faithful representation $\overline{H}_0$ (called holonomy) of $\pi_1(Q)$ in $PSL(2,{\C})$ defined up to conjugation by an element of $PSL(2,{\C})$. It follows from Proposition 3.1.1 of \cite{Cu-S}, that this representation  lifts to a faithful representation denoted by  $H_0:\pi_1(Q)\to SL(2,{\C})$. Note that since $Q$ has finite volume, the representation $H_0$ is necessarily irreducible. Moreover, since $H_0$ is faithful, then for each component $T$ of $\b Q$ and for each element $\alpha\in\pi_1(T)$, the matrix $H_0(\alpha)$ is conjugated to a matrix of the form:
$$\left(\begin{array}{clcr}1 & \lambda \\0 & 1\end{array}\right) \ \mbox{where} \ \lambda\in{\C}^{\ast}$$
We will show here that for each primitive element $\alpha\in\pi_1(T)$, there exists a neighborhood $W$ of $1$ in ${\C}^{\ast}$ such that for all $z\in W$ there exists a representation $\rho:\pi_1(Q)\to SL(2,{\C})$ such that one of the eigenvalues of $\rho(\alpha)$ is equal to $z$.

Denote by $R(\pi_1(Q))$ the affine algebraic variety of representations of $\pi_1(Q)$ in $SL(2,{\C})$ (i.e.\ $R(\pi_1(Q))=\{\rho, \rho:\pi_1(Q)\to SL(2,{\C})\}$) and by $X(\pi_1(Q))$ the space of characters of the representations  
of $\pi_1(Q)$. For each element $g\in\pi_1(Q)$ we denote by $\tau_g$ the map defined by:
$$\tau_g:R(\pi_1(Q))\ni\rho\mapsto\mbox{tr}(\rho(g))\in{\C}$$
Let ${\cal T}$ denote the ring generated by all the functions $\tau_g$ when $g\in\pi_1(Q)$. Since $\pi_1(Q)$ is finetely generated, then so is the ring ${\cal T}$; so we can choose a finite number of elements $\gamma_1,...,\gamma_m$ in $\pi_1(Q)$ such that $\langle\tau_{\gamma_1},...,\tau_{\gamma_m}\rangle={\cal T}$. We define now the map $t$ in the following way:
$$t:R(\pi_1(Q))\ni\rho\mapsto(\tau_{\gamma_1}(\rho),...,\tau_{\gamma_m}(\rho))\in{\C}^m$$
which allows to identify the space of characters $X(\pi_1(Q))$ with $t(R(\pi_1(Q)))$. In particular, if $R_0$ denotes an irreducible component of $R(\pi_1(Q))$ which contains $H_0$, then the space $X_0=t(R_0)$ is an affine algebraic variety called {\it deformation space} of $Q$ {\it near the initial structure} $H_0$. It follows from \cite[Theorem 5.6]{Th2}, or \cite[Proposition 3.2.1]{Cu-S}, that if $Q$ has $k$ boundary components (all homeomorphic to a torus), then dim($X_0$) = dim($R_0)-3\geq k$. We now fix basis of ``meridian-longitude" $(m_i,l_i)$, $1\leq i\leq k$, for each torus $T_1,...,T_k$. This allows us to define a map:$$\mbox{tr}:X_0\to{\C}^k$$ in the following way: let $q$ be an element in $X(\pi_1(Q))$. The above paragraph implies that there exists a representation $H_q$ such that $t(H_q)=q$; then we set $tr(q)=(tr(H_q(m_1)),...,tr(H_q(m_k)))\in{\C}^k$. By construction this map is a well defined polynomial map between the affine algebraic varieties $X(\pi_1(Q))$ and ${\C}^k$. Moreover, if $q_0$ denotes the element of $X_0$ equal to $t(H_0)$, then it follows from the Mostow Rigidity Theorem (see \cite[Chapter C]{Be-P}) that the element $q_0$ is the only point in the inverse image of tr$(q_0)$. Using \cite[Theorem 3.13]{Mu} this implies that dim($X_0)=$ dim$({\C}^k)=k$. Next, applying the Fundamental Openness  Principle (see \cite[Theorem 3.10]{Mu}) we know that there exists a neighborhood $U$ of $q_0$ in $X_0$ such that tr$(U)$ is a neighborhood of $tr(q_0)=(2,...,2)$ in ${\C}^k$ denoted by $V$. Let $f$ be the map defined by $f(z_1,...,z_k)=(z_1+1/z_1,...,z_k+1/z_k)$ and let $W$ be a neighborhood of $(1,...,1)$ in ${\C}^k$ such that $f(W)=$ tr($U)=V$. This proves that for each $k$-uple $\lambda=(\lambda_1,...,\lambda_k)$ of $W$ there exists a representation $H_{\lambda}$ of $\pi_1(Q)$ in $SL(2,{\C})$ such that for each $i\in\{1,...,k\}$ the matrix $H_{\lambda}(m_i)$ has an eigenvalue equal to $\lambda_i$.
\subsubsection{Proof of  Lemma \ref{hyperbolique1}}
Let $H$ be a submanifold of $M^3$ which admits a complete  finite volume  hyperbolic structure $q_0$. We denote by $H_0$ the irreducible holonomy of $\pi_1(H)$ in $SL(2,{\C})$ associated to the complete structure of $H$, by $R_0$ an irreducible component of $R(\pi_1(H))$ which contains $H_0$ and by $X_0$ the component of $X(\pi_1(H))$ defined by $X_0=t(R_0)$ (see paragraph 6.2.2 for definitions). Let $U_1,...,U_r$ be the components of $\b H$ which bound hyperbolic manifolds of $M^3$ along their both sides and let $T_1,...,T_l$ be the components of $\b H$ adjacent to Seifert pieces. For each $T_i$, $i=1,...,l$ (resp.\  $U_j$, $j=1,...,r$), we fix a system of ``longitude-meridian" $(m_i,l_i)$ (resp.\ $(\mu_j,\nu_j)$). Let $\lambda$ be a transcendental number (over ${\Z}$), near of $1$ in ${\C}$ (this is possible since the set of algebraic number over ${\Z}$ is countable). It follows from paragraph 6.2.2 that there is a representation $H_q$ of $\pi_1(Q)$ in $SL(2,{\C})$ satisfying the following equalities:
$$vp(H_q(\mu_j)) = v_j(q) = vp(H_q(\nu_j)) = 1\ \ \mbox{for}\ \ j\in\{1,...,r\}$$
$$vp(H_q(m_i)) = \lambda_i(q) = \lambda \ \ \ \mbox{for} \ i\in\{1,...,l\}$$
where $vp(A)$ denotes one of the eigenvalues of the matrix $A\in SL(2,{\C})$. Thus we get the following equalities:
$$H_q(m_i) = Q_i \left( \begin{array}{clcr}\lambda & 0\\0 &
\lambda^{-1}\end{array} \right) Q_i^{-1},\ \ H_q(l_i) = Q_i\left(
\begin{array}{clcr}\mu_i & 0\\0 & (\mu_i)^{-1}\end{array} \right)Q_i^{-1}$$
for $i\in\{1,...,l\}$ where $\lambda$ is a transcendental number over ${\Z}$ and where the matrix $Q_1,...,Q_l$ are in $SL(2,{\C})$. On the other hand, since $vp(H_q(\mu_j)) = v_j(q) = vp(H_q(\nu_j)) = 1$ for
$j = 1,...,r$, the groups $H_q(\pi_1(U_j))$  are unipotent and isomorphic to
${\Z}\oplus{\Z}$. This implies that:
$$P_jH_q(\mu_j)P_j^{-1} =
\left( \begin{array}{clcr}1 & 1\\0 & 1\end{array} \right), \ \
P_jH_q(\nu_j)P_j^{-1} =  \left( \begin{array}{clcr}1 & \eta_j\\0 & 1\end{array}
\right)\ \
\mbox{for}\  j\in\{1,...,r\}$$  where $\eta_{1},...,\eta_r$  are in
${\C}\setminus{\Q}$ and where the  $P_j$'s  are in $SL(2,{\C})$. Moreover since $\pi_1(H)$ is a finitely generated group, we may choose a finite susbet  $\cal G$ which generates $\pi_1(H)$.
Then consider the subring ${\A}$ of ${\C}$, generated over ${\Z}[\lambda]$   by the following system: $${\Z} \bigcup
\{\mbox{entries\ of\ the\  matrix}\ H_q(g)\ \mbox{for}\ g\in {\cal
G}\}$$ $$\bigcup \{\mbox{entries\ of\ the\  matrix}\ P_j,\ P_j^{-1},\
Q_i,\ (Q_i)^{-1} \} \bigcup \{{\lambda}, {\lambda}^{-1},
{\mu}_i, ({\mu}_i)^{-1},\ \eta_j\}$$
It follows from the above construction that ${\A}$ is a finitely generated ring over  
${\Z}[{\lambda}]$, and ${\Z}[{\lambda}]$  is isomorphic to ${\Z}[X]$ since $\lambda$ is transcendental over ${\Z}$. So using the Noether Normalization Lemma  (see \cite{Ei}, Theorem 3.3) for the ring ${\A}$ over $B_0={\Z}[\lambda]$, we know that there exists a polynomial  $P$ of ${\Z}[X]$ and a finite algebraically free family $\{x_1,...,x_k\}$ over ${\Z}[{\lambda}]$ such that ${\A}$ is integral over $$B =
{\Z}[{\lambda}]\left[\frac{1}{P(\lambda)}\right][x_1,...,x_k]$$
To complete the proof of Lemma \ref{hyperbolique1} we need the following result.
\begin{lem}\label{technique1} Let $n_0>0$ be a fixed integer. Let ${\A}$ be a  subring of ${\C}$ integral over a ring $B$ isomorphic to ${\Z}[{\lambda}]\left[{1}/{P(\lambda)}\right][x_1,...,x_k]$, where $\lambda$ is transcendental over ${\Z}$,  $P$ is a polynomial in ${\Z}[\lambda]$  and  $x_1,...,x_k$ are algebraically free over ${\Z}[\lambda]\simeq{\Z}[X]$. Let $\mu_1,...,\mu_l$ be elements of ${\A}$. Then for all  $n_0\in{\N}$ and for all but finitely many  primes $q=mn_0+1$, there is a finite field ${\F}_q$ of characteristic $q$, an element $\o{c}$ in $({\F}_q)^{\ast}={\F}_q\setminus\{0\}$ of multiplicative order $mn_0$,  elements $\o{\gamma}_1,...,\o{\gamma}_l$ in $({\F}_q)^{\ast}$ and a ring homomorphism $\varepsilon:{\A}\to{\F}_q$ such that:

{\rm(i)}\qua $\l\l\varepsilon(\lambda),\varepsilon(\mu_i)\r\r\ \subset{\F}_q^{\ast}\simeq{\Z}/n{\Z}$ for $i=1,...,l$, where $\l\l g,h\r\r$ is the multiplicative subgroup of ${\A}$ generated by $g$ et $h$ and where  $n=|{\F}_q|-1$, 

{\rm(ii)}\qua $\varepsilon(\lambda)=\o{c}$ and $\varepsilon(\mu_i)=\o{\gamma}_i$ for $i=1,...,l$.
\end{lem}
\begin{proof}[Proof of Lemma  \ref{technique1}] We first prove that for all but finitely many  primes  $q\in{\cal P}_{n_0}$ there exists a ring homomorphism $\varepsilon:B\to{\Z}/q{\Z}$ such that $\varepsilon(\lambda)$ is a generator of the cyclic group $({\Z}/q{\Z})^{\ast}\simeq {\Z}/(q-1){\Z}$. We next show that we can extend  $\varepsilon$ to a homomorphism from ${\A}$ by taking some finite degree extension of ${\Z}/q{\Z}$ and using the fact that  ${\A}$ is integral over $B$. To this purpose we claim that for all but finitely many  primes  $q=mn_0+1$, there is a homomorphism $$\varepsilon: {\Z}[{\lambda}]\left[\frac{1}{P(\lambda)}\right]\to{\Z}/q{\Z}$$ where $\e(\lambda)$ is a generator of the group $(\z{q})^{\ast}$ and where $P=a_0+a_1X+...+a_NX^N$, with integral coefficients. For all sufficiently large primes $q$ we may assume that $(q,a_0)=(q,a_N)=1$. Hence for $q$ sufficiently large the projection  ${\Z}\to{\Z}/q{\Z}$ associates to $P$ a non-trivial polynomial $\overline{P}$ in ${\Z}/q{\Z}[X]$ of degree $N$. 
On the other hand it is well known that   $({\Z}/q{\Z})^{\ast}$ is a cyclic group of order $q-1$, when $q$ is a prime. Thus there exists  $\phi(q-1)=\phi(mn_0)$ elements in ${\Z}/q{\Z}$ generating $({\Z}/q{\Z})^{\ast}$, where $\phi$ is the Euler function. Moreover it is easy to prove that $\lim_{n\to +\infty}\phi(n)= +\infty.$
 Hence for a prime $q$ sufficiently large we get: Card(${\cal G}(({\Z}/q{\Z})^{\ast}))=\phi(q-1)>N\geq$ Card$(\overline{P}^{-1}(0))$ which allows us to choose an element $\o{c}$ in ${\Z}/q{\Z}$ generating $({\Z}/q{\Z})^{\ast}$ and such that $\overline{P}(\o{c})\not=\o{0}$. Hence for all but finitely many  primes  $q=mn_0+1$, we may define a homomorphism $\varepsilon:{\Z}[\lambda]\to{\Z}/q{\Z}$  by setting $\varepsilon(\lambda)=\o{c}$ where $\o{c}$ is a generator of $({\Z}/q{\Z})^{\ast}$, which is possible since $\lambda$ is transcendental over ${\Z}$. Since $\overline{P}(\o{c})\not=0$ we can extend $\varepsilon$ to a homomorphism  $\varepsilon:{\Z}[\lambda][1/P(\lambda)]\to{\Z}/q{\Z}$. Since the elements  $x_1,...,x_k$ are algebraically free over ${\Z}[\lambda]$, we extend the above homomorphism to $B= {\Z}[\lambda][1/P(\lambda)][x_1,...,x_k]$ by choosing arbitrary images for  $x_1,...,x_k$. We still denote by $\varepsilon:B\to{\Z}/q{\Z}$ this homomorphism. Let us remark  that it follows from the above construction that $\lambda$ is sent to an element of multiplicative order $q-1=mn_0$.

We next show  that we can extend   $\e$ to ${\A}$. We first prove that there is an extension of $\e$ to $B[\mu_1,...,\mu_l]$ in such a way that the  $\mu_i$ are sent to non-trivial elements. We assume that there is an integer $0\leq i<l$ such that for all $j\in\{0,...,i\}$ there is a finite field ${\F}^j_q$ of characteristic $q$ which is a finite degree extension of ${\Z}/q{\Z}$ and an extension of $\e$ denoted by $\e^j:B^j=B[\mu_1,...,\mu_j]\to{\F}^j_q$ such that $\e^j(\mu_r)\not=0$ for $r=1,...,j$. Since ${\A}$ is integral over $B$,  there is a polynomial $P_{i+1}=a^{i+1}_0+...+a^{i+1}_{n}X^n$ in $B[X]$ where $a^{i+1}_0$ and $a^{i+1}_{n}$ are invertible in ${\A}$ such that $P_{i+1}(\mu_{i+1})=0$. The homomorphism $\e^i$ associates to $P_{i+1}$ a polynomial $\overline{P}_{i+1}$ which can be assumed to be  irreducible in  ${\F}^i_q[X]$, having a   non-trivial root $x_{i+1}$ in some extension of ${\F}^i_q$. If $\overline{P}_{i+1}$ has no root in ${\F}^i_q$ we take the field extension   ${\F}^{i+1}_q={\F}^i_q[X]/(\overline{P}_{i+1})$ and we set $x_{i+1}=\overline{X}$ where $\overline{X}$ denotes the class of $X$ for the projection ${\F}^i_q[X]\to {\F}^i_q[X]/(\overline{P}_{i+1})$. If $\overline{P}_{i+1}$ has a  non-trivial root  $x_{i+1}$ in ${\F}^i_q$ we set ${\F}^{i+1}_q={\F}^i_q$. This proves, by induction,  that we can extend $\e$ to $B[\mu_1,...,\mu_l]$. To extend $\e$ to ${\A}$ it is sufficient to fix images for its other generators. Since ${\A}$ has a finite number of generators we use the same method as above (using the fact that ${\A}$ is integral over $B$). Let  $\e$ be the homomorphism extended to ${\A}$ and ${\F}_q$ be the extended field. Since $\e(\mu_i)=\o{\gamma}_i\not=0$ for $i=1,...,l$ then $\o{\gamma}_i\in{\F}_q^{\ast}\simeq{\Z}/n{\Z}$ with $n=$\ Card$({\F}_q)-1$, which ends the proof of   Lemma \ref{technique1}.
\end{proof}
\begin{proof}[End of proof of Lemma \ref{hyperbolique1}] 
Let $q$ be a prime satisfying the conclusion of Lemma  6.3. We denote by $\varepsilon:{\A}\to{\F}_q$ the homomorphism given by   Lemma  6.3. This homomorphism combined with the holonomy $H_q$ of $\pi_1(H)$ in $SL(2,{\C})$ induces a  homomorphism $\varphi$ such that the following diagram commutes.
$$\xymatrix{
\pi_1(H) \ar[r]^{H_q} \ar[d]_{\varphi} & SL(2,{\C}) \ar[dl]^{\varrho}\\
SL(2,{\F}_q)}$$
where $\varrho$ is the restriction of the homomorphism   $\varrho:SL(2,{\A})\to SL(2,{\F}_q)$ defined by: $$\left( \begin{array}{clcr}a & b\\c & d\end{array}\right)\mapsto\left( \begin{array}{clcr}\varepsilon(a) & \varepsilon(b)\\ \varepsilon(c) & \varepsilon(d)\end{array}
\right)$$
So we get the following identities: 
$$\varphi(m_i)=\overline{Q}_i\left( \begin{array}{clcr}\o{c} & 0\\0 & \o{c}^{-1}\end{array}\right)\overline{Q}_i^{-1}\ ,\ \varphi(l_i)=\overline{Q}_i\left( \begin{array}{clcr}\o{\gamma}_i & 0\\0 & \o{\gamma}_i^{-1}\end{array}\right)\overline{Q}_i^{-1}\ \ \mbox{for}\ i\in\{1,...,l\}$$
$$\varphi(\mu_j) =
\overline{P}_j\left( \begin{array}{clcr}\overline{1} & \overline{1}\\0 & \overline{1}\end{array} \right)\overline{P}_j^{-1}, \ \
\varphi(\nu_j) =  \overline{P}_j\left( \begin{array}{clcr}\overline{1} & \e(\eta_j)\\0 & \overline{1}\end{array}
\right)\overline{P}_j^{-1}\ \
\mbox{for}\  j\in\{1,...,r\}$$
Let $G_n$ be the subgroup of  $SL(2,{\F}_q)$ defined by:
$$G_n=\left\{\overline{a}=\left( \begin{array}{clcr}\overline{a} & {0}\\0 & \overline{a}^{-1}\end{array} \right)\  \mbox{when}\ \   \overline{a}\in{\F}_q^{\ast}\right\}$$
Since ${\F}_q^{\ast}$ is a cyclic group of order $n$, so is $G_n$. To complete the proof of  (i) it is sufficient to set $g_i=\overline{Q}_i$. To prove (ii), it is sufficient to use the fact that ${\F}_q$ is a field of characteristic $q$ and the form of the elements $\varphi(\mu_j), \varphi(\nu_j)$ for $j=1,...,r$. Indeed this proves that $\varphi(\pi_1(U_j))$ is either isomorphic to ${\Z}/q{\Z}$ or to ${\Z}/q{\Z}\times{\Z}/q{\Z}$ depending on whether the elements  $\overline{1}$ and $\varepsilon(\eta_j)$ are linearly dependant or not. This ends the proof of Lemma 6.2.
\end{proof}
\begin{rem} Lemma \ref{hyperbolique1} can be easily extended to the case of a finite number of complete finite volume hyperbolic manifolds. More precisely, if  $H_1,...,H_{\nu}$  denote $\nu$ hyperbolic submanifolds in $N^3$, we can write  Lemma \ref{hyperbolique1} simultaneously for the $\nu$ submanifolds by choosing the same prime $q\in{\cal P}_{n_0}$, the same group $K$, the same cyclic group $G_n\simeq{\Z}/n{\Z}\subset K$ and the same element $c\in G_n$ of multiplicative order $mn_0$. Hence we get the following corollary.
\end{rem}
\begin{cor}\label{hyperbolique2}  Let $H_1,...,H_{\nu}$ be $\nu$  submanifolds of $N^3$ whose interiors admit a complete finite volume hyperbolic  structure.   Then for any integer  $n_0\geq 1$ and for all but finitely many  primes  $q$ of the  form $mn_0+1$, there exists a finite group $K$, a cyclic subgroup $G_n\simeq{\Z}/n{\Z}$ of $K$, an element $\o{c}\in G_n$ of multiplicative order $mn_0$,  elements $\o{\gamma}^i_j$ in $G_n$, $i=1,...,\nu$, $j=1,...,l_i$, and group homomorphisms $\varphi^{H_i}\co\pi_1(H_i)\to K$, $i=1,...,\nu$ satisfying the following properties:

{\rm(i)}\qua for each $i\in\{1,...,\nu\}$ and $j\in\{1,...,l_i\}$ there is an element $g^i_j\in K$  such that $\varphi^{H_i}(\pi_1(T^i_j))\subset g^i_jG_n(g^i_j)^{-1}\simeq{\Z}/n{\Z}$,

{\rm(ii)}\qua for each $i\in\{1,...,\nu\}$ and $j\in\{1,...,l_i\}$ we have the following equalities: $\varphi^{H_i}(m^i_j)=g^i_j\o{c}(g^i_j)^{-1}$ and $\varphi^{H_i}(l^i_j)=g^i_j\o{\gamma}^i_j(g^i_j)^{-1}$,

{\rm(iii)}\qua for each $i\in\{1,...,\nu\}$ and $j\in\{1,...,r_i\}$ the group $\varphi^{H_i}(\pi_1(U^i_j))$ is isomorphic to either $\z{q}$ or  $\z{q}\times\z{q}$.
\end{cor}
\subsubsection{A finite covering lemma for Seifert fibered manifolds}
In this section we construct a class of finite coverings for Seifert fibered manifolds with non-empty boundary homeomorphic to a product  which allows to extend the {\it hyperbolic} coverings given by Corollary \ref{hyperbolique2}. We show here that these coverings may be extended if some simple combinatorial conditions are checked and we will see that these combinatorial conditions can always be satisfied  up to finite covering over $N^3$. Throughout this paragraph we  consider a Seifert piece $S$ of $N^3$ identified to a product $F\times S^1$, where $F$ denotes an orientable  surface of genus $g\geq 1$ with at least two boundary components. We fix two intergers $n>1$ and $c$ in ${\Z}^{\ast}$ and we denote by $\alpha$ the order of $\o{c}$ in $\z{n}$. Then the main result of this section is the following.
\begin{lem}\label{seifertique} Let $S$ be a Seifert fibered space homeomorphic to $F\times S^1$. We denote by $D^1,...,D^l, G_1,...,G_r$ the components of $\b F$ and we set $d^i=[D^i]\in\pi_1(F)$ and $\delta_j=[G_j]\in\pi_1(F)$ (for a choice of base points). Let $\gamma(S)=\{\gamma^1,...,\gamma^l\}$ be a finite sequence of integers. Then there exists a finite covering $\pi\co\t{S}=\t{F}\times{\bf S}^1\to S=F\times{\bf S}^1$  inducing the trivial covering on the boundary and satisfying the following property: there exists a group homomorphism $\varphi\co\pi_1(\t{S})\to{\Z}/n{\Z}\times G$, where $G$ denotes a finite abelian group such that:

{\rm(i)}\qua for each component  $T^i_j=D^i_j\times{\bf S}^1$ ($j=1,..,deg(\pi)$) of $\pi^{-1}(T^i)=\pi^{-1}(D^i)\times{\bf S}^1$, we have $\varphi(\pi_1(T^i_j))\subset {\Z}/n{\Z}\times\{0\}$ and in particular we have the equalities:
 $\varphi(\t{t})=(\overline{c},0)$ and $\varphi(d^i_j)=(\overline{\gamma}^i,0)$, where $d^i_j=[D^i_j]\in\pi_1(\t{F})$ and where $\t{t}$ denotes the fiber of $\t{S}$,
 
{\rm(ii)}\qua for each component  $U_j$ of $\pi^{-1}(G_j\times{\bf S}^1)$ the group  $\ker(\varphi|\pi_1(U_j))$ is the $\alpha\times\alpha$-characteristic subgroup  in $\pi_1(U_j)$.
\end{lem}
\begin{proof} Let $N_0$ be the integer defined by $N_0=\gamma^1+...+\gamma^l$. Then the proof of Lemma 6.5 is splitted is two cases.

\medskip
{\bf Case 1}\qua We first assume that  $N_0\equiv 0\ (n)$. Then we show in this case that $S$ itself satisfies the conclusion of  Lemma \ref{seifertique}. Recall that with the notations of  Lemma 6.5 the group $\pi_1(S)$ has a presentation:
$$\langle a_1,b_1,...,a_g,b_g,d^1,...,d^l,\delta_{1},...,\delta_r,t:$$
$$[t,d^i]=[t,\delta_j]=[t,a_i]=[t,b_i]=1,  \prod_{i=1}^{i=g}[a_i,b_i]\prod_{j=1}^{j=l}d^j\prod_{k=1}^{k=r}\delta_k=1\rangle$$    with $n\geq 2$ and $r\geq 2$ (Indeed recall that $N$ satisfies the conclusion of Proposition 1.4. In particular, $N$ is a finite covering $P\co N\to N'$ of a Haken manifold $N'$ such that for each canonical torus $T$ of $W_{N'}$ and for each geometric piece $S$ adjacent to $T$ in $N'$ the space $(P|S)^{-1}(T)$ is made of at least two connected components). We show here that we may construct a homomorphism $\varpi\co\pi_1(S)\to\z{n}\times(\z{\alpha})^{r-1}$ such that $\varpi(\langle d^i,t\rangle  )\subset\z{n}\times\{0\}$ and satisfying the following equalities:
\begin{itemize}
{\item $\varpi(d^i)\equiv(\o{\gamma}^i,0)$ for every $i=1,...,l$,}
{\item  $\varpi(t)\equiv (\o{c},0)$ and the group $\ker(\varpi|\langle \delta_j,t\rangle  )$ is the $\alpha\times\alpha$-characteristic subgroup of  $\langle \delta_j,t\rangle  $ for $j=1,...,r$.}
\end{itemize}
Then consider the group $K$ defined by the  following relations:
$$K=\left\langle d^1,...,d^l,\delta_{1},...,\delta_r,t:\ [t,d^i]=[t,\delta_j]=1, \left(\prod_{j=1}^{j=l}d^j\right)=1, \left(\prod_{k=1}^{k=r}\delta_k\right)=1\right\rangle$$
obtained from $\pi_1(S)$ by killing the generators $a_i, b_i$ for $i=1,...,g$ and adding two relations. Denote by $\pi\co\pi_1(S)\to K$ the corresponding projection homomorphism. Then we define a homomorphism $\theta\co K\to{\Z}\times{\Z}^{r-1}$ by setting:
\begin{itemize}
{\item $\theta(d^1)=(\gamma^1,0,...,0),...,\theta(d^{l-1})=(\gamma^{l-1},0,...,0)$,} 
{\item  $\theta(\delta_1)=(0,1,0,...,0),...,\theta(\delta_{r-1})=(0,...,0,1)$ and $\theta(t)=(c,0,...0)$.}
\end{itemize}
Since  $\prod_i d^i=1$ we get: $\theta(d^l)=-(\gamma^1+...+\gamma^{l-1})\times\{0\}\equiv(\gamma^l\ \ (n))\times\{0\}$ and since $\prod_j\delta_j=1$ we have: $\theta(\delta_r)=(0,1,...,1)$. Finally, if  $\lambda\co{\Z}\times{\Z}^{r-1}\to\z{n}\times(\z{\alpha})^{r-1}$ is the canonical epimorphism, then the homomorphism $\phi$ defined by the composition $\lambda\circ\theta\circ\pi$ satisfies the conclusion of Lemma \ref{seifertique}. This ends the proof of  Lemma 6.5 in case 1.

\medskip
{\bf Case 2}\qua We now assume that $N_0=\gamma^1+...+\gamma^l\not\equiv0\ \ (n)$. So there exists an integer $p>1$ (that may be chosen minimal) such that: $(\star\star)$\   $pN_0=p\gamma^1+...+p\gamma^l\equiv0\ \ (n)$.  Let $\pi:\t{S}\to S$ be the finite covering of degree $p$ of $S$ corresponding to the following homomorphism:
 $$\pi_1(S)\stackrel{h}{\to}\langle a_1\rangle  \simeq{\Z}\stackrel{-}\to\frac{\langle a_1\rangle  }{\langle a_1^p\rangle  }\simeq\frac{\Z}{p{\Z}}$$
It follows from the above construction that this covering induces the trivial covering on  $\b S$. So each component   $T$ of $\b S$ has $p$ connected components in its pre-image by  $\pi$. With the same notations as in the above paragraph,  the group   $\pi_1(\t{S})$ has a  presentation:
$$\langle a_1,b_1,...,a_{\t{g}},b_{\t{g}},d^1_1,...,d^1_p,...,d^l_1,...,d^l_p,\t{\delta}_{1},...,\t{\delta}_{\t{r}}:$$    $$\left(\prod_{i=1}^{i=\t{g}}[a_i,b_i]\right).\left(\prod_{i,j}d^i_j\right).\left(\prod_{k=1}^{k=\t{r}}\t{\delta}_k\right)=1\rangle\times\l\t{t}\r$$
Then we show that we can construct a homomrphism $\varpi:\pi_1(\t{S})\to\z{n}\times(\z{\alpha})^{\t{r}-1}$ such that $\varpi(\langle d^i_j,\t{t}\rangle)\subset\z{n}\times\{0\}$ and satisfying the following equalities:
\begin{itemize}
{\item $\varpi(d^i_j)\equiv(\o{\gamma}^i,0)$ for every $j=1,...,p$ and $i=1,...,l$,}
{\item  $\varpi(\t{t})\equiv (\o{c},0)$ and the group $\ker(\varpi|\langle \t{\delta}_j,\t{t}\rangle  )$ is the  $\alpha\times\alpha$-caracteristic subgroup of  $\langle \t{\delta}_j,\t{t}\rangle  $ for $j=1,...,\t{r}$.}
\end{itemize}
Consider now the group  $K$ obtained from  $\pi_1(\t{S})$ by setting: 
$$K=\langle d^1_1,...,d^1_p,...,d^l_1,...,d^l_p,\t{\delta}_{1},...,\t{\delta}_{\t{r}},\t{t}:$$
$$[\t{t},d^i_j]=[\t{t},{\delta}_j]=1, \left(\prod_{i,j}d^i_j\right)=1, \left(\prod_{k=1}^{k=\t{r}}\t{\delta}_k\right)=1\rangle  $$
Let $\pi:\pi_1(\t{S})\to K$ be the  corresponding canonical epimorphism. We define a homomorphism $\theta:K\to\z{n}\times(\z{\alpha})^{\t{r}-1}$ by setting:
\begin{itemize}
{\item $\theta(d^1_1)=(\gamma^1,0,...,0),...,\theta(d^1_p)=(\gamma^1,0,...,0),...,$}
{\item $\theta(d^l_1)=(\gamma^l,0,...,0),...,\theta(d^l_{p-1})=(\gamma^{l},0,...,0)$,} {\item $\theta(\delta_1)=(0,1,0,...,0),...,\theta(\delta_{\t{r}-1})=(0,...,0,1)$ and $\theta(t)=(c,0,...0)$.}
\end{itemize}
Since $\prod_{i,j} d^i_j=1$ we get: $\theta(d^l_p)=-(p\gamma^1+...+(p-1)\gamma^{l})\equiv\gamma^l\ \ (n)$ and since $\prod_j\t{\delta}_j=1$ we have: $\theta(\delta_{\t{r}})=(0,1,...,1)$. Finally if we denote by $\lambda:{\Z}\times{\Z}^{\t{r}-1}\to\z{n}\times(\z{\alpha})^{\t{r}-1}$ the canonical projection then the homomorphism $\phi$ defined by the  composition $\lambda\circ\theta\circ\pi$ satisfies the conclusion of  Lemma 6.5. This completes the proof of Lemma \ref{seifertique}. 
\end{proof}
\subsection{Proof of  Proposition \ref{def}}
Throughout this section $N^3$ will denote a closed Haken manifold with non-trivial Gromov simplicial volume, whose Seifert pieces are product. Let  $n_0\geq 1$ be a fixed integer. We denote by   $H_1,...,H_{\nu}$ the hyperbolic components and by $S_1,...,S_t$ the Seifert pieces of $N\setminus W_N$. We want to apply Corollary \ref{hyperbolique2}  to $H_1,...,H_{\nu}$ uniformly with respect to the integer  $n_0$. To do this we first fix system of ``longitude-meridian" on each boundary component of these manifolds. This choice will be determined in the following way: Let $H$ be a hyperbolic manifold and let $T$ be a component of $\b H$. If $T$ is adjacent on both sides to hyperbolic manifolds we fix a system of ``longitude-meridian" arbitrarily. We now assume that $T$ is adjacent to a Seifert piece in $N$ denoted by $S=F\times S^1$. We identify a regular neighborhood of $T$ with $T\times[-1,1]$, where $T=T\times\{0\}$, $T^-=T\times\{-1\}$ and $T^+=T\times\{+1\}$. We assume that $T^+$ is a component of $\b S$ and that $T^-$ is a component of $\b H$ and we denote by $h_T:T^+\to T^-$ the corresponding gluing homeomorphism. Let $t$ be the fiber of $S$ represented in $T^+$ and let $d$ be the homotopy class of the boundary component of the base $F$ of $S$ corresponding to $T^+$. Then the curves $(t,d)$ represent a system of ``longitude-meridian" for $\pi_1(T^+)$ and allow us to associate to $T^-\subset\b H$ a ``longitude-meridian" system by setting:  $$m=h_T(t) \ \mbox{and} \ l=h_T(d)$$    
We now give some notations: for a hyperbolic manifold $H_i$ of $N$, we denote by $T_1^i,...,T^i_{l_i}$ the components of $\b H_i$ adjacent to a Seifert piece and by $U^i_1,...,U^i_{r_i}$ those which are adjacent on both sides to hyperbolic manifolds. For each $T^i_j$, we denote by $(m^i_j,l^i_j)$ its ``longitude-meridian" system corresponding to the construction described above.

We now describe how the hyperbolic pieces of $N$ allow us to associate, via Corollary 6.4, a sequence of integers $\gamma(S)$ to each Seifert piece of $N$, in the sense of Lemma 6.5. Let $S$ be a Seifert piece in $N$, we denote by $H_1,...,H_m$ the hyperbolic pieces adjacent to $S$ along $\b S$ and we fix a torus $T_1$ in $\b S$ adjacent to $H_1$ (say). It follows from Corollary 6.4 that there exists a homomorphism $\varphi^1\co\pi_1(H_1)\to K$ such that $\varphi^1(\pi_1(T_1))\subset gG_ng^{-1}$, where $G_n$ is a $n$-cyclic subgroup of $K$ and such that $\varphi^1(m)=g\o{c}g^{-1}$, $\varphi^1(l)=g\o{\gamma}_1g^{-1}$ where $\o{c}$ and $\o{\gamma}_1$ are elements of $G_n=\z{n}$ and where $(m,l)$ denotes the ``longitude-meridian'' system of $T_1$. Let $c$ and $\gamma_1$ be representatives in ${\Z}$ of $\o{c}$ and $\o{\gamma}_1$. Then we set $\gamma_1(S)=\gamma_1$. Applying the same method for all tori of $\b S$ which are adjacent to hyperbolic components we get a sequence $\{\gamma_1(S),...,\gamma_{n_i}(S)\}=\gamma(S)$ associated to $S$, when $S$ is a Seifert piece in $N$. We fix a suitable prime $q$ of the form $mn_0+1$ (i.e.\ we choose $q$ sufficiently large) and we apply Corollary 6.4 to the hyperbolic manifolds $H_1,...,H_{\nu}$. This means that for each $i\in\{1,...,\nu\}$ there exists a homomorphism $\varphi^{H_i}\co\pi_1(H_i)\to K$ satisfying the conclusion of Corollary 6.4. This allows us to associate to each Seifert piece $S_1,...,S_t$ an integer sequence $\gamma(S_1),...,\gamma(S_t)$. So the proof of Proposition 1.12 will be splitted in two cases.
\subsubsection{Proof of Proposition 1.12: Case 1}
We first assume that we can apply Lemma 6.5 for each Seifert piece $S$ of $N\setminus W_N$ and the associated integer sequence $\gamma(S)$ (i.e.\ without using a finite covering). It follows from Lemma 6.5, that for each $i\in\{1,...,t\}$, there exists a group homomorphism $\varphi_{S_i}\co\pi_1(S_i)\to\z{n}\times G_i$ satisfying properties (i) and (ii) of this lemma for the sequence $\gamma(S_i)$ with $\alpha=q-1=mn_0$.

It follows from \cite[Lemma 4.1]{He2} or \cite[Theorems 2.4 and 3.2]{Lu} that for each $i\in\{1,...,\nu\}$ (resp.\ $i\in\{1,...,t\}$) there exists a finite group $H$ (resp.\ $L_i$) and a group homomorphism $\hat{\varphi}_{H_i}\co\pi_1(H_i)\to H$ (resp.\ $\hat{\varphi}_{S_i}\co\pi_1(S_i)\to L_i$) which induces the $q\times q$-characteristic covering on $\b H_i$ (resp.\ $\b S_i$). For each $i\in\{1,...,\nu\}$ (resp.\ $i\in\{1,...,t\}$) we consider the homomorphism $\psi_{H_i}$ (resp.\ $\psi_{S_i}$) defined by the following formula:
$$\psi_{H_i}=\varphi^{H_i}\times\hat{\varphi}_{H_i}\co\pi_1(H_i)\to K\times H$$  $$\psi_{S_i}=\varphi^{S_i}\times\hat{\varphi}_{S_i}\co\pi_1(S_i)\to(\z{n}\times G_i)\times L_i$$
where $\varphi^{H_i}$ is given by Corollary 6.4.  The above homomorphisms allow us to associate to each Seifert piece $S$ of $N\setminus W_N$ a finite covering $p_S\co\t{S}\to S$. Define the set ${\cal R}$ by ${\cal R}:=\{p_S\co\t{S}\to S$ when $S$ describes the Seifert pieces of $N\}\cup\{p_H\co\t{H}\to H$ when $H$ describes the hyperbolic pieces of $N\}$. Since for each Seifert piece $S$ of $N$ the homomorphism $\psi_S$ sends the homotopy class of the regular fiber of $S$, denoted by $t_S$, to an element of order $qmn_0$, then to prove Proposition 1.12 it is sufficient to show the following result.
\begin{lem}\label{ouf}
There exists a finite covering $r\co\t{N}\to N$ such that for each component $S$ of $N\setminus W_N$ and for each component $\t{S}$ of $r^{-1}(S)$, the induced covering $r|\t{S}\co\t{S}\to S$ is equivalent to the covering corresponding to $S$ in the set ${\cal R}$. 
\end{lem}  
In the proof of this result, it will be convenient to call a co-dimension 0 submanifold $X_k$ of $N$ a $k-chain$ of $N$ if $X_k$ is a connected manifold made of exactly $k$ components of $N\setminus W_N$. Then we prove Lemma \ref{ouf} by induction on $k$.
\begin{proof}[Proof of Lemma \ref{ouf}]
When $k=1$ this is an obvious consequence of Lemma 6.2, if the 1-chain $X_1$ is hyperbolic or of Lemma 6.5, if $X_1$ is a Seifert piece. Indeed it is sufficient to take the associated homomorphism of type $\psi_H$ or $\psi_S$. We fix now an integer $k\leq t+\nu$ and we set the following inductive hypothesis:

{\sl $({\cal H}_{k-1})$: for each $j<k\leq t+\nu$ and for each $j$-chain $X_j$ of $N$, there exists a finite covering $p_j\co\t{X}_j\to X_j$ such that for each component $S$ of $X_j\setminus W_N$ and for each component $\t{S}$ of $p^{-1}_j(S)$ the induced covering $p_j|\t{S}\co\t{S}\to S$ is equivalent to the covering $p_S$ corresponding to $S$ in the set ${\cal R}$.}

Let $X_k$ be a $k$-chain in $N$. We choose a $(k-1)$-chain denoted by $X_{k-1}$ in $X_k$ and we set $X_1$ the (connected) component of $X_k\setminus X_{k-1}$.

\medskip
{\bf Case 1.1}\qua We first assume that $X_1$ is a Seifert piece of $N$, denoted by $S$. Let $H_1,...,H_m$ be the hyperbolic pieces of $X_{k-1}$ adjacent to $S$ and let $S_1,...,S_k$ be the Seifert pieces of $X_{k-1}$ adjacent to $S$. The hyperbolic manifold $H_i$ is adjacent to $S$ along tori $(T^{i,-}_1,...,T^{i,-}_{\nu_i})\subset\b H_i$ and  $(T^{i,+}_1,...,T^{i,+}_{\nu_i})\subset\b S$ and $S_j$ is adjacent to $S$ along tori  $(U^{j,-}_1,...,U^{j,-}_{n_j})\subset\b S_j$ and  $(U^{j,+}_1,...,T^{j,+}_{n_j})\subset\b S$. With these notations the fundamental group of $S$ has a presentation:
$$\l a_1,b_1,...,a_g,b_g,d^1_1,...,d^1_{\nu_1},...,d^s_1,...,d^s_{\nu_s},\delta^1_1,...,\delta^1_{r_1},...,\delta_1^{\beta},...,\delta^{\beta}_{r_{\beta}} \ :$$
$$\left(\prod_i[a_i,b_i]\right).\left(\prod_{i,j}d^i_j\right).\left(\prod_{i,j}\delta^i_j\right)=1\r\times\l t\r$$  
Where the group $\l t,\delta^i_j\r$ corresponds to $\pi_1(U^{i,+}_j)$ and $\l t,d^i_j\r$ corresponds to $\pi_1(T^{i,+}_j)$. We denote by $h^i_k\co T^{i,+}_k\to T^{i,-}_k$ and by $\varphi^j_k\co U^{j,+}_k\to U^{j,-}_k$ the gluing homeomorphism in $N$ (see figure 5). Let $p_{X_{k-1}}\co\t{X}_{k-1}\to X_{k-1}$ be the covering given by the inductive hypothesis. In particular, for each hyperbolic piece $H_i$ (resp.\ Seifert piece $S_j$) of $X_{k-1}$ and for each component $\t{H}_i$ of $p^{-1}_{X_{k-1}}(H_i)$ (resp.\ $\t{S}_j$ of  $p^{-1}_{X_{k-1}}(S_j)$) the covering $p_{X_{k-1}}|\t{H}_i$ (resp.\ $p_{X_{k-1}}|\t{S}_j$) is equivalent to $p_{H_i}$ (resp.\ $p_{S_j}$) in ${\cal R}$. Denote by $\psi_{H_i}$ (resp.\ $\psi_{S_j}$) the homomorphisms corresponding to $p_{H_i}$ (resp.\ to $p_{S_j}$):
$$\psi_{S_j}=\varphi_{S_j}\times\hat{\varphi}_{S_j}\co\pi_1(S_j)\to(\z{n}\times G_i)\times L_i$$
 $$\mbox{and}\ \ \psi_{H_i}=\varphi^{H_i}\times\hat{\varphi}_{H_i}\co\pi_1(H_i)\to K\times H$$
where $K, H, G_i$ and $L_i$ are finite groups. In particular, we have the following properties $({\cal P}^{i,-}_j)$:
\begin{itemize}
\item[\rm(a)] $\psi_{H_i}|\pi_1(T^{i,-}_j)=\varphi_{H_i}|\pi_1(T^{i,-}_j)\times\hat{\varphi}_{H_i}|\pi_1(T^{i,-}_j)$ is a homomorphism from $\pi_1(T^{i,-}_j)$ to $(g^i_j.G_n.(g^i_j)^{-1})\times(\z{q}\times\z{q})\subset K\times H$ with $g^i_j\in K$ and $\varphi_{H_i}(m^i_j)=(g^i_j.\o{c}.(g^i_j)^{-1},0,0)$ (where $\o{c}$ is an element of order $\alpha=mn_0$ in $G_n$) and $\varphi_{H_i}(l^i_j)=(g^i_j\o{\gamma}^i_j(g^i_j)^{-1},0,0)$ and $\hat{\varphi}_{H_i}(\pi_1(T^{i,-}_j))=\{0\}\times\z{q}\times\z{q}$ for $i=1,...,m$ and $j=1,...,\nu_i$.
\item[\rm(b)] the groups $\ker(\psi_{S_i}|\pi_1(U^{i,-}_j))$ are $q\alpha\times q\alpha$-characteristic in $\pi_1(U^{i,-}_j)$.
\end{itemize} 
We consider the integer sequence $\gamma_{{\cal H}}(S)=\{\gamma^i_j\}_{i,j}$ of liftings in ${\Z}$ of $\{\o{\gamma}^i_j\}_{i,j}$. It follows from the hypothesis of Case 1, that we can apply Lemma 6.5 to $S=F\times S^1$; this means that there exists a homomorphism $\psi_S\co\pi_1(S)\to\z{n}\times G\times L_S$ satisfying the following equalities denoted by $({\cal P}_S)$:
\begin{itemize}
\item[\rm(c)] the group $\ker(\psi_S|\l\delta^i_j, t\r)=\ker(\psi_S|\pi_1(U^{i,+}_j))$ is the chatacteristic subgroup of index $q\alpha\times q\alpha$ in $\l\delta^{i,+}_j,t\r$ for $i=1,...,t$ and $j=1,...,n_i$.
\item[\rm(d)] $\psi_{S}|\pi_1(T^{i,+}_j)=\varphi_{S}|\pi_1(T^{i,+}_j)\times\hat{\varphi}_{S}|\pi_1(T^{i,+}_j)\co\pi_1(T^{i,+}_j)\to\z{n}\times\{0\}\times(\z{q}\times\z{q})\subset\z{n}\times\{0\}\times L_i$ with $\varphi_{S}(d^i_j)=(\o{\gamma}^i_j,0,0)$, $\varphi_{S}(t)=(\o{c},0,0)$ and $\hat{\varphi}_{S}(\pi_1(U^{i,+}_j))=\{0\}\times\z{q}\times\z{q}$ for $i=1,...,t$ and $j=1,...,n_i$.
\end{itemize}

\begin{figure}[ht!]
\centerline{
\relabelbox\small
\epsfbox{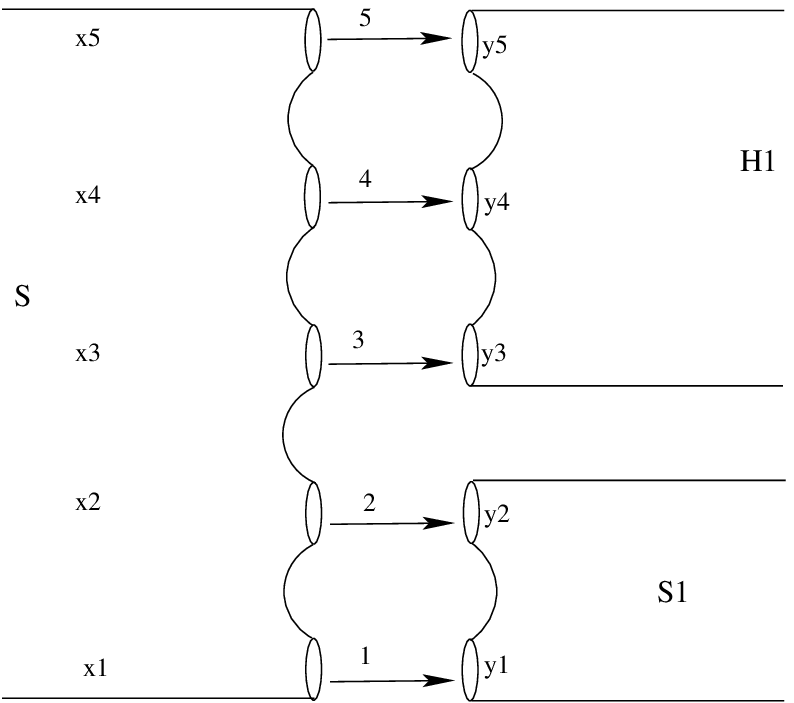}
\relabel {1}{$\varphi^1_2$}
\relabel {2}{$\varphi^1_1$}
\relabel {3}{$h^1_3$}
\relabel {4}{$h^1_2$}
\relabel {5}{$h^1_1$}
\relabel {x1}{$U^{1,+}_2=(\delta^1_2,t)$}
\relabel {x2}{$U^{1,+}_1=(\delta^1_1,t)$}
\relabel {x3}{$T^{1,+}_3=(d^1_3,t)$}
\relabel {x4}{$T^{1,+}_2=(d^1_2,t)$}
\relabel {x5}{$T^{1,+}_1=(d^1_1,t)$}
\relabel {y1}{$U^{1,-}_2$}
\relabel {y2}{$U^{1,-}_1$}
\relabel {y3}{$T_3^{1,-}=(l^1_3,m^1_3)$}
\relabel {y4}{$T_2^{1,-}=(l^1_2,m^1_2)$}
\relabel {y5}{$T_1^{1,-}=(l^1_1,m^1_1)$}
\relabel {S}{$S$}
\relabel {S1}{$S_1$}
\relabel {H1}{$H_1$}
\endrelabelbox}\caption{}
\end{figure}

Denote by $p_S\co\t{S}\to S$ the finite covering corresponding to $\psi_S$, by $\eta_S$ the degree of $p_S$ and by $\eta_{X_{k-1}}$ the degree of $p_{X_{k-1}}$. Let $T^{i,+}_j$ (resp.\ $T^{i,-}_j$) be a torus in $\b S$ (resp.\ in $\b H_i$). It follows from the construction of the coverings $p_{X_{k-1}}$ and $p_S$ that $p_{X_{k-1}}$ (resp.\ $p_S$) is a covering of degree $\eta^{i,-}_j=|\psi_{H_i}(\pi_1(T^{i,-}_j))|=|\psi_{H_i}(\l l^i_j,m^i_j\r)|$ (resp.\ $\eta^{i,+}_j=|\psi_{S}(\pi_1(T^{i,+}_j))|=|\psi_{S}(\l d^i_j,t\r)|$). If $\alpha^{i,+}_j$ (resp.\ $\alpha^{i,-}_j$) denotes the number of connected components of $p_S^{-1}(T^{i,+}_j)$ (resp.\ of $p_{X_{k-1}}^{-1}(T^{i,-}_j)$) we can write:
$$\eta^{i,+}_j\times\alpha^{i,+}_j=\eta_S,\  \eta^{i,-}_j\times\alpha^{i,-}_j=\eta_{X_{k-1}}\  \mbox{and}\ \  \eta^{i,+}_j=\eta^{i,-}_j\eqno{(1)}$$ 
by properties ${\cal P}^{i,-}_j$ and ${\cal P}_S$. For each component $U^{i,-}_j$ of $S_i$ (resp.\ $U^{i,+}_j$ of $S$), the covering $p_{X_{k-1}}$ (resp.\ $p_S$) induces the $q\alpha\times q\alpha$-characteristic covering. If $\beta^{i,+}_j$ (resp.\ $\beta^{i,-}_j$) denotes the number of connected components of $p_S^{-1}(U^{i,+}_j)$ (resp.\ of $p_{X_{k-1}}^{-1}(U^{i,-}_j)$), we can write:
$$q^2\alpha^2\times\beta^{i,+}_j=\eta_S, \ \ \ \ q^2\alpha^2\times\beta^{i,-}_j=\eta_{X_{k-1}}\eqno{(2)}$$
We want to show that there are two positive integers $x$ and $y$ independant of $i$ and $j$ satisfying the following equalities:
$$x\alpha^{i,+}_j=y\alpha^{i,-}_j \ \ \ \ \ x\beta^{i,+}_j=y\beta^{i,-}_j\eqno{(3)}$$   
Using (1), it is sufficient to choose $x$ and $y$ in such a way that $x\eta_S=y\eta_{X_{k-1}}$ which is possible. So we take $x$ copies $\t{S}^1,...,\t{S}^x$ of $\t{S}$ and $y$ copies $\t{X}^1_{k-1},...,$ $\t{X}^y_{k-1}$ of $\t{X}_{k-1}$ with the coverings $p^i_S\co\t{S}^i\to S$ (resp.\ $p^i_{X_{k-1}}\co\t{X}_{k-1}^i\to X_{k-1}$) equivalent to $p_S$ (resp.\ $p_{X_{k-1}}$). Then consider the space ${\cal X}$ defined by
$${\cal X}=\left(\coprod_{i\leq i\leq x}\t{S}^i\right)\coprod\left(\coprod_{1\leq j\leq y}\t{X}_{k-1}^j\right)$$
Note that it follows from the above arguments that the spaces $\coprod_{i\leq i\leq x}\t{S}^i$ and $\coprod_{1\leq j\leq y}\t{X}_{k-1}^j$ have the same number of boundary components. Thus it is sufficient to show that we can glue together the connected components of ${\cal X}$ via those of $(p^i_S)^{-1}(\b S)$ and of $(p^i_{X_{k-1}})^{-1}(\b X_{k-1})$ (see figure 5). To do this,  we fix a component $\t{T}^{i,+}_j$ (resp.\ $\t{T}^{i,-}_j$) of $p_S^{-1}(\t{T}^{i,+}_j)$ (resp.\ $p_{X_{k-1}}^{-1}(\t{T}^{i,-}_j)$) and we proceed as before with the components of $p_S^{-1}(\t{U}^{i,+}_j)$ (resp.\ $p_{X_{k-1}}^{-1}(\t{U}^{i,-}_j)$). Then it is sufficient to prove that there exist homeomorphisms $\t{h}^i_j$ and $\t{\varphi}^i_j$ such that the following diagrams are consistent: 
$$\xymatrix{
\t{T}^{i,+}_j \ar[d]_{{\textstyle(4)} \ \ \ \ p_S|\t{T}^{i,+}_j} \ar[r]^{\t{h}^i_j} & \t{T}^{i,-}_j  \ar[d]^{p_{X_{k-1}}|\t{T}^{i,-}_j}  & & & \t{U}^{i,+}_j \ar[d]_{p_S|\t{U}^{i,+}_j} \ar[r]^{\t{\varphi}^i_j} & \t{U}^{i,-}_j \ar[d]^{p_{X_{k-1}}|\t{U}^{i,-}_j \ \ \ \ {\textstyle(5)}}         \\
{T}^{i,+}_j \ar[r]^{h^i_j} & {T}^{i,-}_j  & & & {U}^{i,+}_j \ar[r]^{\varphi^i_j} & {U}^{i,-}_j}$$
Since the coverings $p_S|\t{U}^{i,+}_j$ and $p_{X_{k-1}}|\t{U}^{i,-}_j$ correspond to the characteristic subgroup of index $q\alpha\times q\alpha$ in $\pi_1(U^{i,+}_j)$ and $\pi_1(U^{i,-}_j)$, it is straightforward that there exists a homeomorphism $\t{\varphi}^i_j$ such that the diagram (5) is consistent (since for each integer $n$, the $n\times n$-characteristic subgroup of $\pi_1(U^{i,-}_j)$ is unique). We now fix a base point $x^+$ (resp.\ $x^-=h^i_j(x^+)$) in $T^{i,+}_j$ (resp.\ $T^{i,-}_j$). So we have $\pi_1(T^{i,+}_j,x^+)=\l d^i_j,t\r$ and $\pi_1(T^{i,-}_j,x^-)=\l l^i_j,m^i_j\r$. By (d), we know that the covering $p_S|\t{T}^{i,+}_j$ corresponds to the homomorphism $\e$ defined by:
$$\e\!=\!\e_i\times\o{\e}_i\!=\!\psi_{S_i}|\l d^i_j,t\r\!=\!\varphi_{S_i}|\l d^i_j,t\r\times\hat{\varphi}_{S_i}|\l d^i_j,t\r\to(\z{n}\times\{0\})\times(\z{q}\times\z{q})$$
$$\mbox{with equalities:} \ \ \psi_{S_i}(d^i_j)=(\varphi_{S_i}(d^i_j),\hat{\varphi}_{S_i}(d^i_j))=((\o{\gamma}^i_j,0),\hat{\varphi}_{S_i}(d^i_j))\eqno{(6)}$$
$$\mbox{and} \ \ \psi_{S_i}(t)=(\varphi_{S_i}(t),\hat{\varphi}_{S_i}(t))=((\o{c},0),\hat{\varphi}_{S_i}(t))$$
It follows from (a) that the covering $p_{X_{k-1}}|\t{T}^{i,-}_j$ corresponds to the homomorphism $\e'\co\l l^i_j,m^i_j\r\to(g^i_jG_n(g^i_j)^{-1})\times(\z{q}\times\z{q})$ defined by:
$$\e'=\e_i'\times\o{\e}_i'=\psi_{H_i}|\l l^i_j,m^i_j\r=\varphi_{H_i}|\l l^i_j,m^i_j\r\times\hat{\varphi}_{H_i}|\l l^i_j,m^i_j\r$$
$$\mbox{with equalities:} \ \ \psi_{H_i}(l^i_j)=(\varphi_{H_i}(l^i_j),\hat{\varphi}_{H_i}(l^i_j))=((g^i_j\o{\gamma}^i_j(g^i_j)^{-1},0),\hat{\varphi}_{H_i}(l^i_j))\eqno{(7)}$$
$$\mbox{and} \ \ \psi_{H_i}(m^i_j)=(\varphi_{H_i}(m^i_j),\hat{\varphi}_{H_i}(m^i_j))=((g^i_j\o{c}(g^i_j)^{-1},0),\hat{\varphi}_{H_i}(m^i_j))$$  where $G_n\simeq\z{n}$.
To prove that the homomorphism $h^i_j$ lifts in the diagram (4) it is sufficient to see that: $(h^i_j)_{\ast}(\ker(\e))=\ker(\e')$. It follows from the above arguments that $\ker(\e)=\ker(\e_i)\cap\ker(\o{\e}_i)$ and $\ker(\e')=\ker(\e_i')\cap\ker(\o{\e}_i')$. We first prove the following equality $(h^i_j)_{\ast}(\ker(\e_i))=\ker(\e_i')$. Using (6) and (7) we know that: 
$$\e_i\co\l d^i_j,t\r\to\z{n} \ \mbox{with} \ \e_i(d^i_j)=\o{\gamma}^i_j \ \mbox{and}\ \e(t)=\o{c} $$
$$\e_i'\co\l l^i_j,m^i_j\r\to g^i_jG_n(g^i_j)^{-1}\simeq\z{n}$$ 
$$\mbox{with}  \ \e_i(l^i_j)=g^i_j\o{\gamma}^i_j(g^i_j)^{-1} \ \mbox{and}\ \e(m^i_j)=g^i_j\o{c}(g^i_j)^{-1} $$
Moreover, since the elements $m^i_j$ and $l^i_j$ have been chosen such that $m^i_j=h^i_j(t)$ and $l^i_j=h^i_j(d^i_j)$, the above arguments imply that $(h^i_j)_{\ast}(\ker(\e_i))=\ker(\e_i')$. Hence it is sufficient to check that $(h^i_j)_{\ast}(\ker(\o{\e}_i))=\ker(\o{\e}_i')$. Since $\ker(\o{\e}_i)$ (resp.\ $\ker(\o{\e}_i')$) is the $q\times q$-characteristic subgroup of $\pi_1(T^{i,+}_j)$ (resp.\ of $\pi_1(T^{i,-}_j)$) this latter equality is obvious. So the lifting criterion implies that there is a homeomorphism $\t{h}^i_j$ such that diagram (4) commutes. Finally the space $\t{N}$ obtained by gluing together the connected components of ${\cal X}$ via the homeomorphisms $\t{\varphi}^i_j$ and $\t{h}^i_j$ satisfies the induction hypothesis $({\cal H}_k)$. This proves Lemma \ref{ouf} in Case 1.1.

\medskip
{\bf Case 1.2}\qua To complete the proof of Lemma \ref{ouf} it remains to assume that the space $X_1$ is a hyperbolic submanifold of $N^3$. In this case the arguments are similar to those of Case 1.1. This ends the proof of Lemma \ref{ouf}.  
\end{proof}
\subsubsection{Proof of Proposition 1.12: Case 2}
We now assume that for some Seifert pieces $\{S_i, i\in I\}$ in $N$, in order to apply Lemma 6.5 we have to take a finite covering of order $p\geq 1$ inducing the trivial covering on the boundary. More precisely, for each $S_i$, $i\in\{1,...,t\}$, we denote by $\gamma(S_i)$ the integer sequence which comes from the hyperbolic coverings via Corollary 6.4 and we denote by $\pi_i\co\t{S}_i\to S_i$ the covering (trivial on the boundary) obtained by applying Lemma 6.5 to $S_i$ with $\gamma(S_i)$. Then we construct a finite covering $\pi\co\t{N}\to N$ such that each component of $\pi^{-1}(S_i)$  is equivalent to the covering $\pi_i\co\t{S}_i\to S_i$ in the following way: for each $i\in\{1,...,t\}$ we denote by $\eta_i$ the degree of $\pi_i$. We define the integer $m_0$ by:
$$m_0=\mbox{l.c.m}(\eta_1,...,\eta_t)$$
 For each $i\in\{1,...,t\}$, we take $t_i=m_0/\eta_i$ copies of $S^i$ denoted by $S^i_1,...,S^i_{t_i}$ and $m_0$ copies of $H_j$ denoted by $H^j_1,...,H^j_{m_0}$ for $j\in\{1,...,m\}$. Since the map $\pi_i$ induces the trivial covering on $\b\t{S}_i$ we may glue together the connected components of the space:
 $$X=\left(\coprod_{1\leq i\leq t}\coprod_{1\leq j\leq t_i} S^i_j\right)\coprod\left(\coprod_{1\leq i\leq t}\coprod_{1\leq j\leq m_0} H^i_j\right)$$
 via the (trivial) liftings of the gluing homeomorphism of the pieces $N\setminus W_N$. This allows us to obtain a Haken manifold $N_1$ which is a finite covering of $N$ and which satisfies the hypothesis of Case 1 (see subsection 6.3.1). It is now sufficient to apply the arguments of subsection 6.3.1 for the induced map $f_1\co M_1\to N_1$. This completes the proof of Proposition 1.12. By paragraph 6.1.1 and paragraph 6.1.2 this completes the proof of Theorem 1.1.  

\medskip
 
{\bf Acknowledgment}\qua The author would like to express his
gratitude to his thesis advisor, professor Bernard Perron, for his
continuing help and encouragement during the process of writing this
article. The author is very grateful to professor Michel Boileau for
numerous and stimulating conversations and to professor Joan Porti for
helpful conversations on deformation theory of geometric structures. I
thank the referee for many useful comments and suggestions.

\Addresses\recd 

\begin{thebibliography}

\bibitem{Be-P} R. Benedetti, C. Petronio, {\it Lectures on Hyperbolic Geometry}, Springer-Verlag, 1992.

\bibitem{Bo-W} M. Boileau, S. Wang, {\it Non-zero degree maps and surface bundles over ${\bf S}^1$}, Differential Geom. 43 (1996), pp. 789-806.

\bibitem{Cu-S} M. Culler, P. Shalen, {\it Varieties of group representations and splittings of 3-manifolds}, Ann. of Math. 117 (1983), pp. 109-146.

\bibitem{De} P. Derbez, {\it Un crit\`ere d'hom\'eomorphie entre vari\'et\'es Haken}, Ph.D. Thesis, Univerit\'e de Bourgogne, 2001.


\bibitem{Do} A. Dold, {\it Lectures on algebraic topology}, Springer-Verlag, 1980.



\bibitem{Ei} R. Eisenbud, {\it Commutative algebra with a view toward algebraic geometry}, Springer-Verlag, 1995.

 
\bibitem{Gro} M. Gromov, {\it Volume and bounded cohomology}, Publi. Math. I.H.E.S. 56 (1982), pp. 5-99.


\bibitem{He1} J. Hempel,  {\it 3-manifolds}, Ann. of Math. Studies 86 Princeton Univ. Press (1976).

\bibitem{He2} J. Hempel,  {\it Residual finiteness for 3-manifolds}, Combinatorial group theory, Ann. of Math. Studies 111 (1987).

\bibitem{I-R} K. Ireland, M. Rosen, {\it A classical introduction to modern number theory}, Springer-Verlag, 1990.
 
\bibitem{Ja} W. Jaco, {\it Lecture on three manifold topology}, Conference board of the Math. Sciences. A.M.S. 43 (1977). 

\bibitem{Ja-S} W. Jaco, P. Shalen, {\it Seifert fibered space in 3-manifolds},
memoirs of the A.M.S. Vol. 21 $n^0$ 220 (1979).

\bibitem{Jo} K. Johannson, {\it Homotopy equivalences of 3-manifolds with boundaries}, Lecture Notes in Math. Vol. 761, Springer-Verlag (1979).



\bibitem{Lu} J. Luecke, {\it Finite cover of 3-manifolds}, Trans. A.M.S. Vol. 310, $n^0$ 1 (1988).

\bibitem{L-W} J. Luecke, Y.Q. Wu, {\it Relative Euler number and finite covers of graph manifolds}, Proceedings of the Georgia Internatinal Topology Conference (1993).

\bibitem{Mu} D. Mumford, {\it Algebraic Geometry I, Complex Projective Varieties}, Springer-Verlag, 1976.





\bibitem{Pe-S} B. Perron, P. Shalen, {\it Homeomorphic graph manifolds: a contribution to the $\mu$ constant problem}, Topology and its Applications, $n^0$ 99 (1), 1999, pp. 1-39.

\bibitem{Pr-S} V.V. Prasolov, A.B. Sossinsky, {\it Knots, Links, Braids and 3-manifolds: An Introduction to the New Invariants in Low-Dimensional Topology}, Trans. of Math. Monographs, Vol 154, 1997. 



\bibitem{Ro} Y. Rong, {\it Degree one maps between geometric 3-manifolds}, Trans A.M.S. Vol. 332 $n^0$ 1 (1992).
 



\bibitem{So} T. Soma, {\it A rigidity theorem for Haken manifolds}, Math. Proc. Camb. Phil. Soc. (1995),  118, pp. 141-160.


\bibitem{Sp} E. Spanier, {\it Algebraic topology}, McGraw-Hill, 1966.

\bibitem{Th1} W.P. Thurston, {\it Hyperbolic structures on 3-manifolds}, Ann. of Math., 124 (1986), pp. 203-246.

\bibitem{Th2} W.P. Thurston, {\it The geometry ant topology of three-manifolds}, Princeton University Mathematics Department (1979).



\bibitem{Wa}  F.\ Waldhausen, {\it On irreducible 3-manifolds which are sufficiently large}, Ann. of Math., Vol 87, 1968, pp. 56-88.
 




\end{thebibliography}
\end{document}